\documentclass{article}

\usepackage{paper_macros}

\begin{document}
\title{Provable Non-Convex Euclidean Distance Matrix Completion: Geometry, Reconstruction, and Robustness}
\author[1]{Chandler Smith\thanks{Corresponding Author, Chandler.Smith@Tufts.edu}}
\author[2]{HanQin Cai}
\author[1]{Abiy Tasissa}
\affil[1]{Department of Mathematics, Tufts University, Medford, MA 02155, USA.}
\affil[2]{Department of Statistics and Data Science and Department of Computer Science, University of Central Florida, Orlando, FL 32816, USA.}
\date{}

\maketitle

\begin{abstract}
The problem of recovering the configuration of points from their partial pairwise distances,  referred to as the Euclidean Distance Matrix Completion (EDMC) problem, arises in a broad range of applications, including sensor network localization, molecular conformation, and manifold learning. In this paper, we propose a Riemannian optimization framework for solving the EDMC problem by formulating it as a low-rank matrix completion task over the space of positive semi-definite Gram matrices. The available distance measurements are encoded as expansion coefficients in a non-orthogonal basis, and optimization over the Gram matrix implicitly enforces geometric consistency through nonnegativity and the triangle inequality, a structure inherited from classical multidimensional scaling. Under a Bernoulli sampling model for observed distances, we prove that Riemannian gradient descent on the manifold of rank-$r$ matrices locally converges linearly with high probability when the sampling probability satisfies $p\geq \cO(\nu^2 r^2\log(n)/n)$, where $\nu$ is an EDMC-specific incoherence parameter. Furthermore, we provide an initialization candidate using a one-step hard thresholding procedure that yields convergence, provided the sampling probability satisfies $p \geq \cO(\nu r^{3/2}\log^{3/4}(n)/n^{1/4})$. A key technical contribution of this work is the analysis of a symmetric linear operator arising from a dual basis expansion in the non-orthogonal basis, which requires analysis of a second order degenerate U-statistic to establish an optimal restricted isometry property in the presence of coupled terms.  Empirical evaluations on synthetic data demonstrate that our algorithm achieves competitive performance relative to state-of-the-art methods. Moreover, we provide a geometric interpretation of matrix incoherence tailored to the EDMC setting and provide robustness guarantees for our method.

\end{abstract}


\section{Introduction} \label{sec: Introduction}

The rapid advancement of technology across various scientific fields has greatly simplified data collection. In many practical applications, however, there are limitations to measurements that can lead to incomplete data.  This can be caused by geographic, climatic, or other factors that determine whether a measurement between two points can be obtained, and as such, some data may be missing \cite{aldibaja2016improving,marti2013multi}. For instance, in protein structure prediction, nuclear magnetic resonance (NMR) spectroscopy experiments yield spectra for protons that are close together, resulting in incomplete known distance information \cite{clore1993exploring}. Similarly, in sensor networks, we may have mobile nodes with known distances only from fixed anchors \cite{boukerche2007localization,kuriakose2014review}. In these and other scenarios, the fundamental problem is determining the configuration of points based on partial information about inter-point distances. This problem is known as the Euclidean Distance Geometry problem, which has numerous applications throughout the applied sciences \cite{biswas2006semidefinite,ding2010sensor,rojas2012distance,porta2018distance,tenenbaum2000global,glunt1993molecular,trosset1997applications,fang2013using,liberti2008branch,einav2023quantitatively}.  

To formulate this problem mathematically, some notation is in order. Let $\{\p_i\}_{i=1}^n\subset\mathbb{R}^r$ denote a set of $n$ points in $\real^{r}$. We define the $n\times r$ matrix $\Pb =[\p_1,\p_2,\cdots,\p_n]^\top$, which has the points as rows. There are two essential mathematical objects related to $\P$. The first object is the Gram matrix $\X \in \mathbb{R}^{n \times n}$, defined as $\X = \Pb \Pb^\top$. By construction, $\X$ is symmetric and positive semi-definite.
The second object is the squared distance matrix $\D \in \mathbb{R}^{n \times n}$, defined entry-wise as $D_{ij} = \Vert \p_i - \p_j \Vert^2_2$. The reason for working with the squared distance matrix instead of the distance matrix will become clear later. Computing $\D$ given $\P$ is conceptually straightforward. However, the inverse problem of determining $\Pb$ from $\D$ is less obvious. To address this problem, we need to precisely define what it means to identify $\P$. Since rigid motions and translations preserve distances, there is no unique $\P$ corresponding to a given squared distance matrix $\D$. From here on, we assume the points are centered at the origin, i.e., for $\one$ as a column vector of ones, $\Pb^\top \bm{1} = \bm{0}$. This implies that $\X \bm{1} = \Pb \Pb^\top \bm{1} = \bm{0}$. We refer to $\Pb$ and $\X$ with this relationship as the centered point and centered Gram matrix, respectively. Since the Gram matrix is invariant under rigid motions, these assumptions allow for a one-to-one correspondence between $\D$ and $\X$.

When we have access to all the distances, a central result in \cite{torgerson1952multidimensional} provides the following one-to-one correspondence between $\D$ and a centered $\X$:
\begin{align}
    \X &= -\frac{1}{2}\J\D\J \label{eqn: D to X},  \\ \
    \D &= \diag(\X)\one^\top +\one\diag(\X)^\top - 2\X, \label{eqn: X to D}
\end{align}
where $\diag(\cdot)$ inputs an $n\times n$ matrix and returns a column vector with the entries along the diagonal, and $\J = \I - \frac{1}{n}\oot$. Once $\X$ is reconstructed using the above formula, $\Pb$ can be computed from the $r$-truncated eigendecomposition of $\X$. It is important to note that, as previously mentioned, $\P$ is unique up to rigid motions. This procedure for computing $\Pb$ from a full squared distance matrix $\D$ is known as classical multidimensional scaling (Classical MDS) \cite{young1938discussion,torgerson1952multidimensional,torgerson1958theory,gower1966some}, and for the truncated eigendecomposition $\X = \U\bm{\Lambda}\U^\top$ with $\U\in\real^{n\times r}$ and $\bm{\Lambda}\in\real^{r\times r}$, 
\begin{equation}\label{eqn: MDS Equation}
    \P = \U\bm{\Lambda}^{1/2}.
\end{equation}
We note that $\X\one = \bm{0}$ also implies that $\U^\top\one=\bm{0}$. In many practical scenarios, the distance matrix may be incomplete, making classical MDS inapplicable for determining the point configuration. However, notice that $\rank(\X) \leq r$, and one can show that $\rank(\D)\leq r+2$ \cite{dokmanic2015euclidean}. This implies that when $r\ll n$, which is often the case in practice (e.g., $r=2$ or $3$ in EDMC), $\X$ and $\D$ are low-rank matrices. This allows us to utilize a rich library of tools from low-rank matrix completion, and moves us to consider the problem of Euclidean Distance Matrix Completion (EDMC). With this in mind, one technique is to directly apply matrix completion techniques on $\D$ \cite{moreira2017novel}. Let $\Omega \subset \{(i,j) \mid 1 \leq i < j \leq n\}$ denote the set of sampled indices corresponding to the strictly upper-triangular part of the distance matrix. Note that, since a distance matrix is hollow and symmetric, it suffices to consider the samples in the upper-triangular part; that is, if $D_{ij}$ is sampled, $D_{ji}$ is also assumed to be sampled. A matrix completion approach would consider the following optimization program to recover $\D$:
\begin{equation}
\begin{split}
\label{eq: Convex MC objective}
    \minimize_{\Z\in \real^{n\times n}}\quad & 
\Vert\Z\Vert_{\ast}\\
  \subjectto\quad & Z_{ij} = D_{ij} \quad \forall (i,j)\in \Omega,
	\end{split}
\end{equation}
where $\Vert\cdot \Vert_{*} = \sum_i \sigma_i $ denotes the nuclear norm, which serves as a convex surrogate for rank \cite{fazel2001rank}. The main idea of these tools is that, under some assumptions, the nuclear norm minimization program reconstructs the true low-rank squared distance matrix exactly with high probability from $\cO(nr\log^2(n))$ randomly sampled entries \cite{candes2005decoding,candes2006robust,candes2009exact,recht2010guaranteed,gross2010note}. Another set of techniques \cite{tasissa2018exact,lai2017solve} focuses on recovering the point configuration by using the Gram matrix as an optimization variable, and using only partial information from the entries in $\D$. Specifically, these works consider the following optimization program for the EDMC problem:
\begin{equation}\label{eq:edg_convex}
\begin{split}
 \minimize_{\X\in \real^{n\times n},\,\X = \X^\top,\,\X\succeq \bm{0},\X\bm{1}=\bm{0}}& \quad  \Vert\X\Vert_{\ast},
\end{split}
\end{equation}
where the constraints follow from the relation of $\X$ and $\D$ in \eqref{eqn: D to X} and \eqref{eqn: X to D}.  Due to the challenge of working with the constraints imposed by distance matrices, i.e., an entrywise triangle inequality that must be satisfied in order to remain a distance matrix, this work will follow the latter approach of optimizing over the Gram matrix. We note that, in contrast to completing the squared distance matrix $\D$, which has rank at most $r+2$, employing a minimization approach based on a Gram matrix that has rank at most $r$ implicitly enforces the constraints of the Euclidean distances. Recent works have indicated that this approach can achieve better sampling complexity than direct distance matrix completion \cite{tasissa2018exact,lai2017solve,Li2024}.

We note that theoretical guarantees for \eqref{eq:edg_convex} have been established in \cite{tasissa2021low,tasissa2018exact}, but these methods still suffer from the lack of scalability of convex techniques. A non-convex Lagrangian formulation was also proposed in \cite{tasissa2018exact}, yielding strong numerical results but lacking local convergence guarantees. The work in \cite{Nguyen2019} uses a Riemannian manifold approach to develop a conjugate gradient algorithm for estimating the underlying Gram matrix. The theoretical analysis therein shows that the squared distance matrix iterates globally converge to the true squared distance matrix at the sampled entries under three assumptions. However, the relationship between the problem parameters, such as the sampling scheme and sampled entries, and the third assumption remains unclear, as noted in Remark III.8 of the paper. In \cite{Li2024}, the authors introduce a Riemannian conjugate gradient method with line search for the EDMC problem. The paper provides a local convergence analysis for the case where the entries of the distance matrix are sampled according to the Bernoulli model, given a suitable initialization. The initialization method used is known as rank reduction, which begins with initial points embedded in a higher-dimensional space than the target dimension. While \cite{Li2024} demonstrates strong empirical results for this initialization via tests on synthetic data for sensor localization, there are no provable guarantees provided for the initialization. 

The work by \cite{li2025euclideandistancematrixcompletion} is thematically closer to ours, as it also considers a similar initialization, and proves linear convergence to
the ground truth. Therein, the core approach is based on a Burer-Monteiro factorization of the Gram matrix, and this leads to a program where the optimization variable is in terms of points. To solve the resulting optimization problem, a simple gradient descent algorithm with a line search followed by a projection is proposed. The projection is done to ensure incoherence, but 
requires as input an incoherence parameter and the maximum singular value of the true Gram matrix. We highlight a few differences in our approach to this work. First, their empirical evaluations assume known incoherence and maximum singular value. We believe 
such sharp estimates may not be necessary, and the impact of not having this information is not clear.  
We also note that we provide robustness guarantees in contrast to  \cite{li2025euclideandistancematrixcompletion}.

The work in \cite{ghosh2024} proposes a non-convex algorithm for the EDMC problem based on the reweighted least squares framework. It considers the case where distance entries are observed uniformly at random and establishes that with $O(\nu r\log{n})$ distance entries, where $\nu$ is the incoherence parameter (see Section \ref{sec:Incoherence_Discussion} for the definition of a weaker form of incoherence used in this paper), are sufficient for local convergence to the ground Gram matrix. However, \cite{ghosh2024} does not provide a provable initialization scheme or robustness guarantees for the proposed algorithm. We note that the analysis in \cite{ghosh2024} achieves optimal sample complexity, matching the lower bound established in \cite{candes2010power}. However, their results rely on a stronger incoherence condition than ours. In fact, under our milder incoherence assumption, their sample complexity aligns with ours up to constant factors.

\subsection{Contributions}
The main contributions of this paper are as follows:

\begin{enumerate} [leftmargin=*]

\item {\bf Geometric Interpretation of EDMC Incoherence:} We provide a geometric interpretation of incoherence within the specific context of the EDMC problem. We derive both lower and upper bounds for this parameter and discuss the geometries that achieve these bounds. Under a random model of the underlying points, we show that the incoherence scales logarithmically, which aligns with the scaling of standard incoherence measures.

 \item {\bf Algorithmic Framework:} We propose a novel non-convex iterative algorithm for the Euclidean Distance Matrix Completion (EDMC) problem based on Riemannian optimization. The algorithm performs first-order updates on the manifold of fixed-rank matrices and enjoys low per-iteration computational complexity. 
    \item {\bf Provable initialization scheme:} We develop a structured initialization procedure from partial distance measurements and establish an explicit error bound between the initialization and the ground truth. The method is simple to implement and only requires available measurements.
    \item {\bf Convergence guarantees, sample complexity requirements, and robustness guarantees:} We provide rigorous analysis establishing high-probability local convergence of the proposed algorithm to the ground truth configuration with near-optimal sample complexity. 
    We also derive sample complexity bounds to ensure that the initialization lies within the basin of attraction and to provide robustness guarantees against bounded noise perturbations of the underlying point cloud.
    \item{\bf Novel Analysis:} We leverage statistical tools to analyze the local behavior of the algorithm, including a restricted isometry property for a symmetric operator with coupled structure. This allows us to certify restricted invertibility of operators in the challenging setting of matrix recovery from non-orthogonal measurements.

\end{enumerate}
To the best of our knowledge, this is the first non-convex algorithm for the EDMC problem that provides provable initialization, provable convergence guarantees, robustness guarantees under noise, and a geometric interpretation of incoherence in the EDMC context.

\subsection{Notation}
The notation used in this paper is summarized in \Cref{tab:notation}. This table provides a general description of the conventions used throughout this paper, but not every assignment is a strict rule. For example, lowercase boldface, such as $\x$, is denoted as reserved for vectors; however, we extensively use the notation $\wa$ and $\va$ for certain matrices, and some vectors in the appendix are in uppercase boldface. If there is any contradiction with \Cref{tab:notation}, the notation should be clear from context.

\begin{table}[ht!]
\centering
\begin{tabular}{ll}
\toprule
\textbf{Symbol} & \textbf{Meaning} \\
\midrule
\textbf{Matrices, Vectors, and Operators} & \\

$\A$, $\B$ & Matrices (uppercase boldface) \\
$\v$ & Vectors (lowercase boldface) \\
$\mathcal{A}$ & Linear operators on matrices (calligraphic) \\
$\mathbb{V}$ & Vector spaces and subspaces (blackboard bold) \\
$\X^\top$ & Transpose of matrix $\X$ \\
$\tr(\X)$ & Trace of matrix $\X$ \\
$\la \A,\B\ra$ & Trace inner product: $\tr(\A^\top\B)$ \\
$\delta_{ij}$ & Kronecker delta \\
$X_{ij}$ & $(i,j)$-th entry of matrix $\X$ \\
$\mathcal{A}^\ast$ & Adjoint of operator $\mathcal{A}$ \\
$\one$ & Column vector of ones (size determined by context) \\
$\bm{0}$ & Zero vector or zero matrix (depending on context) \\
$\e_i$ & Standard basis vector: $1$ at $i$-th position, zeros elsewhere \\
$\e_{ij}$ & Standard matrix basis: $1$ at $(i,j)$, zeros elsewhere \\
$\Vec(\Y)$ & Column stack of matrix $\Y$ into $\real^{n^2}$ \\
$\odot$ & Hadamard (entrywise) product \\
$\mathcal{I}$ & Identity operator on matrices \\
$\I$ & Identity matrix \\
$\A\succeq\B$ & Loewner ordering: $\A-\B$ is positive semi-definite \\
$\X=\U\D\U^\top$ & Thin spectral decomposition of symmetric rank-$r$ matrix \\
 $H_{\alphab\betab}$&$\langle\wa,\wb\ra$\\
 $H^{\alphab\betab}$&$\la \va,\vb\ra$\\
 
\midrule
\textbf{Norms and Spectral Quantities} & \\

$\Vert\x\Vert_2$ & Euclidean ($\ell_2$) norm of vector $\x$ \\
$\Vert\X\Vert_\fro$ & Frobenius norm of matrix $\X$ \\
$\Vert\X\Vert$ & Operator norm (largest singular value) \\
$\Vert\X\Vert_\infty$ & Max absolute entry of $\X$ \\
$\Vert\X\Vert_\ast$ & Nuclear norm: $\sum_i\sigma_i(\X)$ \\
$\Vert\mathcal{A}\Vert$ & Operator norm of $\mathcal{A}$: $\sup_{\Vert\X\Vert_\fro=1}\Vert\mathcal{A}(\X)\Vert_\fro$ \\
$\Vert\X\Vert_{2\to\infty}$ & $2\to\infty$ norm: $\max_i\Vert\X\e_i\Vert_2$\\
$\lambda_{\max}(\X), \lambda_{\min}(\X)$ & Max/min eigenvalues of $\X$ \\
$\lambda_1(\X)\geq\cdots\geq\lambda_r(\X)$ & Ordered non-zero eigenvalues of rank-$r$ matrix, when clear\\
& from context, $\X$ is omitted \\
$\sigma_r(\Y)$ & $r$-th singular value of matrix $\Y$ \\
$\kappa$ & Condition number: $\Vert\X\Vert/\sigma_r(\X)$ \\

\midrule
\textbf{Sets and Indexing} & \\

$\mathbb{I}$ & Universal set of indices $\{(i,j):1\leq i<j\leq n\}$ \\
$\Omega$ & Random subsets of $\mathbb{I}$ \\
$\emptyset$ & Empty set \\
$\x_i$, $\x^i$ & $i$-th row and $i$-th column of $\X$, respectively,\\
& represented as column vectors. \\

\midrule
\textbf{Manifolds and Geometry} & \\

$\mfr$ & Manifold of rank-$r$ matrices \\
$\mf$ & General smooth manifolds \\
$\T$, $\T_l$ & Tangent space at $\X\in\mfr$ and at $l$-th iterate $\X_l\in\mfr$ \\
$\nabla f$ & Euclidean gradient of $f\in C^1(\Rnn)$ \\
$\mathrm{grad}\,f$ & Riemannian gradient of $f\in C^1(\mfr)$ \\
\bottomrule
\end{tabular}
\vspace{2mm}
\caption{Summary of notation used throughout the paper.}
\label{tab:notation}
\end{table}

\subsection{Organization}
The organization of this paper is as follows. In \Cref{sec: background}, we discuss the requisite background information necessary to understand the work done in this paper. This consists of a brief discussion of low-rank matrix completion and a discussion of EDMC, with further background on dual bases and first-order Riemannian methods found in Appendix~\ref{appendix: background}. \Cref{sec: R_omega alg section} is a discussion of our proposed methodology for solving the EDMC problem using geometric low-rank matrix completion ideas in the developed dual basis framework. \Cref{sec:Incoherence_Discussion} gives a detailed discussion of the EDMC-specific incoherence condition. \Cref{sec: Analysis results} discusses the underlying assumptions, convergence analysis and initialization guarantees of the proposed algorithm with most proofs deferred to the Appendices. Next, we include robustness guarantees in \Cref{subsec: robustness}. The convergence analysis leverages the discussed dual basis structure, with properties proven in Appendix~\ref{appendix: dual basis}, to get local convergence guarantees, discussed in more detail in Appendices~\ref{appendix: RIP}~and~\ref{appendix: local convergence}. We additionally provide initialization and robustness guarantees in this section, with relevant proofs in Appendices~\ref{appendix: initialization} and \ref{appendix: robustness}. \Cref{sec: Related work} discusses related geometric approaches in matrix completion, relevant work done in EDMC, and a more detailed discussion of geometric approaches to EDMC. \Cref{sec: Numerics} discusses the numerical results of this algorithm, and compares its efficacy to another algorithm in the literature. We conclude the paper in \Cref{sec: Conclusion} with a brief discussion of the work and possible future research directions.

\section{Preliminary Material}\label{sec: background}
In this section, we will provide some minor background necessary to understand the work done in the following sections. A discussion of dual bases in linear algebra and first-order Riemannian methods can be found in Appendix~\ref{appendix: background}.

\subsection{Matrix Completion}\label{subsec: matrix completion}
This work is related to low-rank matrix completion, where a subset of the entries of a low-rank ground truth matrix $\X$ is observed. Consider $\X$ as an $n \times n$ matrix for simplicity, with $\Omega\subset [n]\times[n]$ representing the set of observed indices. Here, a sampling operator $\Po:\Rnn\to\Rnn$ is introduced, which aggregates the observed entries of $\X$ projected onto specific basis elements $\e_{ij}$:
\begin{equation}\label{eqn: Po definition}
    \Po(\X) = \sum_{(i,j)\in\Omega}\langle \X,\e_{ij}\rangle \e_{ij}.
\end{equation}
If $\Omega$ does not contain any repeated indices, $\Po$ is an orthogonal projection operator. The standard low-rank matrix completion problem can be phrased as
\begin{equation*}
    \minimize_{\Y\in\Rnn} ~\rank(\Y) ~ \subjectto ~ \Po(\Y) = \Po(\X).
\end{equation*}
As minimizing the rank directly is generally a challenging problem \cite{candes2009exact,meka2008rank}, relaxations of this problem are often considered. For details on the complexity class of rank constrained problems, we refer the reader to \cite{bertsimas2022mixed}. 
Exact recovery of $\X$ from $\Po(\X)$ using a convex relaxation to the nuclear norm, such as the objective described in \eqref{eq: Convex MC objective},  is a well-studied problem \cite{candes2006robust,recht2011simpler,gross2011recovering}. This problem is at the core of matrix completion literature, and has inspired work in the completion of distance matrices \cite{lai2017solve,tasissa2018exact}. However, solving the convex problem is expensive for large matrices, which has led to the consideration of non-convex methodologies to solve the underlying problem. One approach that has received a great deal of attention is the Burer-Monteiro factorization approach, pioneered for semi-definite methods in \cite{burer2003nonlinear}, whereby a low-rank matrix $\X\in\Rnn$ can be factored into a product $\X = \A\B^\top$ for $\A,\B\in\real^{n\times r}$. Minimizing $\Vert \Po(\X) - \Po(\A\B^\top)\Vert_\fro^2$ is a common approach, and is often dealt with using alternating minimization methods in both the noiseless and noisy cases \cite{jain2013low,Hardt2014,zhang2016provable,Chen2020}. Beyond standard matrix completion, these methods have also been applied to other structured problems \cite{abbasi2023fast,abbasi2025efficient,netrapalli2013phase}.

In \cite{wei2020guarantees}, the authors develop a gradient descent algorithm to solve the low-rank matrix completion problem, allowing the reconstruction of a ground truth matrix $\X$ with known rank $r$ from partial access to entries, by leveraging the Riemannian structure of the set of fixed-rank matrices. The objective function used in \cite{wei2020guarantees} is as follows:
\begin{equation}\label{eqn: P_omega qf objective}
    \minimize_{\Y\in\Rnn} ~\langle \Y-\X,\Po(\Y-\X)\rangle ~ \subjectto ~ \rank(\Y) = r.
\end{equation}
We provide more details on \cite{wei2020guarantees} and other approaches to matrix completion in \Cref{sec: Related work}.

One of the main statistical approaches to analyzing matrix completion problems is through studying the behavior of the sampling operator $\Po$ restricted to a feasible space for recovery. We formalize this by fixing a rank-$r$ ground truth matrix $\X$. As the set of rank-$r$ matrices for a fixed $r$ forms a manifold, denoted $\mfr$, we can consider the tangent space $\T$ at $\X$ to aid our understanding of $\Po$.
Explicitly, for $\X = \U\bm{\Lambda}\U^\top$, we can characterize $\T$ as
\begin{equation}\label{eqn: T defn}
\T = \{\U \Z^\top + \Z\U^\top ~\vert~ \Z\in\real^{n\times r}\}.
\end{equation}
Intuitively, restricting $\Po$ to $\T$ and measuring the deviation of this operator from the identity measures how well $\Po$ preserves information associated with $\X$ upon measurement, and whether or not $\X$ is uniquely recoverable given the information accessed. Mathematically, this manifests in proving statements such as 
\[
\Vert\Pt\Po\Pt-c\Pt\Vert\leq \varepsilon_0,
\]
for some constant $c>0$ and some small $\varepsilon_0>0$, which depends on both the number of samples and intrinsic properties of the ground truth matrix $\X$\cite{recht2011simpler}. This property is known as the Restricted Isometry Property (RIP), and variants of this property are critical to the low-rank matrix completion and compressive sensing literature\cite{Wright_Ma_2022,foucart2013mathematical}.

\subsection{Dual Basis Approach to EDMC}\label{subsec: Dual Basis EDMC}

In the EDMC problem, using the relation \eqref{eqn: X to D}, we can relate each entry of the squared distance matrix to the Gram matrix
as follows: $D_{ij} = X_{ii} + X_{jj} - X_{ij}-X_{ji}$. We describe here briefly the dual
basis approach introduced in \cite{tasissa2018exact}. Given $\alphab = (i, j), i < j$,
we define the matrix $\wa$ as follows:
\begin{equation}\label{eqn: wa definition}
    \wa = \e_{ii}+\e_{jj}-\e_{ij}-\e_{ji}.
 \end{equation}
If we consider the set $\mathbb{I}=\{(i,j),1\leq i<j\leq n\}$, it can be checked that the set $\{ \wa \}_{\ai}$ is a non-orthogonal basis for the subspace of symmetric matrices with zero row sum, denoted $\mathbb{S} = \{\Y\in\Rnn~\vert~\Y=\Y^\top, \Y\one=\bm{0}\}$. In fact, for any two pairs of indices $\alphab,\betab\in\mathbb{I}$, we have:
\begin{equation*}
    \langle\wa,\wb\rangle = \begin{cases}
        4 & \alphab=\betab;\\
        1 & \alphab\neq\betab,~\alphab\cap\betab\neq\emptyset;\\
        0 & \alphab\cap\betab = \emptyset.
    \end{cases}
\end{equation*}
It can also easily be verified that the dimension of the linear space $\mathbb{S}$ is $L=n(n-1)/2$.
Using this basis, we can realize each entry of the squared distance matrix as the trace inner product of the Gram matrix with the basis. Formally, $D_{ij} = \langle \X, \wa \rangle$ for $\alphab = (i,j)$. Further, we can introduce the dual basis
to $\{ \wa \}_{\ai}$, denoted as $\{ \va \}_{\ai}$, and represent any centered Gram matrix $\X$ using
the following expansion:
\[
\X = \sum_{\alphab} \langle \X\,,\wa\rangle \va.
\]
The advantage of the dual basis representation is that it allows us to recast the EDMC problem as a low-rank matrix recovery problem where we observe a subset
of the expansion coefficients. In \cite{tasissa2018exact}, this dual basis formulation has been used to provide theoretical guarantees for the convex program given in \eqref{eq:edg_convex}. 

The direct form of the dual basis, based on its definition, relies on the inverse of a matrix of size $L\times L$, which requires solving a large linear system. This can be explicitly constructed as follows. Define the matrix $\H\in\real^{L\times L}$ entrywise as $H_{\alphab\betab} = \langle\wa,\wb\ra$. It is a classical result that the inverse of $\H$ can be defined entrywise as $(\H^{-1})_{\alphab\betab} = H^{\alphab\betab} = \la\va,\vb\ra$. Throughout the remainder of the paper, the superscript $H^{\alphab\betab}$ will denote the entries of $\H^{-1}$, and the subscript $H_{\alphab\betab}$ will denote the entries of $\H$. As such, we can construct any dual basis element $\va$ via
\[
\va = \sum_{\betab\in\univ}H^{\alphab\betab}\wb.
\]
To make use of the dual basis approach both in theory and applications, one of the first steps is to have a representation of the dual basis that is easier to use. In \cite{lichtenberg2023dual}, it was shown that the dual basis admits a simple explicit form
\begin{equation}\label{eqn: v_a form}
    \va = -\frac{1}{2}\left(\a\b^\top + \b\a^\top\right),
\end{equation}
where $\a = \e_i - \frac{1}{n}\one$ and $\b = \e_j - \frac{1}{n}\one$ for $\alphab = (i,j)$.
We now highlight a few operators that are related to the dual basis approach. The first one is the restricted sampling operator. Similar to low-rank matrix completion, let $\Omega\subset \univ$, and define
$\Ro:\mathbb{S}\to\mathbb{S}$ as follows:
\begin{equation*}
    \Ro(\cdot) = \sum_{\alphab\in\Omega}\langle\cdot,\wa\rangle\va.
\end{equation*}
From the bi-orthogonality relationship of the dual basis, it follows that $\Ro^2 = \Ro$ if $\Omega$ does not have repeated indices. The adjoint of $\Ro$ admits the following form:
\begin{equation*}
    \Ro^\ast(\cdot) = \sum_{\alphab\in\Omega}\la \cdot,\va\ra\wa.
\end{equation*}
Due to the lack of self-adjointness, $\Ro$ without repeated indices in $\Omega$ is not an orthogonal projection operator, and is instead an oblique projection operator. In \cite{Smith2023}, $\Ro(\X)$ is related to the sampling operator $\Po(\D)$ as follows:
\begin{equation}\label{eqn: Ro and Po equiv}
    \Ro(\X) = -\frac{1}{2} \J\Po(\D)\J,
\end{equation}
where $\J$ is as defined in \Cref{sec: Introduction}. The next operator is the restricted frame operator $\Fo:\mathbb{S}\to\mathbb{S}$, first studied in \cite{tasissa2018exact}, and defined as
\begin{equation}\label{eqn: F_omega}
    \Fo(\cdot) = \sum_{\alphab\in\Omega}\la\cdot,\wa\ra\wa.
\end{equation}
This operator is self-adjoint, positive semi-definite, but unlike $\Ro$, does not reference the dual basis. We note that this operator is central to the analysis of the algorithm in \cite{Li2024}.

\section{The Riemannian Dual Basis Approach to EDMC} \label{sec: R_omega alg section}

With the goal of translating the standard matrix completion problem to Gram matrix completion of a ground truth matrix $\X\in \S$, where $\S= \{\Y\in\Rnn \vert \Y = \Y^\top, \Y \one = \bm{0}\}$, the most direct adaptation of the work conducted in \cite{wei2020guarantees} would be defining an objective function by analogy to \eqref{eqn: P_omega qf objective} as follows:
\begin{equation*}
    \minimize_{\Y\in\mathbb{S}} ~\langle \Y-\X,\Ro(\Y-\X)\rangle ~\subjectto ~\mathrm{rank}(\Y) = r.
\end{equation*}
However, a notable challenge arises: computing the Euclidean gradient of the objective function necessitates unavailable information in the form $\langle \X,\va\rangle$ from $\Ro^\ast(\X)$ as
\begin{equation*}
    \nabla_{\Y} \left(\langle \Y-\X,\Ro(\Y-\X)\rangle\right) = \Ro(\Y-\X)+\Ro^\ast(\Y-\X),
\end{equation*}
where $\nabla_{\Y}$ denotes the gradient with respect to $\Y$. This quantity is inaccessible under the given problem statement, since each $\va$ depends on all $\wa$ through the relation $\va = \sum_{\ai}H^{\alphab\betab}\wb$. Consequently, the computation $\langle\X,\va\rangle =  \sum_{\ai}H^{\alphab\betab}\left\langle\X,\wa \right\rangle $ requires complete knowledge of all pairwise distances. Therefore, evaluating $\Ro^\ast(\X)$ necessitates access to all information. To circumvent this difficulty, there has been exploration into self-adjoint alternatives to $\Ro$ \cite{Tasissa2021,tasissa2018exact,Smith2023}. The novel surrogate introduced in this work, denoted $\Mo$, allows for the definition of an objective function in analogy to \eqref{eqn: P_omega qf objective}.

Let $\Omega\subset \univ$ be a random subset of $\univ$ sampled with uniform Bernoulli probability $p\in [0,1]$. We now define $\Mo:\S\to\S$ as follows:
\begin{equation}\label{eqn: Mo definition}
    \Mo(\cdot) = \sum_{\alphab,\betab\in\Omega}C_{\alphab\betab}\wa\la\va,\vb\ra\la\cdot,\wb\ra,
\end{equation}
where $C_{\alphab\alphab} = p$ for all $\alphab$, and $C_{\alphab\betab} = 1$ for all $\alphab\neq \betab$. Previous literature introduced an unscaled form of this operator, i.e., $C_{\alphab\alphab} = 1$, which is equivalent to $\RRo$ \cite{Smith2023}. However, $\RRo$ does not concentrate around $p^2\Id$. The concentration around $p^2\Id$ is a necessary condition for convergence of the algorithm presented in \cite{Smith2023}, and to rectify this discrepancy, a different operator must be considered. The diagonal re-scaling present in $\Mo$ enforces that $\Ebb[\Mo] = p^{2}\mathcal{I}$, with proof provided in \Cref{lem: Expectation of Mo}.  Additionally, $\Mo$ is self-adjoint, allowing us to define the following objective function for the EDMC problem using this operator:
\begin{equation}\label{eqn: M_omega objective function}
    \minimize_{\Y\in\mathbb{S}} ~\frac{1}{2}\langle \Y-\X,\Mo(\Y-\X)\rangle ~\subjectto ~\mathrm{rank}(\Y) = r.
\end{equation}
This object is a true quadratic form with a symmetric operator, and its Euclidean gradient is given solely by $\Mo(\Y-\X)$. As such, it can be approached identically to \eqref{eqn: P_omega qf objective} following the principles outlined in Appendix~\ref{appendix: background}. To perform this first-order retraction method, we need to define the tangent space at each iterate $\X_l\in\mfr$. Given $\X_l = \U_l\bm{\Lambda}_l\U_l^\top$ for $\U_l \in \real^{n\times r}$, we can define $\T_l$ as in \eqref{eqn: T defn}. The projection onto $\T_l$, denoted $\Ptl$, can be given by 
\[
\Ptl = \Pul\Y + \Y\Pul - \Pul\Y\Pul,
\]
where $\Pul = \U_l\U_l^\top$, allowing for easy computation of the Riemannian gradient of the objective function. Now that the tangent space is defined, we define the retraction map, known as the hard thresholding operator $\mathcal{H}_r:\T_l\to\mfr$, as follows:
\begin{equation}\label{eqn: hard threshold operator}
    \mathcal{H}_r(\Y) = \sum_{i=1}^r \lambda_i(\Y) \U_i\U_i^\top,
\end{equation}
where $\U_i$ is the $i$-th eigenvector of $\Y$ corresponding to the eigenvalue with the $i$-th largest magnitude $\lambda_i(\Y)$. We note that for matrices $\Y$ with $\rank(\Y)\geq r$, we have $\rank\left(\mathcal{H}_r(\Y)\right) = r$. We can now define \Cref{alg: M_omega descent}, the main object of study in this work:

\begin{algorithm}
\caption{Dual Basis Riemannian EDMC (DBRE)}\label{alg: M_omega descent}
\begin{algorithmic}[1]
\STATE \textbf{Input:} $\Omega$, $\{\la \X,\wa\ra\}_{\alphab\in\Omega}$, $\rank(\X)=r$
\STATE\textbf{Initialization:} $\X_0 =\U_0\D_0\U_0^\top$
\FOR{$l = 0,1,\cdots$}
    \STATE  $\G_l = \Mo(\X-\X_l)$
    \STATE  $\alpha_l = \frac{\left\Vert \Ptl\G_l\right\Vert_\fro^2}{\langle \Ptl\G_l,\Mo\Ptl\G_l\rangle }$
    \STATE $\W_l = \X_l + \alpha_l \Ptl\G_l$
    \STATE $\X_{l+1} = \mathcal{H}_r(\W_l)$
\ENDFOR
\STATE\textbf{Output:} $\X_\mathrm{rec}$: Estimated Gram matrix.
\end{algorithmic}
\end{algorithm}

We note that this algorithm is a standard projected gradient descent scheme on $\mfr$; the novelty presented in this work is in the construction of the measurement operator $\Mo$ to build the proper objective function, which enables the convergence analysis of Algorithm~\ref{alg: M_omega descent}, presented in \Cref{sec: Analysis results}. 

In the approach seen in \Cref{alg: M_omega descent}, the thin spectral decomposition in the gradient descent scheme is the most expensive, especially when $n$ is large. As described previously, the authors in \cite{wei2020guarantees} found an efficient way to reduce the computational complexity of this decomposition from $\mathcal{O}(rn^2)$ to $\mathcal{O}(r^3+r^2n)$, substantially reducing the cost per iteration, which we implement as well. We note that in \Cref{alg: M_omega descent}, the reconstruction of the ground truth Gram matrix $\X$ is equivalent to the reconstruction of $\D$, as there is a one-to-one correspondence between $\X$ and $\D$ through \eqref{eqn: X to D}.

\begin{remark}
    We wish to provide an interpretation of the operators $\Mo$ and $\RRo$, and why they are considered. Although $\Fo$ concentrates around its expectation through a Bernstein inequality---see \Cref{lem: Concentration of Fo}---its expectation is not the identity operator $\mathcal{I}$, and as such $\mathcal{F}_{\univ}\neq \Id$. In fact, the spectrum of $\FI$ is known to be equivalent to the spectra of $\H$, and thus $\lambda_{\max}(\FI) = 2n$ \cite{lichtenberg2023dual}. As such, neither $\Vert\Fo-p\mathcal{I}\Vert$ nor $\Vert\Pt\Fo\Pt-p\Pt\Vert$ are small. In particular, $\Vert\Pt\Fo\Pt-p\Pt\Vert = \cO(1)$, but as the operator needs to exhibit RIP for the analysis of an algorithm like \Cref{alg: M_omega descent}, this bound is not tight enough to guarantee convergence. We can instead consider re-scaling the geometry of the linear space $\S$ through some function $f$ so that $\Ebb[\Fo(f(\cdot))]\approx p\mathcal{I}$. First, define $\So:\Rnn\to \real^{L}$ as $(\So(\X))_{\alphab} = \la\X,\wa\ra$ for $\alphab\in\Omega$, and $0$ otherwise. As such, one can show that $\Fo = \So^\ast\So$, where $\So^\ast:\real^L\to\real^{n\times n}$ is the adjoint of $\So$. To re-scale $\Fo$, one can instead consider $\So^\ast \H^{-1}\So$. This rescaling is done with $\H^{-1}$ to make it so that $\So^\ast\H^{-1}\So = \Id$ when $\Omega = \univ$. One can compute out $\So^\ast \H^{-1}\So$ and show that
    \[
    \So^\ast \H^{-1}\So= \RRo.
    \]
     Further investigation in \Cref{lem: Expectation of Mo} validates the necessity of considering a diagonally reweighted variant of $\RRo$ to ensure concentration around $\mathcal{I}$, resulting in $\Mo$. 
\end{remark}

\subsection{Implementation Efficiency}

We use recent advances in Riemannian optimization from
\cite{wei2020guarantees} and \cite{Smith2023} to develop an efficient implementation of the proposed algorithm. Computation of $\Ro(\X)$  and $\Mo(\X)$ can be done efficiently, with minimal complexity per iteration. For $\Ro$, a given iterate $\X_l$ can be easily translated to its distance matrix $\D_l$ via \eqref{eqn: X to D}, and through \eqref{eqn: Ro and Po equiv}, $\Ro(\X_l)$ can be computed in $\mathcal{O}(m)$ operations, for $\vert\Omega\vert = m$. First, we note that $\Mo(\X_l) = \RRo(\X_l) - \Vert\va\Vert_\fro^2(1-p)\Fo(\X_l)$ and $\Vert\va\Vert_\fro$ is constant for all $\alphab\in\univ$ (\Cref{lem: H and H^-1 eigvals}). This can be seen as follows:
\begin{align}
    \Mo(\cdot) &= \sum_{\alphab,\betab\in\Omega}C_{\alphab\betab}\wa\la\va,\vb\ra\la\cdot,\wb\ra\nonumber\\ \nonumber
    & = \sum_{\substack{\alphab,\betab\in\Omega\\\alphab=\betab}}C_{\alphab\betab}\wa\la\va,\vb\ra\la\cdot,\wb\ra+ \sum_{\substack{\alphab,\betab\in\Omega\\\alphab\neq\betab}}C_{\alphab\betab}\wa\la\va,\vb\ra\la\cdot,\wb\ra \\ \nonumber
    &= \sum_{\alphab\in\Omega}C_{\alphab\alphab}\wa\la\va,\va\ra\la\cdot,\wa\ra + \sum_{\substack{\alphab,\betab\in\Omega\\\alphab\neq\betab}}C_{\alphab\betab}\wa\la\va,\vb\ra\la\cdot,\wb\ra \\ \nonumber
    & = p\Vert\va\Vert_\fro^2\Fo(\cdot) + \sum_{\substack{\alphab,\betab\in\Omega\\\alphab\neq\betab}}\wa\la\va,\vb\ra\la\cdot,\wb\ra \\ \nonumber
    & = p\Vert\va\Vert_\fro^2\Fo(\cdot) + \sum_{\alphab,\betab\in\Omega}\wa\la\va,\vb\ra\la\cdot,\wb\ra  - \sum_{\alphab\in\Omega}\la\cdot,\wa\ra\la\va,\va\ra\wa\\
    & = p\Vert\va\Vert_\fro^2\Fo(\cdot) +\RRo(\cdot) - \Vert\va\Vert_\fro^2\Fo(\cdot), \label{eq:form_of_Momega}
\end{align}
as expected. It is known that $\RRo(\X_l)$ is $ \cO(m)$ sparse and requires $\cO(m)$ operations to compute \cite{Smith2023}. The argument is outlined as follows. Let $\mathcal{T}:\Rnn\to\Rnn$ denote the map defined by \eqref{eqn: X to D} and let $\mathcal{T}^\ast$ denote its adjoint. It was shown in \cite{Smith2023} that, correcting for a previously missing minus sign, for a Gram matrix $\X$,
\[
\RRo(\X)  = -\frac{1}{4}\mathcal{T}^\ast\left(\Po\left(\J\Po(\mathcal{T}(\X))\J\right)\right).
\]
For any matrix $\Y$, both $\Po(\Y)$ and $\J\Po(\Y)\J$ are computable in $\cO(m)$ operations. The accessible information in the EDMC problem is of the form $\Po(\mathcal{T}(\X))$. Furthermore, $\mathcal{T}^\ast\left(\Po(\Y)\right)$ for any $\Y$ is computable in $\cO(m)$ operations as well. 

Next, $\Fo(\X_l)$ is efficiently computable in $\cO(m)$ operations as each matrix $\wa$ has $4$ non-zero entries, allowing for easy computation given $\{\langle\X,\wa\rangle\}_{\alphab\in\Omega}$. As such, $\Fo(\X_l)$ is $\mathcal{O}(m)$ sparse. Using the fact that $\RRo(\cdot)$  and $\Fo(\cdot)$ are sparse, it can be easily argued that the sum of the three terms in \eqref{eq:form_of_Momega} preserves a common sparsity pattern, and it can be computed in $\cO(m)$ operations. Therefore $\Mo(\X_l)$, and thus $\G_l$ in Step 4 of \Cref{alg: M_omega descent}, is computable in $\cO(m)$ operations. 

Step 5 can be computed in $\cO(n^2)$ operations, as $\Pt\G_l$ is a dense matrix. 
Some calculations yield that Steps 6 and 7 can be computed with $n^2r+\mathcal{O}(nr^2+r^3)$ \cite{wei2020guarantees}, giving a total cost per iteration of $n^2r + \mathcal{O}(m+nr^2+r^3)$. Note that the dominant cost is $n^2r$, which is less expensive than computing Step 7 using the truncated singular value decomposition directly. Although both approaches have the same asymptotic complexity, the latter incurs a significantly higher constant factor (e.g., a factor of $6$ or $14$ depending on the choice of algorithm; see, for example, Figure 8.6.1 in \cite{golub_loan}).

\section{Geometric Interpretation of EDMC Incoherence}
\label{sec:Incoherence_Discussion}

In pathological cases, the ground truth matrix $\X$ may exhibit a sparse representation in the basis $\{\wa\}_{\ai}$, which could lead to challenges in its recovery from sampled measurements. While the concept of incoherence is well-established in the standard matrix completion literature, the condition specific to the EDMC problem slightly differs in structure and admits a natural geometric interpretation. This section is devoted to a detailed examination of this geometric perspective. We will state more formally the incoherence assumptions in \Cref{sec: Analysis results}, but we will first introduce one of the conditions below. We say that a rank-$r$ Gram matrix $\X\in\Rnn$ is $\nu$-incoherent with respect to $\{\wa\}_{\ai}$ if the following statement holds:

\begin{equation}\label{eqn: wa incoherence equation}
    \max_{\ai}\Vert\,\,\Pu\wa\Vert_\fro\leq\sqrt{\frac{4\nu r}{n}}.
\end{equation}
We remark that the above is inspired by the standard incoherence condition, which, up to a scaling factor, states that
\begin{equation}\label{eqn: standard incoherence main text}
    \max_{i\in[n]}\Vert\,\,\Pu\e_i\Vert_2\leq \sqrt{\frac{4\nu r}{n}}.
\end{equation}
The standard incoherence assumption, shown in \eqref{eqn: standard incoherence main text}, is prevalent throughout the matrix completion literature and is a measure of ``entrywise diffuseness'' in the ground truth matrix. Further discussion of standard matrix incoherence can be seen in \cite{candes2009exact}.

The incoherence condition introduced in \eqref{eqn: wa incoherence equation} can be interpreted in terms of the underlying point cloud data. For the specific case of the EDMC problem, \eqref{eqn: wa incoherence equation} can be expanded as follows for $\alphab = (i,j)$ with $i<j$: 
\begin{align}
    \Vert\Pu\wa\Vert_\fro^2 & = \langle\Pu\wa,\Pu\wa\rangle \nonumber\\
    & = \tr\left(\wa\wa\U\U^\top\right)\nonumber\\
    &=\tr\left( 2\wa\U\U^\top\right)\nonumber\\
        & =2 \left\langle \e_{ii}-\e_{ij}+\e_{jj}-\e_{ji},\U\U^\top\right\rangle\nonumber\\
    &=2\left((\U\U^\top)_{ii} + (\U\U^\top)_{jj} - (\U\U^\top)_{ij} - (\U\U^\top)_{ji}\right)\nonumber\\
    &= 2\left({\u_i}^\top\u_i + {\u_j}^\top\u_j -{\u_j}^\top\u_i -{\u_i}^\top\u_j \right)\nonumber\\
    &= 2\left({\u_i}^\top\left(\u_i-\u_j\right) +{\u_j}^\top\left(\u_j-\u_i\right)\right)\nonumber\\ 
    & = 2\left(\u_i-\u_j\right)^\top\left(\u_i-\u_j\right).\label{eqn: Puwa and u_i-u_j}
\end{align}
The incoherence condition can then equivalently be stated as 
\begin{equation} \label{eq:coherence_u}
\max_{(i,j):i<j}\,\,
\left(\u_i-\u_j\right)^\top\left(\u_i-\u_j\right) \leq2\frac{\nu r}{n}.
\end{equation}
The next Lemma provides the lower and upper bounds for $\nu$. 
\begin{lem}
For the incoherence condition in \eqref{eq:coherence_u}, $\nu$ is bounded below by $1+\frac{1}{n-1}$ and above by $2\frac{n}{r}$. 
\end{lem}

\begin{figure}[ht!]
  \centering
 
    \includegraphics[scale=1.5]{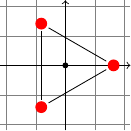}
    \caption{Visualizing the rows of $\U$ that lead to a minimal incoherence parameter.}
        \label{fig:low_coherence}

    \end{figure}

\begin{figure}
    \centering
    \includegraphics[ scale=0.75]{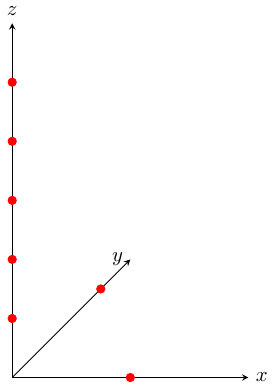}
   \caption{Example of a set of points with a maximal incoherence parameter.}
  \label{fig:high_coherence}

\end{figure}

\begin{proof}

We consider $\sum_{(i,j):i<j}  \left(\u_i-\u_j\right)^\top\left(\u_i-\u_j\right)$.
Note that $\sum_{i} (\u_i)^\top\u_i = \tr(\U\U^\top)=\tr(\U^\top\U)=r$. Since we assume centered configurations, $\U^\top\bm{1}=\bm{0}$. It then follows that $\sum_{i,j} (\u_i)^\top\u_j= (\sum_{i} \u_i)^\top \sum_{j} \u_j = 0$. Using these two relations, we obtain:
\[
\sum_{(i,j):i<j}  \left(\u_i-\u_j\right)^\top\left(\u_i-\u_j\right) = \frac{1}{2}
\sum_{(i,j)} \left(\u_i-\u_j\right)^\top\left(\u_i-\u_j\right) = \frac{1}{2}2nr = nr.
\]
The above equality shows that the sum of $L$ terms is $nr$. Therefore, the maximum summand must be at least $\frac{nr}{L}$. In particular, we have:
\[
\max_{(i,j):i<j} \left(\u_i-\u_j\right)^\top\left(\u_i-\u_j\right) \geq \frac{nr}{L} = 2\frac{r}{n}\left(1+\frac{1}{n-1}\right).
\]
Therefore, the minimum value of the incoherence parameter $\nu$ is $1+\frac{1}{n-1}$. To find the maximum value of the incoherence, we use the parallelogram inequality, $\left(\u_i-\u_j\right)^\top\left(\u_i-\u_j\right) \leq 2||\u_i||^2+2||\u_j||^2 \leq 4$. Therefore, the upper bound for $\nu$ is $2\frac{n}{r}$.
\end{proof}
\begin{remark}
To show that the lower bound for the incoherence can be attained, we consider the following example:
\[
\U = \sqrt{\frac{2}{3}} \begin{bmatrix}
1 & 0 \\
-\frac{1}{2} & \frac{\sqrt{3}}{2} \\
-\frac{1}{2} & -\frac{\sqrt{3}}{2}
\end{bmatrix}.
\]
Up to the scaling factor of $\sqrt{\frac{2}{3}}$, the rows of $\U$ correspond to the vertices of an equilateral triangle inscribed in the unit circle. It can be easily verified that this attains the lower bound on incoherence. For the upper bound,  a simple example is the matrix $\U\in \real^{n\times 3}$, where the first two columns are the standard basis vectors $\e_1$ and $\e_2$, respectively, and the third column is a unit vector which is zero in its first two entries. Any set of points generated from this $\U$ lies 
entirely along the 
z-axis, except for two points, which lie on the x- and 
y-axes, respectively. Figures \ref{fig:low_coherence} and \ref{fig:high_coherence}  provide a visual illustration of these examples.

\end{remark}


Next, we aim to state the incoherence condition in terms of the points. To do this, we assume that $\rank(\P)=r$. Using the fact that $\U\U^\top$ is a projection operator onto the column space of $\P$, it follows that the unique orthogonal projector onto this range is defined as:
$\bm{U}\bm{U}^\top = \bm{P}(\bm{P}^\top \bm{P})^{-1}\bm{P}^\top.$
By substituting this definition directly into the expanded trace formulation, we notice the following:
\begin{align*}
 &(\bm{U}\bm{U}^\top)_{ii} + (\bm{U}\bm{U}^\top)_{jj} - (\bm{U}\bm{U}^\top)_{ij} - (\bm{U}\bm{U}^\top)_{ji}\\
&= \bm{e}_i^\top \bm{P}(\bm{P}^\top \bm{P})^{-1}\bm{P}^\top \bm{e}_i + \bm{e}_j^\top \bm{P}(\bm{P}^\top \bm{P})^{-1}\bm{P}^\top \bm{e}_j - 2\bm{e}_i^\top \bm{P}(\bm{P}^\top \bm{P})^{-1}\bm{P}^\top \bm{e}_j \\
&= \bm{p}_i^\top (\bm{P}^\top \bm{P})^{-1}\bm{p}_i + \bm{p}_j^\top (\bm{P}^\top \bm{P})^{-1}\bm{p}_j - 2\bm{p}_i^\top (\bm{P}^\top \bm{P})^{-1}\bm{p}_j \\
&= (\bm{p}_i - \bm{p}_j)^\top (\bm{P}^\top \bm{P})^{-1} (\bm{p}_i - \bm{p}_j).
\end{align*}
The norm generated by $(\P^\top\P)^{-1}$ is known in the statistics literature as the Mahalanobis distance \cite{mahalanobis2018generalized}. As we can see above, the Mahalanobis distance emerges algebraically in this expression, 
establishing the direct equivalence between the projected Frobenius norm of the measurement operator and the scaled pairwise Mahalanobis distance of the spatial coordinates.

This indicates that our incoherence condition can be reinterpreted as
\begin{equation}\label{eq:main_coherence_points}
\max_{(i,j):i<j}\,\, \left(\p_i-\p_j\right)^\top
    (\P^\top\P)^{-1}
    \left(\p_i-\p_j\right) \leq 2\frac{\nu r}{n}.
\end{equation}
We first start our interpretation for the case where $(\P^\top\P)^{-1}$ is the identity matrix. In this setting, for any pair $(i,j)$, the expression $\left(\p_i-\p_j\right)^\top\left(\p_i-\p_j\right)$ is the squared Euclidean
distance between the points $\p_i$ and $\p_j$. Hence, the incoherence can be directly linked to
the maximum distance among the points. In the general case, the above quantity suggests that incoherence serves as a measure of how the displacement vectors $\p_i-\p_j$ align with the principal components of the embedding. In particular, for a fixed choice of $(\P^\top\P)^{-1}$, varying the matrix $\U$ leads to different sets of points. If the displacement vectors tend to align with directions corresponding to the smallest principal components (i.e., those with the lowest variance), the incoherence is expected to be high. Conversely, if they align more with the dominant components (those with the highest variance), the incoherence tends to be low. In essence, high incoherence indicates that certain pairs of points are stretching significantly in directions where the embedding space has low variance. Using the variational characterization, note that \[\left(\p_i-\p_j\right)^\top
    (\P^\top\P)^{-1}
    \left(\p_i-\p_j\right)\leq \lambda_1((\P^\top\P)^{-1}) \left(\p_i-\p_j\right)^\top \left(\p_i-\p_j\right).\] Noting that 
$\lambda_1((\P^\top\P)^{-1}) = \frac{1}{\lambda_r}$, we can also state the incoherence condition in \eqref{eq:main_coherence_points} as:
\begin{equation}\label{eqn: simplified incoherence}
    \left(\p_i-\p_j\right)^\top\left(\p_i-\p_j\right)\leq 2\frac{\nu r}{n}\lambda_r.
\end{equation}
We note that these statements are not equivalent, merely that this simpler statement implies the original incoherence condition. Continuing with the simplified incoherence condition in \eqref{eqn: simplified incoherence}, we seek to derive an upper bound on $\nu$ in terms of other geometric properties of $\P$, or spectral properties of $\X$.
\noindent
First, notice that
\begin{align*}
    \left(\p_i-\p_j\right)^\top\left(\p_i-\p_j\right) 
    &\leq \max_{ij}\left\Vert\p_i-\p_j\right\Vert^2_2 = \Vert\D\Vert_\infty.
\end{align*}
\noindent
As we seek a constant $\nu$ such that \eqref{eqn: simplified incoherence} is satisfied for all $(i,j)\in\mathbb{I}$, we can see that this will be satisfied if
\[
2\frac{\nu r}{n}\lambda_{r}(\X)\leq\Vert\D\Vert_\infty.
\]
\noindent
This yields the following upper bound for $\nu$, in terms of a geometric constant and a spectral constant:
\begin{equation}\label{eqn: incoherence upper bound}
    \nu \leq \frac{n}{2r}\frac{\Vert\D\Vert_\infty}{\lambda_{r}(\X)}.
\end{equation}
\noindent
In the following lemmas, we show that data drawn from sub-Gaussian isotropic distributions exhibit the property that $\nu\ll \cO(n)$ with high probability, indicating that the recovery results derived in \Cref{sec: Analysis results} will hold with reasonable sample complexity for a wide category of data distributions.

\begin{lem}\label{lem: vershynin covariance estimate}\cite[Page 31]{Vershynin_2012}
    Let  $\{\p_i\}_{i=1}^n\sim \mu$ where $\mu$ is a sub-Gaussian probability measure defined on $\real^r$, and let $\P = [\p_1\, \cdots\,\p_n]^\top\in\real^{n\times r}$. Define the covariance matrix of $\mu$ as $\bm{\Sigma}$. If $n\geq C(t/\varepsilon)^2r$ for some constant $C>0$, then with probability at least $1-2\exp(-t^2n)$ 
    \[
    \left\Vert\frac{1}{n}\P^\top\P-\bm{\Sigma}\right\Vert\leq\varepsilon\Vert\bm{\Sigma}\Vert.
    \]
\end{lem}

\noindent Let us now assume that $\mu$ is both sub-Gaussian and isotropic, that is, $\bm{\Sigma} = \I$. Furthermore, as we are interested in point clouds satisfying $\P^\top\one=\bm{0}$, we consider mean-zero distributions. As such, we provide the following lemma:

\begin{lem}\label{lem: hwi proof for distance calc}
    Let $\{\p_i\}_{i=1}^n\subset\real^r$ be a collection of points drawn i.i.d. from an isotropic sub-Gaussian distribution $\mu$. Furthermore, let $\Ebb[\p_i]=\bm{0}$, and assume each coordinate of $\p_i$ is independent. Let $\Vert\p_i\Vert_{\psi_2}\leq K$, where $\Vert\cdot\Vert_{\psi_2}$ is the sub-Gaussian norm. Then for every $i,j\in[n]\times[n]$ with $i\neq j$, with probability at least $1-Cn^{-2}$, 
    \[
   \left\vert \left\Vert\p_i-\p_j\right\Vert_2^2-2r\right\vert\leq 4K^2\sqrt{r}\log{n},
    \]
    where $C>0$ is an absolute constant.
\end{lem}
\begin{proof}
We first notice that, as $\p_i$ and $\p_j$ are independent and $\Ebb[\p_i]=\zeros$,
\begin{align*}
    \Ebb\left[\left(\p_i-\p_j\right)^\top\left(\p_i-\p_j\right)\right] & =  \Ebb\left[\left\Vert\p_i\right\Vert^2_2\right]-\Ebb\left[{\p_j}^\top\p_i\right]-\Ebb\left[{\p_j}^\top\p_i\right]+\Ebb\left[{\p_j}^\top\p_j\right]\\
    & = \Ebb\left[\left\Vert\p_i\right\Vert^2_2\right] + \Ebb\left[\left\Vert\p_j\right\Vert^2_2\right]- 2\Ebb\left[{\p_i}\right]^\top\Ebb\left[\p_j\right] \\
    & = \Ebb\left[\left\Vert\p_i\right\Vert^2_2\right] + \Ebb\left[\left\Vert\p_j\right\Vert^2_2\right]\\
    & = 2\Ebb\left[\left\Vert\p_i\right\Vert^2_2\right]\\
    & = 2\Ebb\left[\tr\left({\p_i}{\p_i}^\top\right)\right]\\
    & = 2~\tr\left(\Ebb\left[\p_i{\p_i}^\top\right]\right)\\
    & = 2~\tr(\I) = 2r,
\end{align*}
where the second and fourth lines follow from the independence of $\p_i$ and $\p_j$, the third line follows from the fact that $\Ebb[\mu]=0$, and the seventh line follows from the fact that $\mu$ is isotropic, i.e. $\bm{\Sigma} = \I$. 
    The lemma statement follows from a simple application of the Hanson-Wright inequality, seen in \Cref{thm: Hanson Wright}. First, given 
    vectors $\p_i,\p_j$, notice that
    \begin{equation*}
        \begin{pmatrix}
             {\p_i} & {\p_j}
        \end{pmatrix}
        \underbrace{\begin{pmatrix}
            \I &-\I\\
            -\I & \I
        \end{pmatrix}}_{\A}
        \begin{pmatrix}
             {\p_i} \\ {\p_j} 
        \end{pmatrix} = \left(\p_i-\p_j\right)^\top\left(\p_i-\p_j\right).
    \end{equation*}
    Previously, we have shown that $\Ebb\left[(\p_i-\p_j)^\top(\p_i-\p_j)\right]=2r$. Furthermore, $\Vert\A\Vert_\fro^2=4r$, and $\Vert\A\Vert\leq 2$ by the Gershgorin circle theorem. The result follows from an application of \Cref{thm: Hanson Wright}.
\end{proof}

\noindent We can use the above lemma to prove an upper bound on the incoherence of an isotropic, sub-Gaussian distribution, presented in the following lemma:
\begin{lem}\label{lem: incoherence ub lemma}
    Let  $\{\p_i\}_{i=1}^n\sim \mu$ where $\mu$ is an isotropic sub-Gaussian probability measure defined on $\real^r$ with sub-Gaussian norm $\Vert\mu\Vert_{\psi_2}\leq K$, and let $\P = [\p_1\, \cdots\,\p_n]^\top\in\real^{n\times r}$. Then with probability at least $1-Dn^{-2}$ for an absolute constant $D>0$, the incoherence parameter $\nu$ of $\X = \P\P^\top$ can be upper bounded by
    \[
    \nu \leq \cO\left(K^2\frac{\log{n}}{\sqrt{r}}\right)
    \]
\end{lem}
\begin{proof}

From \Cref{lem: hwi proof for distance calc}, we have that $\Vert\D\Vert_\infty\leq 2r+4K^2\sqrt{r}\log{n}$ with high probability for a sub-Gaussian isotropic distribution. Furthermore, from \Cref{lem: vershynin covariance estimate}, $\lambda_r(\X) \geq cn$ for some $\cO(1)$ constant $c>0$ with high probability. As such, we can see that the incoherence constant can be upper-bounded using \eqref{eqn: incoherence upper bound} by
\begin{align*}
    \nu&\leq \frac{n}{2r}\frac{2r+4K^2\sqrt{r}\log{n}}{cn}\\
    &\leq 1+\frac{2K^2\log{n}}{c\sqrt{r}}\\
    &=\cO\left(K^2\frac{\log{n}}{\sqrt{r}}\right),
\end{align*}
thus concluding the proof.
\end{proof}

We note that \Cref{lem: incoherence ub lemma} shows, with high probability, that the incoherence constant remains in a regime where it does not degrade the recovery guarantees established in \Cref{sec: Analysis results} for data generated from sub-Gaussian distributions. We note that this result is very similar to the condition derived in \cite{candes2009exact} for the incoherence of matrices in the random orthogonal model.

We comment that the analysis in this section exclusively pertained to data generated from sub-Gaussian isotropic measures for expositional simplicity. These techniques can be extended to anisotropic sub-Gaussian measures, and one can show the resulting bound is $\nu = \cO(\kappa+K^2\log{n}/\sqrt{r})$, where $\kappa$ is the condition number of $\X$. We provide a proof of this result in \Cref{lem: anisotropic incoherence}.

\begin{remark}
We now provide a geometric interpretation of \eqref{eqn: standard incoherence main text} for positive semidefinite matrix completion. In this case, matrices admit a decomposition $\X = \P\P^\top = \U\bm{\Lambda}\U^\top$, so expanding \eqref{eqn: standard incoherence main text} we obtain:
\begin{align*}
    \Vert\Pu\e_i\Vert_2^2 &= \e_i^\top\Pu\e_i= 
     = \tr(\e_i^\top(\P(\P^\top\P)^{-1}\P\e_i) 
     = \p_i^\top(\P^\top\P)^{-1}\p_i
\end{align*}
As such, standard incoherence for positive semidefinite matrix completion represents the maximum standard Mahalanobis distance of each point, rather than the pairwise Mahalanobis distance between points in the geometric framework.

\end{remark}

\subsection{Finer Interpretation of EDMC Incoherence and Applications}
Throughout this work, we have treated incoherence as an index-by-index bound; that is to say, we only consider terms such as $\Vert\Pu\wa\Vert_\fro^2$. We wish to investigate this in more detail now. The main technical problem that the incoherence assumption provides a solution for is in the variance estimations used in concentration inequalities, such as in \Cref{thm: RIP of Mo}, for example. This variance estimate comes from a Gershgorin-style upper bound on the matrix $\tilde{\H} = [\la\Pu\wa,\Pu\wb\ra]\in\real^{L\times L}$, seen in \Cref{lem: Bound for largest eigval of H tilde}. The eigenvalue bound leverages the fact that, if $\alphab\cap\betab=\emptyset$, $\la\Pu\wa,\Pu\wb\ra = 0$, and for the other terms we use \Cref{asp: Incoherence assumption} in tandem with Cauchy-Schwarz to get a uniform bound on the non-zero entries. This yields an upper bound that is used to estimate the variance term in the concentration inequalities. We argue here that a more fine-grained representation of incoherence could potentially sharpen incoherence results and lead to more geometrically-optimal sampling strategies in the future.

For the Gershgorin estimate, we need to estimate $\vert\la\Pu\wa,\Pu\wb\ra\vert$ for all non-zero entries of $\tilde{\H}$. Without loss of generality, we assume that $\alphab = (i,j)$ and $\betab = (i,k)$ for $i,j,k\in[n]$. Following a nearly identical chain of computations as in the derivation of \eqref{eq:coherence_u}, one can show that 
\[
\vert\la\Pu\wa,\Pu\wb\ra\vert = \left\vert\left(\p_i-\p_j\right)^\top(\P^\top\P)^{-1}\left(\p_i-\p_k\right)\right\vert.
\]
This interpretation indicates that what might be more relevant to variance minimization is sampling more orthogonal angles with respect to a whitened dataset, rather than just considering lengths. This could lead to more optimal non-uniform sampling techniques for solving the EDMC problem.

\section{Theoretical Analysis}\label{sec: Analysis results}
In this section, we will provide the main results of this work, which are the local convergence and recovery guarantees for \Cref{alg: M_omega descent}, presented in Theorems~\ref{thm: Mo Local convergence} and ~\ref{thm: Mo Recovery Guarantee}. Prior to this, we formally state our incoherence assumptions, expanding upon the assumption first described in \Cref{sec:Incoherence_Discussion}:

\begin{asp}[Incoherence assumption]\label{asp: Incoherence assumption}
    Let $\X\in\Rnn$ be a rank-$r$ matrix with eigenvalue decomposition $\X = \U\D\U^\top$. We assume that $\X$ is $\nu$-incoherent to the basis $\{\wa\}_{\alphab\in\mathbb{I}}$ and $\nu$-incoherent to its dual basis $\{\va\}_{\alphab\in\mathbb{I}}$; that is, there exists a constant $\nu\geq 1$ such that for all $\alphab = (i,j)\in\mathbb{I}$:
    \begin{equation}\label{eq: Incoherence equations}
        \left\Vert\Pu\wa\right\Vert_\fro\leq\sqrt{\frac{\nu r}{2n}},\quad \mathrm{and} \quad \left\Vert\Pu\va\right\Vert_\fro\leq\sqrt{\frac{\nu r}{2n}}.
    \end{equation}
    In addition to the above, we require that
        \begin{equation}\label{eq: Pt incoherence equations}
        \left\Vert\Pt\wa\right\Vert_\fro\leq\sqrt{\frac{\nu r}{2n}},\quad \mathrm{and} \quad\left\Vert\Pt\va\right\Vert_\fro\leq\sqrt{\frac{\nu r}{2n}}.
    \end{equation}
\end{asp}
 Notice that the two definitions in \eqref{eq: Incoherence equations} and \eqref{eq: Pt incoherence equations} are equivalent up to a small constant, as 
\begin{align*}
    \left\Vert\Pt\wa\right\Vert_\fro = \left\Vert \Pu\wa + \wa\Pu - \Pu\wa\Pu\right\Vert_\fro\leq3\left\Vert\Pu\wa\right\Vert_\fro,
\end{align*}
where the first inequality follows from the triangle inequality and the self-adjointness of $\Pu\wa$, and because
\begin{align*}
    \Vert\Pu\wa\Vert_\fro = \Vert\Pu\Pt\wa\Vert_\fro\leq\Vert\Pt\wa\Vert_\fro,
\end{align*}
where the equality follows from the definition of $\Pu$ and $\Pt$, and the inequality follows from Cauchy-Schwarz. As such, we pick a $\nu$ large enough such that the inequalities in \eqref{eq: Incoherence equations} and \eqref{eq: Pt incoherence equations} hold. We note that the constant difference in the condition stated above and in \Cref{sec:Incoherence_Discussion} is merely a matter of mathematical convenience. We also note that these incoherence conditions are similar to those seen in matrix completion with respect to the standard basis \cite{recht2011simpler}, as well as completion with respect to other bases \cite{gross2011recovering,tasissa2018exact}.

\begin{remark}\label{rem: incoherence remark}
    We want to note that $\nu$-incoherence with respect to $\wa$ in both $\eqref{eq: Incoherence equations}$ and $\eqref{eq: Pt incoherence equations}$ implies, at worst, $4\nu$-incoherence with respect to $\va$. As such, we choose a $\nu$ large enough so that both $\X$ is $\nu$-incoherent with respect to $\wa$ and $\va$. See \Cref{lem: Incoherence equivalence} for details.
\end{remark}

We provide one further assumption for this work. As we are typically interested in large $n$, assuming that $n\geq 3$ produces uniform results for several numerical bounds in the appendix, and is formally stated as an assumption.
\begin{asp}\label{asp: size}
    For the given ground truth rank-$r$ matrix $\X\in\Rnn$, we assume that $n\geq 3$.
\end{asp}

Throughout the remainder of this work, we will assume that our ground truth matrix $\X\in\mathbb{S}$ satisfies \Cref{asp: Incoherence assumption} with constant factor $\nu$. As in \cite{wei2020guarantees}, we identify a neighborhood in $\mfr$ around which any initial guess in this neighborhood converges linearly to the true solution with high probability using  \Cref{alg: M_omega descent}.

\subsection{Local Convergence Analysis}

The most critical property for a sampling operator to possess in matrix completion theory is the restricted isometry property, briefly discussed in \Cref{subsec: matrix completion}. This property roughly states that, when restricted to the local structure (or tangent space) around the true low-rank matrix, the partial observations preserve enough information to allow for faithful algorithmic recovery. We state this more formally with the following theorem: 
\begin{thm}[RIP of $\Mo$]\label{thm: RIP of Mo}
    Let $\X\in\S$ be the ground truth, rank-$r$, $\nu$-incoherent Gram matrix with tangent space $\T$ in $\mfr$. Let $\Omega$ be sampled from $\univ$ via a Bernoulli sampling process with parameter $p\geq \frac{128}{3}\beta\frac{\log{n}}{{n}}$ for any $\beta>3$. With probability at least $1-C'n^{3-\beta}$, where $C'$ is an absolute constant, we have that for some absolute numerical constant $C>0$ independent of $n$ or $\X$ that
    \[
    p^{-2}\Vert\Pt\Mo\Pt - p^2\Pt\Vert \leq 7\sqrt{\frac{\nu^2 r^2 \beta \log{n}}{pn}} +C\sqrt{\beta}\nu r\frac{\log{n}}{pn}.
    \]
  Furthermore, for any $\varepsilon_0$, if $p\geq \frac{C''\nu^2 r^2}{\varepsilon_0^2 }\frac{\beta\log{n}}{n}$ for some sufficiently large numerical constant $C''>0$, then
    \[
    p^{-2}\Vert\Pt\Mo\Pt - p^2\Pt\Vert \leq \varepsilon_0.
    \]
\end{thm}
\begin{proof}[Proof sketch]
    This proof works by decomposing $\Mo$ into diagonal and off-diagonal components. We recognize that estimating deviations from expectation in the off-diagonal terms in $\Vert\Pt\Mo\Pt-p^2\Pt\Vert$ is a second-order degenerate U-statistic---see\cite{lee2019u} for details on the theory of U-statistics---allowing the application of the decoupling inequality from \cite{de2012decoupling}, with the relevant theorem reproduced in \Cref{thm: Decoupling}. Once the random variables are decoupled, we can conditionally apply a non-commutative Bernstein inequality, reproduced in \Cref{thm:Bernstein}, and use the law of total probability to obtain a high-probability bound for the total deviation. The concentration of the diagonal terms in the full expression is equivalent to the concentration of $p\Vert\Pt\Fo\Pt-p\Pt\FI\Pt\Vert$, which is more simply analyzed through \Cref{thm:Bernstein} alone. See Appendix~\ref{subsubsec: proof of Mo RIP} for details. We also note that the last sample complexity statement supersedes the original sample complexity statement for sufficiently large $C''$.
\end{proof}

\begin{remark}
    We note that this result is given in terms of the Bernoulli sampling probability, rather than the more traditional number of samples with replacement seen in the matrix completion literature. To provide a more direct comparison, and remarking that $\Ebb[\vert\Omega\vert]=m$ and $p = \frac{m}{L}$, we have that for a sufficiently large constant $C>0$ that 
    \[
    m\geq C\frac{\nu^2 r^2}{\varepsilon_0^2} \beta n\log{n},
    \]
    gives $\varepsilon_0$-RIP of $\Mo$. We again note that, due to using the weaker \Cref{asp: Incoherence assumption} instead of the incoherence assumption in \cite{ghosh2024}, this is optimal up to constant factors and equivalent to the RIP established in \cite{ghosh2024}.
\end{remark}

Now that RIP is established, we can prove local convergence of \Cref{alg: M_omega descent}. This theorem describes a high-probability guarantee that \Cref{alg: M_omega descent} exhibits linear convergence in an attractive basin near the solution, provided that $\Mo$ exhibits RIP.

\begin{figure}[ht!]
\centering
\begin{tikzpicture}[
    node distance=1cm,
    every node/.style={rectangle, draw, minimum width=.8cm, minimum height=.8cm, text centered, font =\small},
    arrow/.style={thick,->,>=stealth}
    ]

    \node (RIP)[text width = 2.5cm]{\Cref{thm: RIP of Mo} $\Mo$ RIP};
    \node (waprop)[left of=RIP,xshift=-2cm,text width = 2.5cm]{Appendix~\ref{appendix: dual basis} Properties of the dual basis };
    \node (HWI)[below of=waprop,yshift=-1.3cm,text width = 2.5cm]{\Cref{thm: Decoupling} Decoupling of second-order degenerate U-statistics};
    \node (LRIP)[right of=RIP,xshift=3cm,text width=3cm]{\Cref{lem: Local RIP Mo} Local RIP of $\Mo$};
    \node (step)[below of =LRIP,yshift=-.3cm,text width = 2.5cm]{\Cref{lem: Mo Stepsize} Stepsize bound};
    \node (I3)[right of =LRIP,xshift=2.5cm]{$I_3$ bound};
    \node (I1)[below of=I3,yshift = -.3cm]{$I_1$ Bound};
    \node (I2)[below of =I1,yshift=-.3cm]{$I_2$ bound};
    \node (proj)[below of =step,yshift=-.3cm,text width = 2.8cm]{\Cref{lem: Projection Bounds} Projection bounds \cite{wei2020guarantees}};
    \node (1SHT)[below of= proj,text width = 2.3cm,yshift=-1cm]{\Cref{lem: Mo 1SHT initialization} One-step hard threshold bound};
    \node (convergence)[right of=I3, xshift = 2cm,text width = 2.3cm]{\Cref{thm: Mo Local convergence} Local Convergence of \Cref{alg: M_omega descent}};
    \node (Ber)[left of= 1SHT,xshift=-3cm,text width=2.5cm] {\Cref{thm:Bernstein} Non-commutative Bernstein Inequality};
    \node (1SHT cpxty)[below of=convergence,text width = 2.5cm,yshift = -3.6cm]{\Cref{thm: Mo Recovery Guarantee} One-step hard thresholding initialization guarantees};

    \draw [arrow] (Ber) -- (RIP);
    \draw [arrow] (waprop) -- (RIP);
    \draw [arrow] (RIP) -- (LRIP);
    \draw [arrow] (LRIP) -- (step);
    \draw [arrow] (step) -- (I1);
    \draw [arrow] (proj) -- (I2);
    \draw [arrow] (LRIP) -- (I3);
    \draw [arrow] (HWI) -- (RIP);

    \draw [arrow] (I1) -- (convergence);
    \draw [arrow] (I2) -- (convergence);
    \draw [arrow] (I3) -- (convergence);
    \draw [arrow] (1SHT) -- (1SHT cpxty);
    \draw [arrow] (convergence) -- (1SHT cpxty);
    \draw [arrow] (Ber) -- (1SHT);
\end{tikzpicture}

\caption{This diagram is a schematic of the overall proof of convergence for \Cref{alg: M_omega descent}. Arrows indicate how results depend on one another, and how they link together to form the overall proof of convergence. Not every exact dependency is shown in this figure for legibility purposes, instead focusing on the key pieces of the overall flow of the argument.}\label{fig: Mo Convergence proof schematic} 
\end{figure}
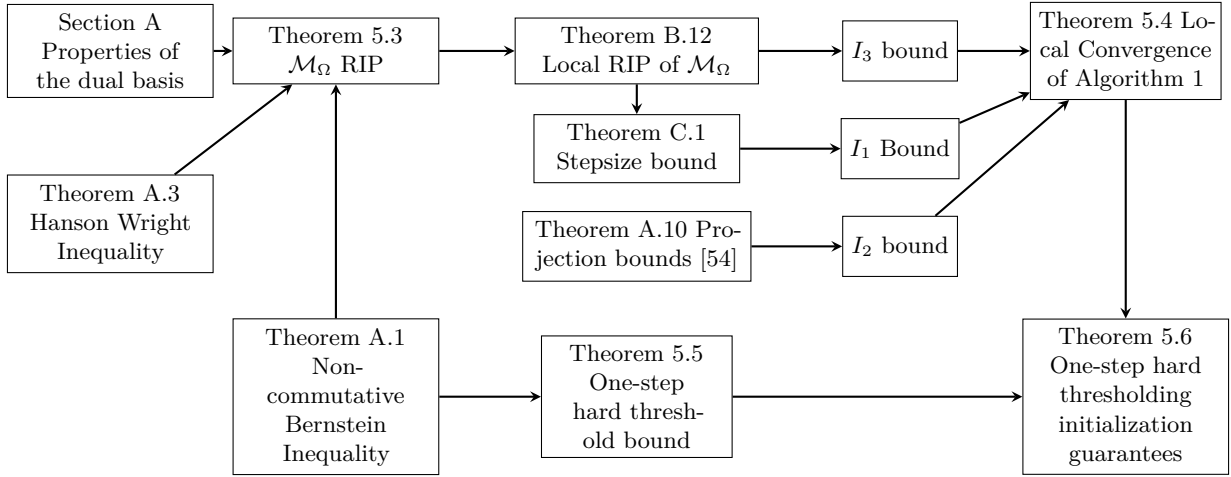

\begin{thm}[Local Convergence of \Cref{alg: M_omega descent}]\label{thm: Mo Local convergence}
 Let $\X\in\Rnn$ be the ground truth rank-$r$, $\nu$-incoherent matrix and let $\T$ be the tangent space of $\mfr$ at $\X$. Suppose that $p\geq C\frac{\beta \log{n} }{n}$ for some sufficiently large constant $C>0$ and any $\beta>3$. Further assume that
 \begin{align}
      p^{-2}\Vert\Pt\Mo\Pt-p^2\Pt\Vert&\leq\varepsilon_0,\label{eqn: Mo RIP assumption} \\
      \Vert\Mo\Vert&\leq  50p^{3/2}\sqrt{\beta n\log{n}},\label{eqn: M_o assumption}\\
        \Vert\Mo\Pt\Vert&\leq 25p^{3/2}\sqrt{{\nu r \beta\log{n}}},\label{eqn: MoPt assumption}\\
      \Vert\Mo\Ptl\Vert&\leq 100 p^{3/2}\sqrt{\beta n\log{n}}\frac{\Vert\X_l-\X\Vert_\fro}{\lambda_r(\X)}+25p^{3/2}\sqrt{{\nu r \beta\log{n}}},\label{eqn: MoPtl Assumption}\\
     \frac{\Vert\X_0-\X\Vert_\fro}{\lambda_r(\X)}&\leq\frac{\varepsilon_0 p^{1/2}}{50\left(\beta n\log{n}\right)^{1/4}},\label{eqn: Mo neighborhood assumption}
    \end{align}
    where $\varepsilon_0$ is a constant satisfying
        \[
        \delta=\frac{18\varepsilon_0}{1-4\varepsilon_0}<1.
        \]
        Then \Cref{alg: M_omega descent} converges linearly as the iterates satisfy 
        \[
        \Vert\X_{l} - \X\Vert_\fro \leq \delta^l\Vert\X_0-\X\Vert_\fro.
        \]
\end{thm}
\begin{proof}[Proof Sketch]
    We first note that all of the above assumptions, save for \eqref{eqn: Mo neighborhood assumption}, hold with high probability for $p\geq C\frac{\nu^2 r^2}{\varepsilon_0^2}{\frac{\beta\log{n}}{n}}$, where $C>0$ is an absolute constant. See Appendix~\ref{appendix: RIP} for explicit proof details on Assumptions \eqref{eqn: Mo RIP assumption}, \eqref{eqn: M_o assumption}, \eqref{eqn: MoPt assumption}, and \eqref{eqn: MoPtl Assumption}. \eqref{eqn: Mo RIP assumption}, proven to hold in \Cref{thm: RIP of Mo}, is the most stringent assumption and is what gives the theoretical sample complexity result. All other conditions hold with probability at least $1-C'n^{3-\beta}$ for absolute constants $C'>0$, $\beta>3$, in this sampling regime. 

    \noindent
    Once the requisite operator bounds are established by sufficiently high sample complexity, the proof of the theorem begins first by simple linear algebra, as we have
    \begin{align*}
        \Vert \X_{l+1} -\X\Vert_\fro &= \Vert\X_{l+1} -\W_l -\X + \W_l\Vert_\fro\\
        &\leq \Vert \X_{l+1} - \W_l\Vert_\fro + \Vert \X-\W_l\Vert_\fro\\
        &\leq 2\Vert \W_l-\X\Vert_\fro,
    \end{align*}
    where the last inequality follows from $\X_{l+1}$ being the best rank-$r$ approximation to $\W_l$ by Eckart-Young-Mirsky \cite{eckart1936approximation}. Next, plugging in $\W_l = \X_l + \alpha_l\Ptl\G_l$, we see that
    \begin{align*}
        \Vert \X_{l+1}-\X\Vert_\fro&\leq 2\left\Vert \X_l + \alpha_l\Ptl\G_l - \X\right\Vert_\fro\\
        &=2\Vert \X_l-\X -\alpha_l\Ptl\Mo(\X_l-\X)\Vert_\fro\\
        &\leq \underbrace{2\Vert (\Ptl - \alpha_l\Ptl\Mo\Ptl)(\X_l-\X)\Vert_\fro}_{I_1}\\
        &\quad+ \underbrace{2\Vert(\mathcal{I}-\Ptl)(\X_l-\X)\Vert_\fro}_{I_2}\\
        &\quad+ \underbrace{2\vert\alpha_l\vert\Vert\Ptl\Mo(\mathcal{I}-\Ptl)(\X_l-\X)\Vert}_{I_3}.
    \end{align*}
The remainder of the proof is in the bounding of $I_1$, $I_2$, and $I_3$. $I_1$ is proven by showing that in a neighborhood of the solution, defined by \eqref{eqn: Mo neighborhood assumption}, a local form of RIP for $\Mo$ holds if \eqref{eqn: Mo RIP assumption} is true; the proof of RIP is in Appendix \ref{subsubsec: proof of Mo RIP}, and the local form of RIP is proven in \Cref{lem: Local RIP Mo}. This proof leverages the assumptions made in \eqref{eqn: MoPt assumption}, and \eqref{eqn: MoPtl Assumption}, both of which are proven in \Cref{lem: PtMo Bounds}. $I_2$ follows from the neighborhood assumption of \eqref{eqn: Mo neighborhood assumption} in tandem with \Cref{lem: Projection Bounds}, and $I_3$ follows from bounds on the step size (seen in \Cref{lem: Mo Stepsize}), the assumption in \eqref{eqn: MoPt assumption}, and \Cref{lem: Projection Bounds}. We want to note that \Cref{alg: M_omega descent} converges deterministically provided the stated assumptions hold, and the assumptions themselves hold with high probability.
The technical details of the proof of convergence are deferred to the appendix; see Appendix~\ref{subsubsec: Mo Local convergence} for the full proof of \Cref{thm: Mo Local convergence}. See \Cref{fig: Mo Convergence proof schematic} for a diagram of the main dependencies for the convergence proof.
\end{proof}

\subsection{Initialization Results}
In this section, we outline our initialization guarantees for \Cref{alg: M_omega descent}. Given that the convergence of this algorithm is only local, initialization is important to consider in the context of sample complexity. The simplest initialization, a hard thresholding to $\mfr$ of the measured information, provides a reasonable starting point and is described in \Cref{alg: Initialization}, where $\mathcal{H}_r$ is as defined in \eqref{eqn: hard threshold operator}. The following sections describe how close a one-step hard-thresholding initialization will be to the ground truth for \Cref{alg: M_omega descent}. Following this, and in tandem with \Cref{thm: Mo Local convergence}, we show recovery guarantees for \Cref{alg: M_omega descent}.

\begin{algorithm}
\caption{1 Step Hard Thresholding Initialization}\label{alg: Initialization}
\begin{algorithmic}[1]
\STATE \textbf{Input:} $\Omega$, $\{\la \X,\wa\ra\}_{\alphab\in\Omega}$, $\rank(\X)=r$
    \STATE  $\Ro(\X) =\sum_{\alphab\in\Omega}\la\X,\wa\ra\va$
    \STATE $\X_0 = \frac{1}{p} \mathcal{H}_r(\Ro(\X))$
    \STATE\textbf{Output:} $\X_0$: 1 Step Hard Thresholding Initialization.
\end{algorithmic}
\end{algorithm}

\begin{lem}\label{lem: Mo 1SHT initialization}
    Under a Bernoulli sampling parameter $p\geq\frac{128\beta \log{n}}{3n}$, then with probability at least $1-2n^{1-\beta}$ we have for $\X_0 = p^{-1}\mathcal{H}
    _r(\Ro(\X))$ that
        \begin{equation*}
        \Vert \X_0-\X\Vert_\fro\leq\sqrt{2r}\Vert\X_0-\X\Vert \leq \sqrt{\frac{2\beta nr\log{n}}{3p}}\max_{\ai}\la\X,\wa\ra \leq \sqrt{\frac{ \beta \nu^2 r^3 \log(n)}{24pn}}\Vert \X\Vert.
    \end{equation*}
\end{lem}
\begin{proof}
    See Appendix~\ref{appendix: initialization}.
\end{proof}

\begin{thm}[Recovery Guarantee for \Cref{alg: M_omega descent}]\label{thm: Mo Recovery Guarantee}
    For $p\geq\max\left\{\frac{20\kappa\nu r^{3/2}}{\varepsilon_0}\frac{\beta^{3/4}\log^{3/4}(n)}{n^{1/4}},C\frac{\nu^2 r^2}{\varepsilon_0^2}{\frac{\beta\log{n}}{n}}\right\}$, where $\kappa$ is the condition number of $\X$, $\beta>3$, and with $\varepsilon_0<\frac{1}{22}$ for some sufficiently large constant $C>0$, then with probability $1-C'n^{3-\beta}$ where $C'$ is an absolute constant, \Cref{alg: M_omega descent} recovers the ground truth matrix $\X$ when initialized by \Cref{alg: Initialization}.
\end{thm}
\begin{proof}
    This result is a consequence of \Cref{lem: Mo 1SHT initialization} and the local neighborhood assumption in \eqref{eqn: Mo Alg Local Neighborhood}. We can see this by increasing the sample complexity $p$ to a sufficiently large value such that the initialization is smaller than the local neighborhood assumption. Furthermore, for a high enough sample complexity, all of the assumptions in \Cref{thm: Mo Local convergence} hold with high probability. The result then follows from the union bound.
\end{proof}

\begin{remark}
    For \Cref{alg: M_omega descent}, we use a Bernoulli sampling model with parameter $p$, while other matrix completion methodologies use a uniform at random with replacement model. To provide a more direct sample complexity comparison, let $m = \Ebb[\vert\Omega\vert]$ under a Bernoulli model. This implies that $p = \frac{m}{L}$. \Cref{thm: Mo Recovery Guarantee} therefore implies that, if
    \[
    m\geq \max\left\{\frac{20\nu r^{3/2} \kappa \beta^{3/4}}{\varepsilon_0} n^{7/4}\log^{3/4}(n),C\frac{\nu^2 r^2 \beta}{\varepsilon_0^2}n\log{n}\right\},
    \]
    for some sufficiently large constant $C>0$, \Cref{alg: M_omega descent} recovers $\X$.
\end{remark}

\begin{remark}
    We note here that a more delicate initialization through a resampling technique, such as the one in \cite{wei2020guarantees}, could likely reduce the sample complexity from $p\gtrsim \frac{\log^{3/4}{n}}{n^{1/4}}$ to $p\gtrsim \frac{\log{n}}{n}$. Further investigation of initialization has been omitted from this work due to space constraints, but it is an area of interest for future research.
\end{remark}

\section{Robustness Guarantees}\label{subsec: robustness}

In many applications, the distance matrix may be corrupted, and understanding the sources of this corruption is central to designing robust recovery algorithms \cite{khan2009diland,guo2023perturbation,biswas2006semidefinite_noisy,tasissa2025robustnodelocalizationrough,Kundu_2025,kundu2025dualbasisapproachstructured}. Broadly, there are two main causes. First, even if distance measurements are perfectly accurate, the underlying point configuration may itself be perturbed due to physical factors. For instance, sensors placed in dynamic environments, such as the ocean, may drift over time. In such cases, the observed distances correspond to a perturbed version of the true point set. Second, the points themselves may be fixed, but the distance measurements are noisy. This can arise from various sources: sensor imprecision, environmental interference, or limited measurement resolution. In practice, both types of corruption may occur simultaneously. However, in this paper, we focus on the first scenario: perturbations in the point configuration. This assumption simplifies the analysis, since the resulting distance matrix remains a valid Euclidean distance matrix, and avoids challenges associated with arbitrary noise patterns that could violate geometric consistency. We believe that this setting is relevant to settings where 
environmental drift is more dominant than measurement noise. Moreover, the developed technical analysis for this setting could potentially serve as a foundation for future extensions to more general noise models.

In this section, we will provide robustness results for \Cref{alg: M_omega descent}. To begin, we assume the following: for a given point matrix $\P$, we denote $\hat{\P} = \P+\N$, where $\N\in\real^{n\times r}$ is a random matrix. We denote $\hat{\X} = \hat{\P}\hat{\P}^\top$.
To demonstrate robustness to noise, we first show that, under bounded noise satisfying $\|\N\|_{\infty} \leq C$ for some constant $C$, the following holds:
\[
\|\X - \Xh\|_{\fro} \leq C\,\cO\!\big(\sqrt{n\|\X\|}\big) + C^2 n r,
\]
as established in \Cref{lem: X perturbation bound}.
This quantifies the amount that the drift noise can perturb the dataset by. Then, we show that for small enough noise, the incoherence of the perturbed matrix is not drastically changed, and is bounded by an $\cO(1)$ constant in \Cref{lem: coherence parameter after perturbation}. This result is important to demonstrate that the sample complexity required for recovery is not pushed out of a feasible range. With these results established, we can combine them together to show that for bounded noise, given a sufficiently large Bernoulli sample complexity $p$, that $\Xh$ is recovered with \Cref{alg: M_omega descent} with high probability and that the relative error between $\Xh$ and $\X$ is small. This result is formally stated in the following theorem:

\begin{thm}[Robustness Guarantee for \Cref{alg: M_omega descent}]\label{thm: Robustness guarantee for Mo}
Let $\P \in \real^{n \times r}$ and define $\Ph = \P + \N$, where $\N \in \real^{n \times r}$ is a random matrix with $\Ebb[N_{ij}] = 0$. Assume that the following quantity, which admits a signal-to-noise ratio (SNR) interpretation $\frac{\|\X\|^{1/2}}{\|\N\|_{\infty}}$,
satisfies
\[
\frac{\|\X\|^{1/2}}{\|\N\|_{\infty}} \geq \frac{32\sqrt{2}\,\beta \kappa n \log n}{3 \nu^{1/2}},
\]
for some $\beta > \max\left\{3, \frac{3r}{8 \log n}\right\}$, where $\kappa$ denotes the condition number of $\X$. 

Assume that the measured distances are of the form $\la \Xh, \wa \ra$ and are sampled according to a Bernoulli model with parameter $p$ satisfying
\[
p \geq C \nu^2 r^2 \frac{\beta \log n}{n},
\]
where $C > 3.92 \times 10^{4}$ is an absolute constant. Furthermore, suppose that \Cref{alg: M_omega descent} is initialized at a point $\X_0$ satisfying the conditions of \Cref{thm: Mo Local convergence} with $\varepsilon_0 < \frac{1}{22}$.

Then, with probability at least $1 - C' n^{3 - \beta}$, where $C'$ is an absolute constant, \Cref{alg: M_omega descent} recovers $\Xh$, and
\[
\|\X - \Xh\| \leq \left( \sqrt{\frac{\nu}{8n}}+ \frac{9 \nu r}{2048 \kappa \beta^2 n \log^2 n} \right) \lambda_r(\X).
\]

\end{thm}
\begin{proof}[Proof of \Cref{thm: Robustness guarantee for Mo}]
    This result follows first from Lemmas~\ref{lem: X perturbation bound} and \ref{lem: coherence parameter after perturbation}. From here, the sample complexity guarantee of \Cref{thm: RIP of Mo}, with sufficiently large constant $C$ as stated above, coupled with the high probability guarantees of Assumptions \eqref{eqn: Mo RIP assumption}-\eqref{eqn: MoPtl Assumption} in \Cref{thm: Mo Local convergence}, gives the desired result, provided an initialization satisfying~\eqref{eqn: Mo neighborhood assumption}. See \Cref{appendix: RIP} for details on the guarantees of Assumptions \eqref{eqn: Mo RIP assumption}-\eqref{eqn: MoPtl Assumption}.
\end{proof}

This result shows that recovery under noise depends intrinsically on the underlying geometry of the object. In particular, highly degenerate objects with large condition numbers can tolerate only low levels of noise, as captured by the ratio $\|\X\|^{1/2} / \|\N\|_{\infty}$ in \Cref{thm: Robustness guarantee for Mo}. 
Moreover, as this ratio increases, the allowable perturbation in the incoherence parameter also increases, which in turn leads to a higher sample complexity required for recovery. While this dependence may initially seem counterintuitive, the interplay between this noise-to-signal ratio and the sample complexity shows that the most recoverable objects under noise are those that are both well-conditioned and highly incoherent, as a large $\lambda_r(\X)$ facilitates stronger recovery. As a corollary to \Cref{thm: Robustness guarantee for Mo}, we obtain the following relative error bound:
\begin{cor}
Under the assumptions of \Cref{thm: Robustness guarantee for Mo}, the relative recovery error of \Cref{alg: M_omega descent} satisfies, with probability at least $1 - C' n^{3 - \beta}$ for some absolute constant $C'$,
\[
\frac{\|\X - \Xh\|_{\fro}}{\|\X\|_{\fro}} \leq  \sqrt{\frac{\nu}{4n}}+ \frac{9 \sqrt{2}\nu r}{2048 \kappa \beta^2 n \log^2 n} .
\]
\end{cor}
\begin{proof}
    This follows from \Cref{thm: Robustness guarantee for Mo} as 
    \[
    \frac{\Vert\X-\Xh\Vert_\fro}{\Vert\X\Vert_\fro} \leq \frac{\sqrt{2r}\Vert\X-\Xh\Vert}{\sqrt{r}\lambda_r(\X)} = \sqrt{2}\frac{\Vert\X-\Xh\Vert}{\lambda_r(\X)}.
    \]
\end{proof}

\section{Related Work}\label{sec: Related work}

\subsection{A Riemannian Approach to Matrix Completion}\label{subsec: RO and MC}

A notable non-convex approach to the low-rank matrix completion problem is to utilize prior knowledge regarding the rank of $\X$. Some techniques leverage the fact that the set of fixed-rank matrices forms a Riemannian manifold, turning the constrained low-rank matrix completion program into an unconstrained optimization task over a manifold. These methodologies lose convexity, however, and generally only local convergence guarantees can be established, done by proving the existence of attractive basins around solutions. Various retraction-based methodologies have been used with differing metrics and geometric structures \cite{vandereycken2012lowrank,mishra2013fixedrank,Boumal2015,dai2010,keshavan2010matrix,keshavan2010matrix2,wei2020guarantees}. The analysis conducted by \cite{wei2020guarantees} stands out for its interpretation of its first-order method as an iterative hard-thresholding algorithm with subspace projections and efficient numerical implementation. This implementation is done by reducing the hard thresholding step from a thin eigenvalue decomposition of an $n\times n$ matrix to a thin QR decomposition followed by a full eigenvalue decomposition of a far smaller $2r\times 2r$ matrix. The convergence analysis in this work builds on the analysis done in \cite{wei2020guarantees}, and as such, a brief exposition of their work is provided.

As mentioned in \Cref{sec: background}, the authors of \cite{wei2020guarantees} developed a gradient descent algorithm to solve the low-rank matrix completion problem by considering a first-order Riemannian gradient method on the manifold of fixed-rank matrices, $\mfr$. The objective function used in \cite{wei2020guarantees} can be seen in \eqref{eqn: P_omega qf objective}. The authors used a uniform sampling at random with replacement model for recovering a subset of the indices of the ground truth matrix. This is standard practice in existing matrix completion literature, as much of the analysis relies on concentration inequalities for sums of random matrices to get high probability guarantees. It follows that \eqref{eqn: P_omega qf objective} is not equivalent to $\Vert \Po(\X-\M)\Vert_\fro^2$ when indices in $\Omega$ repeat, as $\Po^2 \neq \Po$ when this occurs. This is distinct from \cite{vandereycken2012lowrank}, which minimized the Frobenius norm difference between the observed entries of the low-rank matrices to solve the problem. Additionally, \cite{vandereycken2012lowrank} demonstrates that the limit of their proposed algorithm agrees with the ground truth in the revealed entries when projected onto the tangent space of the ground truth. However, as the sampling operator has a non-trivial null space, noted in \cite{vandereycken2012lowrank}, this does not necessarily guarantee identification of the ground truth. In contrast, \cite{wei2020guarantees} establishes linear convergence to the ground truth solution in a local neighborhood of the ground truth, with high probability. After defining \eqref{eqn: P_omega qf objective}, \cite{wei2020guarantees} constructs a Riemannian gradient descent procedure similar to the retraction procedure described in \Cref{subsec: intro to RO} for its solution.

In addition to this approach, the work in \cite{wei2020guarantees} considered two initialization schemes. One is a simple one-step hard threshold onto $\mfr$, and is given by $\X_0 = \frac{n^2}{m}\mathcal{H}_r(\Po(\M))$. Additionally, a more delicate initialization can be considered by partitioning the set $\Omega$ into $S$ equally sized subsets, and performing one Riemannian gradient descent step for each subset. This Riemannian resampling initialization breaks the dependence on each iterate from the previous, and provides a more reliable initialization for large enough sample sizes.

\subsection{Decoupling and Matrix Completion}

This section highlights prior works using decoupling across various domains. In matrix completion, the seminal work of \cite{candes2009exact} employs decoupling to analyze nuclear norm minimization. Within their duality analysis framework in convex optimization, a critical step involves bounding the operator norm of a term in the adjoint of the sampling operator. This is achieved by decomposing a summation into diagonal and off-diagonal components, the latter of which is controlled via decoupling (see Section 6.2). Another use of decoupling in matrix completion includes \cite{li2024pairwise}, which addresses data missing not at random. There, decoupling is used to establish the statistical properties of a pseudolikelihood estimator, specifically, to prove the asymptotic recovery of the low-rank parameter matrix (Theorem 5.3). Beyond matrix completion, decoupling has been applied to low-rank covariance estimation \cite{wang2022low}, concentration inequalities for random tensors \cite{bamberger2022hanson}, K-means clustering for non-Euclidean data \cite{chen2021hanson}, and the invertibility of random submatrices \cite{CHRETIEN20121479}. 
While our work shares thematic similarities with these approaches, our proof of \Cref{thm: RIP of Mo} specifically reformulates the isometry property into a decoupling problem. Furthermore, we emphasize that a successful application of this technique requires deriving non-trivial, fine-grained linear algebraic bounds related to the dual basis construction.

\subsection{Euclidean Distance Matrix Completion Algorithms}

To solve the EDMC problem, various algorithms have been developed. Among them, one prominent family of algorithms is based on semi-definite programming (SDP), which leverages the connection between squared distance matrices and Gram matrices. To provide a concrete example of this approach, we briefly outline the method proposed in \cite{alfakih1999solving}. Consider the matrix $\bm{V} \in \mathbb{R}^{n\times (n-1)}$, whose columns form an orthonormal basis for the space $\{\bm{z} \in \mathbb{R}^{n} : \bm{z}^\top\bm{1} = 0\}$. The operator $\mathcal{K}$ is defined as:
\[
\mathcal{K}(\bm{X}) =  \text{diag}(\bm{X})\bm{1}^\top + \bm{1}\text{diag}(\bm{X})^\top - 2\bm{X}.
\]
This definition of the operator $\mathcal{K}(\bm{X})$ is equivalent to the mapping of the Gram matrix to the squared Euclidean distance matrix, as expressed in \eqref{eqn: X to D}. In \cite{alfakih1999solving}, the optimization program is based on the operator $\mathcal{K}_{\bm{V}}(\bm{X})$, which is defined as $\mathcal{K}_{\bm{V}}(\bm{X}) = \bm{V}\bm{X}\bm{V}^\top$. The optimization problem in \cite{alfakih1999solving} can then formulated as follows:
\begin{equation*}
\begin{split}
 \minimize_{\bm{X} \in \mathbb{R}^{(n-1) \times (n-1)},\,\,\bm{X} = \bm{X}^\top, \,\bm{X} \succeq \bm{0}} & \quad  \sum_{(i,j) \in \Omega} \left[(\mathcal{K}_{\bm{V}}(\bm{V}\bm{X}\bm{V}^\top))_{ij} - D_{ij}\right]^2.
\end{split}
\end{equation*}
We refer the reader to \cite{alfakih1999solving} for theoretical and numerical aspects of the above optimization program.
Given that standard SDP formulations can be computationally intensive, distributed and divide-and-conquer methods have also been explored. For additional SDP-based formulations of the EDMC problem and their applications to molecular conformation and sensor network localization, we refer the reader to \cite{biswas2006semidefinite, biswas2004semidefinite, biswas2008distributed, leung2010sdp, alipanahi2012protein, guo2023perturbation}.

In the context of protein structure determination, various algorithmic approaches to EDMC have been developed. One notable example is the EMBED algorithm \cite{havel1991evaluation,more1999distance,crippen1988distance}, which comprises three main steps \cite{havel1998distance}. The first step, known as bound smoothing, involves generating lower and upper bounds for all distances by extrapolating from the available limits of known distances. The second step is the embed step, where distances are sampled from these bounds to form a full distance matrix from which an initial estimate of the protein structure is obtained. The final step involves refining this initial structure by minimizing an energy function using non-convex optimization methods. Another approach to structure prediction is the discretizable molecular distance geometry framework, which can be formulated as a search in a discrete space followed by a Branch-and-Prune method \cite{lavor2012recent,lavor2012discretizable}. 

Another category of approaches to the EDMC problem involves initially estimating a smaller portion of the point cloud and then using this initial estimate to incrementally reconstruct the rest of of the structure. These methods are referred to as geometric build-up algorithms \cite{wu2007updated,dong2003geometric,sit2009geometric}. The algorithm proposed in \cite{hendrickson1995molecule} addresses the molecular conformation problem by adopting a divide-and-conquer strategy, where a sequence of smaller optimization problems is solved instead of solving a single global optimization problem. 

Next, we highlight algorithms that estimate the underlying points through non-convex optimization.  These utilize a combination of methods such as majorization, alternating projection, global continuation (transforming the optimization problem to a function with few local minimizers), and an asymmetric projected gradient descent scheme \cite{leeuw1977application,
more1997global,glunt1993molecular,fang2012euclidean,ghosh2024,li2025euclideandistancematrixcompletion}. One of particular interest is the iteratively re-weighted least squares (IRLS) methodology. This technique relies on computing smoothed \textit{log-det} objectives at each iterate of the continuous non-convex rank minimization problem, along with a least squares computation at each step. This algorithm relies on RIP of an operator related to $\Ro$, established for $\vert\Omega\vert\gtrsim \frac{\nu r}{\varepsilon_0^2}n\log{n}$ given a stronger incoherence assumption than used in this paper, and exhibits provable quadratic convergence in a local neighborhood around the solution provided RIP holds. No initialization guarantees are provided, however. 

Certain nonconvex EDG algorithms have been shown to have better performance when the problem is formulated in a dimension higher than the true rank of the underlying points \cite{fang2012euclidean, tasissa2018exact}. This overparameterization has previously been shown to enhance numerical performance in sensor network localization problems \cite{tang2023, lei2023}. However, to the best of our knowledge, theoretical guarantees for such overparameterization in EDG problems remain largely unexplored. 
A recent study \cite{criscitiello2025sensornetworklocalizationbenign} conducts a landscape analysis of a nonconvex optimization problem for classical MDS and identifies dimensional regimes that lead to benign optimization landscapes.

We note that the above discussion does not comprehensively cover all EDMC algorithms, and we refer readers to \cite{liberti2014euclidean,dokmanic2015euclidean} for a more detailed overview.

\subsubsection{Related Geometric Approaches to EDMC}
The main perspective taken in this paper is in line with the low-rank matrix completion approach, albeit not one that employs the trace heuristic seen in \cite{tasissa2018exact,biswas2006semidefinite,Javanmard_2012}. This work is more in line with non-convex approaches based on optimizing over a Riemannian manifold \cite{Nguyen2019,Parhizkar2013}, and extends the Riemannian approach of \cite{wei2020guarantees} to the EDMC basis case. 

A recent work in \cite{Li2024} adopts a similar approach to us and considers solving the EDMC problem through Riemannian methods as well. In this work, the authors use a Riemannian conjugate method paired with an inexact line search method to minimize the following s-stress objective function:
\begin{equation}\label{eqn: s-stress fn}
    \minimize_{\Y\in\R^{n\times d}} ~\frac{1}{2}\Vert \W\odot \Po(g(\Y\Y^\top)-\D_e)\Vert_\fro^2,
\end{equation}
where $g$ is the map defined by \eqref{eqn: X to D}, $\W$ is a weight matrix to model noisy entries, and $\odot$ is the Hadamard product, and $\Po$ is defined as in \eqref{eqn: Po definition}. The analysis in \cite{Li2024} centers around the minimization of the s-stress function in \eqref{eqn: s-stress fn} using a generalization of a Hager-Zhang line search method to a Riemannian quotient manifold. The main result in this work is that there exists an attractive basin for \eqref{eqn: s-stress fn} that, with high probability, gives linear convergence to the ground truth provided an initialization in the basin. This result requires a Bernoulli sample complexity $p>C\frac{(\nu r)^3\log{n}}{n}$, where $\nu$ is the incoherence of the ground truth matrix and $r$ is the rank. In contrast, our method also shows linear convergence in a local neighborhood and describes a strong initialization candidate for the noiseless EDMC recovery problem with provable high probability guarantees. We also provide robustness analysis for an EDMC problem perturbed by noise, and provide provable guarantees as well.

\section{Numerical Results} \label{sec: Numerics}
All the experiments were conducted in MATLAB. The code used for the following experiments can be found in the GitHub repository at \url{https://github.com/chandlersmith2/NonConvexEDMC}.

\subsection{Synthetic Data Experiment}
In this section, we test the proposed algorithm on synthetic data. Various two- and three-dimensional datasets were used, and are referred to in \Cref{tab: synthetic IPM means table} with their corresponding sizes. The goal of \Cref{alg: M_omega descent} is to recover the full set of points $\Pb$ up to orthogonal transformation by sampling the entries above the diagonal of $\D$ uniformly with replacement, with a total of $\gamma L$ entries chosen for $\gamma\in[0,1]$. The algorithm reconstructs the Gram matrix $\X_\mathrm{rec} = \P\P^\top$, from which $\P$ can be recovered using \eqref{eqn: MDS Equation}. The comparison referenced in \Cref{tab: synthetic IPM means table} is the relative error between the recovered matrix $\X_\mathrm{rec}$ and the ground truth matrix $\X$ in Frobenius norm. Each run was terminated at either 1000 iterations or or when the relative Frobenius norm difference between successive iterates fell below $10^{-5}$.. This experiment was initialized using the one-step-hard-thresholding method described in \Cref{alg: Initialization}.

\begin{table}[ht!]
\centering
\caption{\textsc{Relative recovery error $\left\Vert \X-\X_\mathrm{rec}\right\Vert_\fro/\left\Vert \X\right\Vert_\fro$ between the recovered Gram matrix and the true Gram matrix averaged over 25 trials using \Cref{alg: M_omega descent} when initialized by \Cref{alg: Initialization}. \vspace{0.05in}}} \label{tab: synthetic IPM means table}
\begin{tabular}{l|c c c c c c} 
\toprule
\diagbox[width=12em]{Dataset}{$\gamma$}&
   {10\%}&  {7\%}& {5\%} & {3\%} & {2\%} & {1\%}  \\
\midrule
 Sphere (3D, $n=1002$)& $3.38\times 10^{-7}$ &  $4.61\times 10^{-7}$ & $6.12 \times 10^{-7}$ & $1.48\times 10^{-6}$ & $8.40\times 10^{-3}$ & $6.81\times 10^{-1}$ \\
 
Cow (3D, $n=2601$)& $4.41 \times 10^{-7}$& $5.24\times 10^{-7}$&  $6.04\times 10^{-7}$&  $9.14\times 10^{-7}$&  $2.47\times 10^{-4}$&  $5.71\times 10^{-3}$\\

Swiss Roll (3D, $n=2048$)& $3.85\times 10^{-7}$&  $4.70\times 10^{-7}$  &  $5.81\times 10^{-7}$ &  $9.47\times 10^{-7}$  &  $1.56\times 10^{-6}$  &  $6.40\times 10^{-2}$\\
\bottomrule
 \end{tabular}
\end{table}

\begin{figure}
    \centering
    \includegraphics[width=.8\linewidth]{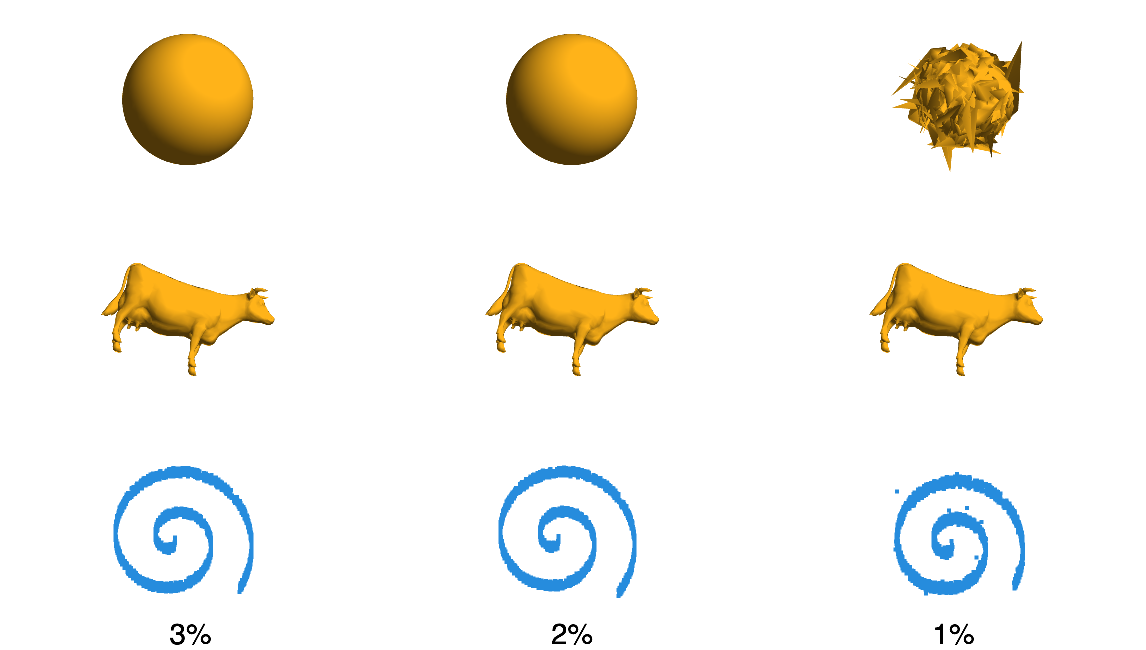}
    \caption{Reconstruction of the synthetic datasets referenced in \Cref{tab: synthetic IPM means table}. From left to right, the Bernoulli parameter is $0.03$, $0.02$, and $0.01$.}
    \label{fig: synthetic data figure}
\end{figure}

We note that the recovery completely fails for the sphere at $1\%$ sampling, while recovery is partially successful for the other two datasets. This is because the other datasets are larger while maintaining the same rank, allowing for better scaling in the low sampling regime. In \Cref{fig: synthetic data figure}, we show an image of the reconstruction of the figures described in \Cref{tab: synthetic IPM means table}.

\subsection{Comparison to Existing Methods}
\begin{figure*}[!t]
\centering
\subfloat[DBRE]{\includegraphics[width=2in]{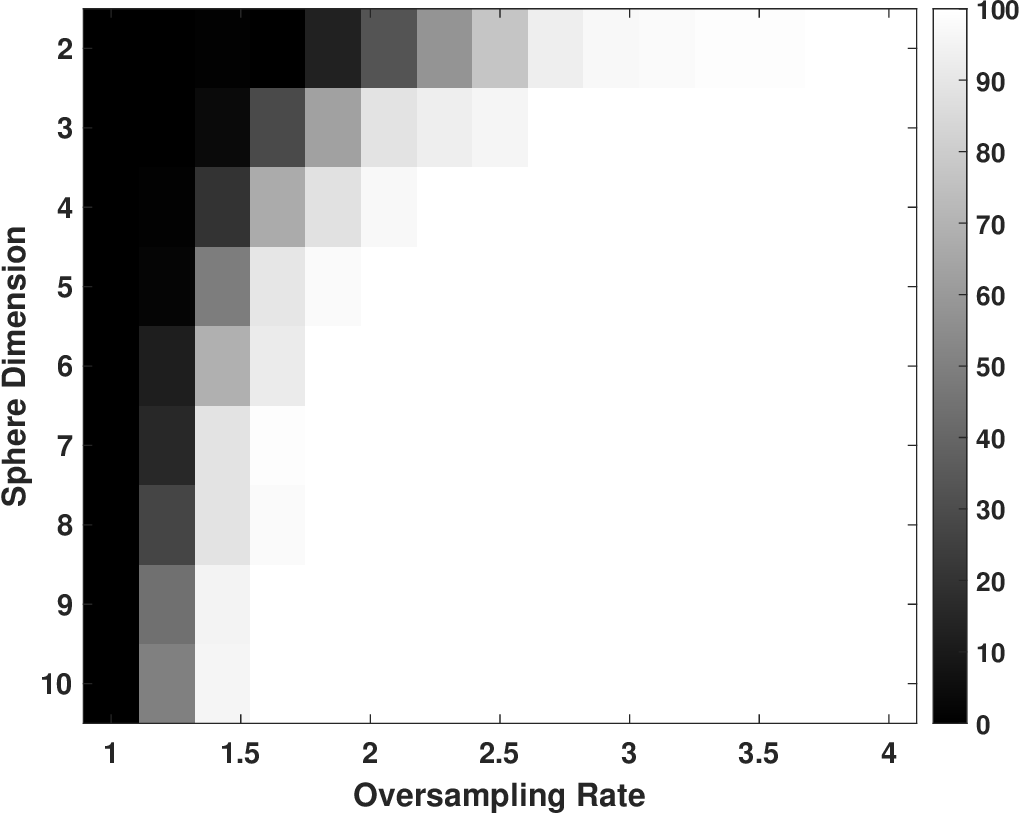}
\label{fig_first_case}}
\hfil
\subfloat[MatrixIRLS-EDMC]{\includegraphics[width=2in]{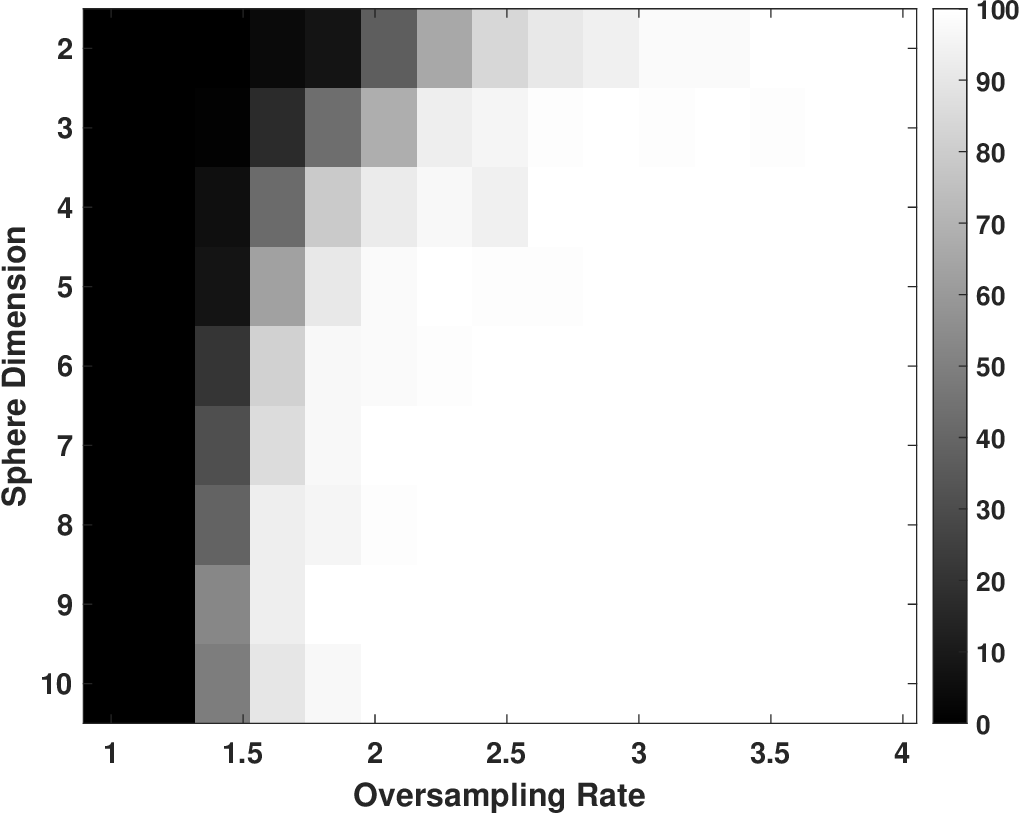}
\label{fig_second_case}}
\hfil
\subfloat[Augmented Lagrangian]{\includegraphics[width=2in]{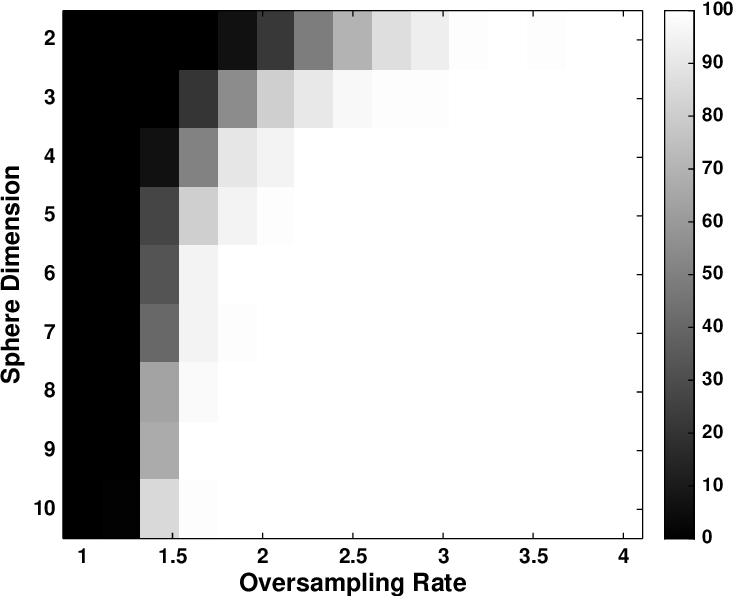}}
\caption{Oversampling ratio $\rho$ versus dimension $r$ for $100$ points on the uniform distribution on $S^{r-1}$. Each parameter was tested 100 times.}
\label{fig: oversampling on spheres}
\end{figure*}

We provide additional experiments to compare the efficacy of our algorithm DBRE to another provably convergent non-convex EDMC algorithm\cite{ghosh2024}, which we refer to as MatrixIRLS-EDMC. For the first experiment, let $r\in[2,10]$, and consider $n=100$ points sampled from $\mathrm{Unif}(S^{r-1})$, the uniform distribution on the sphere embedded in $r$ dimensions. As the number of degrees of freedom in a rank-$r$ $n\times n$ matrix is $nr - \frac{r(r-1)}{2}$, define the oversampling ratio $\rho$ as 
\[
\rho = \frac{pL}{nr - \frac{r(r-1)}{2}},
\]
as $\Ebb[\vert\Omega\vert] = pL$ for Bernoulli random sampling with parameter $p$. In \Cref{fig: oversampling on spheres}, we compare the oversampling ratio versus the dimension of the sphere in a transition plot. Black indicates complete failure, classified as a relative Gram matrix error larger than $10^{-3}$, and white indicates success. Each of these squares was run for 100 trials using DBRE, MatrixIRLS-EDMC, and the Augmented Lagrangian technique from \cite{tasissa2018exact}. DBRE was initialized with 10 iterations of the Augmented Lagrangian algorithm, and MatrixIRLS-EDMC was initialized using a least-squares methodology described in \cite{ghosh2024}.
As these experiments indicate, performance between the Augmented Lagrangian, MatrixIRLS for EDMC, and \Cref{alg: M_omega descent} is comparable, with \Cref{alg: M_omega descent} performing slightly better overall, but most noticeably in the higher rank regime. 

Next, we compare the performance of DBRE, the Augmented Lagrangian, and the MatrixIRLS methods on the synthetic datasets seen in \Cref{fig: synthetic data figure}. We initialize both DBRE and the Augmented Lagrangian method with 20 iterations of the Augmented Lagrangian method, and we use the default settings for the MatrixIRLS algorithm. DBRE was run for 1000 iterations, or until a $10^{-8}$ relative difference between iterates was achieved.

\begin{table}[ht!]
\centering
\caption{\textsc{Relative recovery error $\left\Vert \X-\X_\mathrm{rec}\right\Vert_\fro/\left\Vert \X\right\Vert_\fro$ between the recovered Gram matrix and the true Gram matrix averaged over 25 trials comparing DBRE, MatrixIRLS, and the Augmented Lagrangian approach}. \vspace{0.05in}}
\label{tab: synthetic IPM means all algs table}
\begin{tabular}{l|c c c c c c} 
\toprule
\diagbox[width=7em]{Dataset}{$\gamma$}&
{10\%}& {7\%}& {5\%} & {3\%} & {2\%} & {1\%} \\
\midrule
\multicolumn{7}{c}{\textbf{DBRE}} \\
\midrule
Sphere&
$6.38\times 10^{-8}$ &
$5.72\times 10^{-8}$ &
$7.27\times 10^{-8}$ &
$1.58\times 10^{-7}$ &
$5.08\times 10^{-7}$ &
$6.47\times 10^{-1}$ \\

Cow&
$4.92\times 10^{-8}$ &
$5.67\times 10^{-8}$ &
$6.51\times 10^{-8}$ &
$9.93\times 10^{-8}$ &
$1.55\times 10^{-7}$ &
$6.28\times 10^{-7}$ \\

Swiss Roll&
$4.24\times 10^{-8}$ &
$4.88\times 10^{-8}$ &
$6.45\times 10^{-8}$ &
$1.01\times 10^{-7}$ &
$1.71\times 10^{-7}$ &
$2.29\times 10^{-6}$ \\
\midrule
\multicolumn{7}{c}{\textbf{Augmented Lagrangian}} \\
\midrule
Sphere&
$2.33\times 10^{-7}$ &
$2.13\times 10^{-7}$ &
$3.53\times 10^{-7}$ &
$9.29\times 10^{-7}$ &
$4.44\times 10^{-6}$ &
$1.05\times 10^{-1}$ \\

Cow&
$1.97\times 10^{-7}$ &
$3.03\times 10^{-7}$ &
$3.86\times 10^{-7}$ &
$5.69\times 10^{-7}$ &
$7.51\times 10^{-7}$ &
$2.30\times 10^{-7}$ \\

Swiss Roll&
$1.22\times 10^{-6}$ &
$3.27\times 10^{-6}$ &
$2.42\times 10^{-6}$ &
$4.16\times 10^{-6}$ &
$7.35\times 10^{-6}$ &
$2.45\times 10^{-6}$ \\
\midrule
\multicolumn{7}{c}{\textbf{MatrixIRLS}} \\
\midrule
Sphere &
$1.30\times 10^{-7}$ &
$1.33\times 10^{-7}$ &
$1.36\times 10^{-7}$ &
$3.81\times 10^{-6}$ &
$3.76\times 10^{-2}$ &
$6.09\times 10^{-1}$ \\

Cow &
$4.48\times 10^{-7}$ &
$5.33\times 10^{-7}$ &
$6.12\times 10^{-6}$ &
$3.75\times 10^{-3}$ &
$3.53\times 10^{-2}$ &
$1.12\times 10^{-1}$ \\

Swiss Roll &
$1.10\times 10^{-7}$ &
$1.23\times 10^{-7}$ &
$1.28\times 10^{-7}$ &
$5.14\times 10^{-7}$ &
$2.78\times 10^{-3}$ &
$4.78\times 10^{-1}$\\
\bottomrule
 \end{tabular}
\end{table}

We note that for an initialization at the same point, both DBRE and the Augmented Lagrangian algorithm perform nearly identically in their recovery. MatrixIRLS is outperformed by both of the other algorithms in this experiment, but this is likely a result of comparing the default initialization to the Augmented Lagrangian initialization. Previous work \cite{ghosh2024} shows similar performance between the Augmented Lagrangian and MatrixIRLS. Both \Cref{tab: synthetic IPM means all algs table} and \Cref{fig: oversampling on spheres} both indicate the comparable to state-of-the-art performance of DBRE.

\subsection{Experiments on Noisy Distance Measurements}
\begin{figure}
    \centering
    \includegraphics[width=0.6\linewidth]{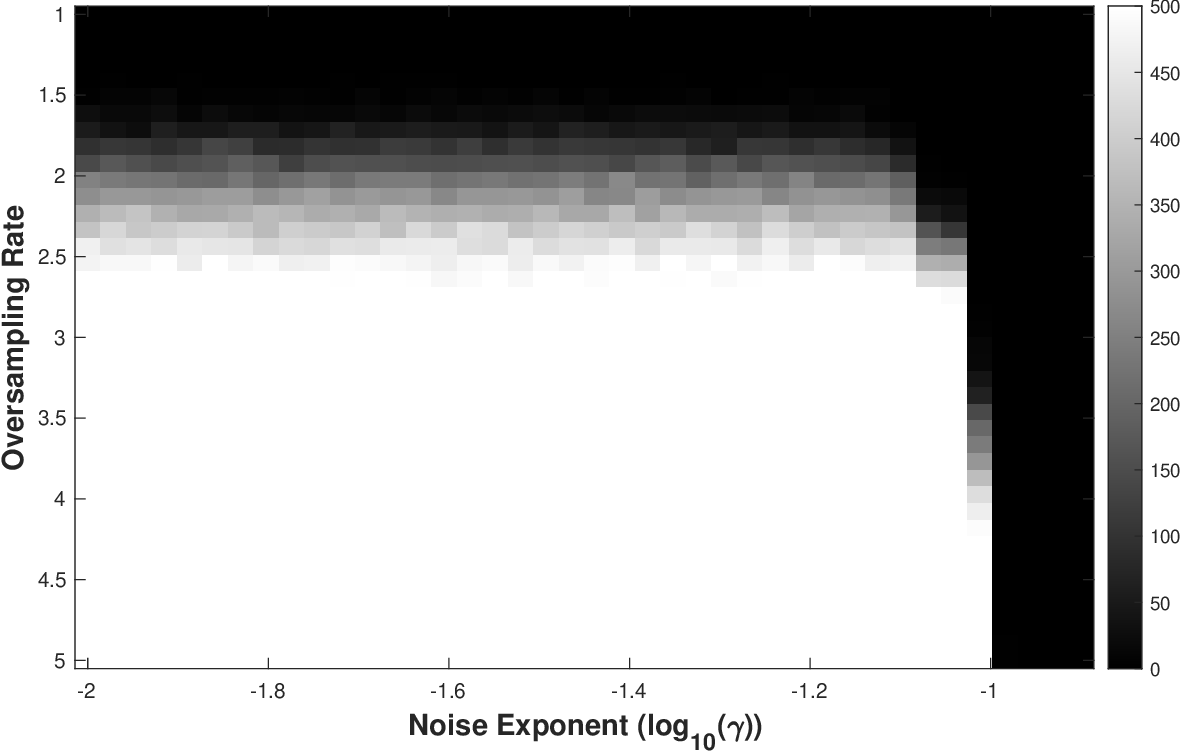}
    \caption{Oversampling ratio $\rho$ versus noise level $10^{\gamma}$ for $100$ points drawn i.i.d. from $\mathrm{Unif}(S^2)$. Each parameter was tested 500 times.}
    \label{fig: noise experiment}
\end{figure}

Finally, we also ran an experiment with noise following the model in \Cref{subsec: robustness} using \Cref{alg: M_omega descent}. Let $\{\p_i\}_{i=1}^{100}\sim\mathrm{Unif}(S^{2})$ be drawn i.i.d. and where $\P = [\p_1\cdots\p_{100}]^\top\in\real^{100\times 3}$. We perturb $\P$ with a bounded, centered noise matrix $\N\in\real^{100\times 3}$ with $\Vert\N\Vert_\infty \leq 10^{\gamma}$ for $\gamma\in[-2,-1]$. Similar to the previous experiment, we set the oversampling ratio $\rho\in[1,5]$. We set the success threshold at $10^{-2}$ relative difference, a relaxed value from previous experiments due to the addition of noise. \Cref{fig: noise experiment} shows the results over 500 trials.

\Cref{fig: noise experiment} indicates that recovery up to tolerance predominantly gets worse with added noise, although there is some clear dependence on the size of the noise impacting the reconstruction of the ground truth. This is most likely due to an increase in the incoherence of the dataset, requiring higher sample complexities to reconstruct. However, the noise level is still the dominant factor, and after a large enough noise value, reconstruction up to a certain tolerance is no longer viable.

We also present recovery results for the spherical point configuration shown in \Cref{fig: synthetic data figure} perturbed by Gaussian noise with varying standard deviations. For this experiment, the sampling rate $\gamma$ varies from $0.05$ to $0.01$, and the standard deviation of the Gaussian perturbations $\sigma$ ranges from $10^{-6}$ to $1$. The relative error between the Gram matrix recovered using \Cref{alg: M_omega descent} for 500 iterations, denoted $\hat{\X}_\mathrm{rec}$, and the unperturbed Gram matrix $\X$ is presented in \Cref{tab: noisy data}, averaged over 10 trials. For this experiment, \Cref{alg: M_omega descent} was initialized with \Cref{alg: Initialization}.

\begin{table}[ht!]
\centering
\caption{\textsc{Relative recovery error $\left\Vert \X-\hat{\X}_\mathrm{rec}\right\Vert_\fro/\left\Vert \X\right\Vert_\fro$ between the recovered perturbed Gram matrix and the true Gram matrix averaged over 10 trials using \Cref{alg: M_omega descent} for the Sphere dataset referenced in \Cref{tab: synthetic IPM means table} and \Cref{fig: synthetic data figure}. $\sigma$ is the standard deviation of the perturbing noise, and $\gamma$ is the sampling rate. \vspace{0.05in}}} \label{tab: noisy data}
\sisetup{
    scientific-notation = true,
    retain-zero-exponent = true,
    exponent-product = \times,
    round-mode = figures,
    round-precision = 2
}
\begin{tabularx}{\linewidth}{c | *{6}{>{\centering\arraybackslash}X}}
\toprule
\diagbox[width=3em]{$\gamma$}{$\sigma$} 
& $1=10^{0}$  & $10^{-1}$ & $10^{-2}$&  $10^{-3}$ & $10^{-4}$  & $10^{-5}$ \\
\midrule
0.05 
& $9.68\times 10^{-1}$
& $2.39\times 10^{-1}$ 
& $2.45\times 10^{-3}$
& $2.46\times 10^{-3}$
& $2.44\times 100^{-4}$
& $2.42\times 10^{-5}$ \\
0.04 
& $9.67\times  10^{-1}$ 
& $2.40 \times 10^{-1}$ 
& $2.45 \times 10^{-2}$ 
& $2.49\times 10^{-3}$
& $2.43 \times 10^{-4}$
& $2.46\times 10^{-5}$ \\
0.03 
& $9.68\times 10^{-1}$ 
& $2.39\times 10^{-1}$ 
& $2.44\times 10^{-2}$ 
& $2.53\times 10^{-3}$ 
& $1.03\times 10^{-2}$ 
& $2.44 \times 10^{-5}$  \\
0.02 
& $9.70 \times 10^{-1}$ 
& $3.02\times 10^{-1}$ 
& $6.21\times 10^{-2}$ 
& $1.62 \times 10^{-2}$ 
& $4.02 \times 10^{-2}$ 
& $3.83 \times 10^{-2}$  \\
0.01 
& $1.08$
& $1.35$ 
& $1.39$ 
& $1.34$ 
& $1.36$ 
& $1.38$ \\
\bottomrule
\end{tabularx}
\end{table}

The experimental results in \Cref{tab: noisy data} also indicate that the magnitude of the drift noise is the dominant factor in the reconstructability of the point cloud. While there are sample complexity effects, which are noticeable for $\gamma = 0.03$, the dominant effect is the noise level itself, further confirming the conclusion drawn from \Cref{fig: noise experiment}.

\section{Conclusion and Future Work}\label{sec: Conclusion}

In this work, we proposed a novel Riemannian gradient descent approach for solving the EDMC problem using a matrix completion approach on the manifold of rank-$r$ matrices in \Cref{alg: M_omega descent}. In a local neighborhood, we proved that \Cref{alg: M_omega descent} exhibits linear convergence with high probability. To the authors' knowledge, this is the first work to provide both initialization guarantees and robustness guarantees for a non-convex approach to the EDMC problem. The convergence analysis of \Cref{alg: M_omega descent} was predicated on a statistical understanding of coupled terms in a random operator, and required novel analysis to the matrix completion literature. For our method, we provided numerical results to underline its efficacy, and \Cref{alg: M_omega descent} performs comparably to other state-of-the-art non-convex methods. Additionally, we provided robustness analysis and corresponding convergence guarantees. Finally, we provided a novel interpretation of incoherence in the EDMC setting, highlighting potential areas for development of non-uniform sampling methods in this field. This is a primary avenue of future interest, as improving the sample complexity through geometrically-optimal sampling schemes would represent a noteworthy development in the EDMC literature.

\section*{Acknowledgments}
Abiy Tasissa and Chandler Smith acknowledge partial support from the National Science Foundation through grant DMS-2208392. 
HanQin Cai acknowledges partial support from the National Science Foundation through grant DMS-2304489. HanQin Cai and Abiy Tasissa appreciate insightful discussions with Prof. Rongjie Lai and Prof. Jian-Feng Cai at the initial stages of this work. The authors would also like to thank Dr. Christopher Criscitiello and Dr. Andrew D. McRae for their invaluable feedback on the manuscript; in particular, we thank them for identifying a sub-optimal bound in a prior version. 

\bibliographystyle{IEEEtran} 
\bibliography{IEEEabrv,manuscript_main_nonconvex_edg}

\appendices
\section{Properties of the Dual Bases and Concentration Inequalities}\label{appendix: dual basis}
This section of the appendix develops technical results concerning the specific dual bases $\{\wa\}_{\ai}$ and $\{\va\}_{\ai}$, as well as related constructions derived from them. These results are used to establish various auxiliary lemmas throughout the paper and play a particularly important role in the proof of \Cref{thm: RIP of Mo}. We also prove that $\mathbb{E}[\Mo] = \Id$. In addition, for completeness, we collect several probabilistic tools that are used repeatedly in the analysis. These include  vector, and noncommutative (operator) Bernstein inequalities, the Hanson--Wright inequality, a decoupling theorem for degenerate U-statistics, and the Davis--Kahan Theorem.

\begin{thm}[Operator Bernstein Inequality\cite{recht2011simpler,tropp2012user}]\label{thm:Bernstein}
Let $\X_1,\dots,\X_m\in\Rnn$ be independent, mean-zero random matrices. Assume that $\Vert \X_i\Vert\leq c$ almost surely for all $i$, and define $\sigma^2\geq \max\left\{\left\Vert\sum_{i=1}^m\mathbb{E}\left(\X_i\X_i^\top\right)\right\Vert,\left\Vert\sum_{i=1}^m\mathbb{E}\left(\X_i^\top\X_i\right)\right\Vert\right\}$. Then, for all $t>0$,
\[
\bb{P}\left(\left\Vert\sum_{i=1}^m \X_i\right\Vert>t\right)\leq 2n \exp\left(-\frac{t^2/2}{\sigma^2 + ct/3}\right).
\]
In particular,
\begin{equation}\label{eqn:Bernstein}
    \bb{P}\left(\left\Vert\sum_{i=1}^m \X_i\right\Vert>t\right)\leq 2n \exp\left(-\frac{3t^2}{8\sigma^2}\right) \quad \text{for } t \leq \frac{\sigma^2}{c},
\end{equation}
and
\begin{equation*}
    \bb{P}\left(\left\Vert\sum_{i=1}^m \X_i\right\Vert>t\right)\leq 2n \exp\left(-\frac{3t}{8 c}\right) \quad \text{for } t \geq \frac{\sigma^2}{c}.
\end{equation*}

\end{thm}

\begin{thm}[Vector Bernstein Inequality \cite{gross2011recovering}]\label{thm: Vector Bernstein}
Let $\X_1,\dots,\X_n \in \mathbb{R}^d$ be independent, mean-zero random vectors. Suppose that there exist constants $R > 0$ and $\sigma^2 > 0$ such that $\|\X_i\|_2 \leq R$ almost surely for all $i$, and $\sum_{i=1}^n \mathbb{E}\big[\|\X_i\|_2^2\big] \leq \sigma^2$. Then, for all $t > 0$ satisfying $t\leq \frac{\sigma^2}{R}$,
\begin{align*}
\mathbb{P}\left(\left\|\sum_{i=1}^n \X_i\right\|_2\geq t\right)\leq \exp\left(\frac{-t^2}{4\sigma^2}\right).
\end{align*}
\end{thm}

\begin{thm}[Hanson-Wright Inequality~\cite{rudelson2013hansonwrightinequalitysubgaussianconcentration}]\label{thm: Hanson Wright}
    Let $\Y = (Y_1 \cdots Y_n)\in\real^n$ be a random vector with $Y_i$ independent and $\mathbb{E}[Y_i] = 0$, and $\Vert Y_i\Vert_{\psi_2}\leq K$ for some $K\geq 0$, where $\Vert\cdot\Vert_{\psi_2}$ is the sub-Gaussian norm. Additionally, let $\A\in\Rnn$. Then, for every $t\geq 0$,
    \[
    \mathbb{P}\left[\left\vert \Y^\top\A\Y - \mathbb{E}\left[\Y^\top\A\Y\right]\right\vert> t\right] \leq 2\exp\left(-c\min\left\{\frac{t^2}{K^4 \Vert\A\Vert_\fro^2},\frac{t}{K^2\Vert\A\Vert}\right\}\right).
    \]
\end{thm}

\begin{thm}[Matrix--valued Decoupling Inequality {\cite[Theorem 3.4.1]{de2012decoupling}}]\label{thm: Decoupling}
Let $\X_1,\dots\,\X_n$ be independent random variables in $\real^m$, and let $\X_1',\dots,\X_n'$ be an independent copy  defined on the same probability space. Let $\{h_{ij}\}_{1 \leq i \ne j \leq n}$ be a collection of measurable kernels
$h_{ij} : \real^m \times\real^m\to \real^{m\times m}$ satisfying
\[
\mathbb{E}\big[h_{ij}(\X_i,\X_j) \mid \X_i\big] = \bm{0},
\qquad
\mathbb{E}\big[h_{ij}(\X_i,\X_j) \mid \X_j\big] = \bm{0},
\]
for all $i \ne j$. Define the coupled and decoupled chaos sums as
$$
S := \sum_{i \ne j} h_{ij}(\X_i,\X_j),
\qquad
S' := \sum_{i \ne j} h_{ij}(\X_i,\X_j').
$$
Then there exists a universal constant $C > 0$ such that, for all $t > 0$,
$$
\mathbb{P}\big( \|S\| \geq t \big)
\;\leq\;
C \, \mathbb{P}\big( \|S'\| \geq t / C \big).
$$
The constant $C$ is universal and independent of the kernels $h_{ij}$, the variables $\{\X_i\}_{i=1}^n$, and the dimensions $n$ and $m$.
.
\end{thm}

\noindent \textbf{Example}: To provide a concrete example of how to use the above result, we consider Bernoulli random variables $\{\xi_{\alphab}\}_{\alphab \in I}$ that are $1$ with probability $p$ and $0$ otherwise. To establish the main theorem in this manuscript, one quantity of central interest (see \ref{lem: Concentration of Fo}) is 
\[
(\xi_{\alphab} - p)(\xi_{\betab} - p)\, H^{\alphab\betab} \U_{\alphab} \V_{\betab}^\top.
\]
This term is a second-order matrix chaos and can be written in the form required by the decoupling result  To see this, for $\alphab \neq \betab$ define the matrix-valued kernel
\[
h_{\alphab\betab}(\xi_{\alphab},\xi_{\betab})
:= (\xi_{\alphab} - p)(\xi_{\betab} - p)\, H^{\alphab\betab} \U_{\alphab} \V_{\betab}^\top.
\]
By independence of the Bernoulli variables and centering, we have
\[
\mathbb{E}\big[h_{\alphab\betab}(\xi_{\alphab},\xi_{\betab})\mid \xi_{\alphab}\big] = \bm{0},
\qquad
\mathbb{E}\big[h_{\alphab\betab}(\xi_{\alphab},\xi_{\betab})\mid \xi_{\betab}\big] = \bm{0},
\]
so the kernel satisfies the conditions required by Theorem~A.4. Therefore, introducing an independent copy $\{\xi'_{\betab}\}$, the decoupled chaos
$
\sum_{\alpha\neq\beta} (\xi_{\alphab} - p)(\xi'_{\betab} - p)\, H^{\alphab\betab} \U_{\alphab} \V_{\betab}^\top
$ satisfies
\[
\mathbb{P}\!\left(
\left\|
\sum_{\alphab\neq\betab} (\xi_{\alphab} - p)(\xi_{\betab} - p)\, H^{\alphab\betab} \U_{\alphab} \V_{\betab}^\top
\right\| \geq t
\right)
\leq
C\,
\mathbb{P}\!\left(
\left\|
\sum_{\alphab\neq\betab} (\xi_{\alphab} - p)(\xi'_{\betab} - p)\, H^{\alphab\betab} \U_{\alphab} \V_{\betab}^\top
\right\| \geq t/C
\right).
\]

\begin{thm}[Davis--Kahan $\sin \Theta$ theorem (rank-$r$ case) ~\cite{daviskahan1970},~\cite{yu2015useful}]\label{thm: Davis Kahan}]
Let $\X, \hat \X \in \Rnn$ be symmetric positive semidefinite matrices of rank $r$. Let $\V, \hat \V \in \mathbb{R}^{n \times r}$ have orthonormal columns spanning the top-$r$ eigenspaces of $\X$ and $\hat \X$, respectively. Then
\[
\|\sin\Theta(\V, \hat \V)\| \leq \frac{\|\X - \hat \X\|}{\lambda_r(\X)}.
\]
\end{thm}

One result that will be used throughout this work is a technique for constructing eigenvalue bounds through a vectorization technique. This result is as follows.
\begin{lem}[[Quadratic Form and Gram Matrix]\label{lem: vectorization}
Let $\{\Z_k\}_{k=1}^m$ be a basis for a subspace $\mathbb{V}\subset\Rnn$, and let $\G = [\la \Z_i,\Z_j\ra]\in\real^{m\times m}$. Then for any $\Y\in\Rnn$,
\[
\max_{\Vert\Y\Vert_\fro=1}\sum_{k=1}^m \la\Y,\Z_k\ra^2 = \lambda_{\mathrm{max}}(\G).
\]
\end{lem}

\begin{proof}
Using $\la \Y,\Z_k\ra = \Vec(\Y)^\top \Vec(\Z_k)$, we obtain
\[
\sum_{k=1}^m \la\Y,\Z_k\ra^2
= \Vec(\Y)^\top\Big(\sum_{k=1}^m \Vec(\Z_k)\Vec(\Z_k)^\top\Big)\Vec(\Y)
= \Vec(\Y)^\top \Z_{\mathbb{V}}\Z_{\mathbb{V}}^\top \Vec(\Y),
\]
where $\Z_{\mathbb{V}}=[\Vec(\Z_1)\ \cdots\ \Vec(\Z_m)]$. Hence
\[
\max_{\Vert\Y\Vert_\fro=1}\sum_{k=1}^m \la\Y,\Z_k\ra^2
= \lambda_{\mathrm{max}}(\Z_{\mathbb{V}}\Z_{\mathbb{V}}^\top)
= \lambda_{\mathrm{max}}(\Z_{\mathbb{V}}^\top \Z_{\mathbb{V}})
= \lambda_{\mathrm{max}}(\G),
\]
where we used the variational characterization of the largest eigenvalue and the identity $\lambda_{\mathrm{max}}(\A\A^\top)=\lambda_{\mathrm{max}}(\A^\top\A)$.
\end{proof}

\begin{lem}[Gram Matrix Bound under Incoherence]\label{lem: Bound for largest eigval of H tilde}
Let $\tilde{\H} = [\langle \Pu\wa,\Pu\wb\rangle]\in\mathbb{R}^{L\times L}$, where $U$ is the row/column space of the true solution $\X = \U\D\U^\top$ of rank $r$, and $\Pu$ is the projection onto $U$. Then
\[
\lambda_\mathrm{max}(\tilde{\H}) \leq \nu r.
\]
\end{lem}

\begin{proof}
By incoherence,
\[
|\langle \Pu\wa,\Pu\wb\rangle|
\leq \|\Pu\wa\|_\fro \|\Pu\wb\|_\fro
\leq \frac{\nu r}{2n}.
\]
Moreover, since $\Pu = \U\U^\top$, for $\alphab \cap \betab = \emptyset$,
\[
\langle \Pu\wa,\Pu\wb\rangle
= \tr(\wa\Pu\Pu\wb)
= \tr(\wb\wa\Pu)
= 0,
\]
as $\wa\wb = \bm{0}$ when $\alphab \cap \betab = \emptyset$. Thus each row of $\tilde{\H}$ has at most $2n-3$ nonzero entries. A Gershgorin argument combined with the entrywise bound yields $
\lambda_\mathrm{max}(\tilde{\H}) \leq \nu r$. 
\end{proof}

\begin{lem}[Tangent Projection Decomposition and Energy Bound]\label{lem: Pu wa inner product identity}
For any symmetric $\X\in\Rnn$ and any $\wa\in\{\wb\}_{\betab\in\mathbb{I}}$,
\[
\langle \Pt\X,\wa\rangle
= \langle \X,\Pu\wa\rangle + \langle \X\Pup,\Pu\wa\rangle.
\]
Additionally, if $\|\X\|_\fro=1$,
\[
\sum_{\ai}\la\X,\Pt\wa\ra^2 \leq 4\lambda_{\mathrm{max}}(\tilde{\H}) = 4\nu r.
\]
\end{lem}

\begin{proof}
Using symmetry and cyclicity of the trace,
$
\langle \X\Pu,\wa\rangle = \langle \Pu\X,\wa\rangle$. It follows then that
\begin{align*}
\langle \Pt\X,\wa\rangle
&= \langle \Pu\X+\X\Pu-\Pu\X\Pu,\wa\rangle \\
&= 2\langle \Pu\X,\wa\rangle - \langle \Pu\X\Pu,\wa\rangle \\
&= \langle \Pu\X,\wa\rangle + \langle \Pu\X\Pup,\wa\rangle \\
&= \langle \X,\Pu\wa\rangle + \langle \X\Pup,\Pu\wa\rangle.
\end{align*}
For the second claim, by $(a+b)^2 \leq 2a^2+2b^2$,
\[
\sum_{\ai}\la\X,\Pt\wa\ra^2
\leq 2\sum_{\ai}\la\X,\Pu\wa\ra^2
+ 2\sum_{\ai}\la\X\Pup,\Pu\wa\ra^2.
\]
Since $\Pup$ is an orthogonal projection, the two terms are bounded identically, giving
\[
\sum_{\ai}\la\X,\Pt\wa\ra^2
\leq 4 \max_{\Vert\X\Vert_\fro=1}\sum_{\ai}\la\X,\Pu\wa\ra^2.
\]
The conclusion follows from Lemmas~\ref{lem: vectorization} and \ref{lem: Bound for largest eigval of H tilde}.
\end{proof}

\begin{lem}[Gram Matrix Spectrum and Dual Basis Structure {\cite{lichtenberg2023dual}}]\label{lem: H and H^-1 eigvals}
Let $\H\in\real^{L \times L}$ be the Gram matrix associated with $\{\wa\}_{\alphab\in\bb{I}}$, defined by $H_{\alphab\betab} = \la \wa,\wb\ra$, and let $\H^{-1} = [H^{\alphab\betab}]$. Then
\[
\lambda_\mathrm{max}(\H) = 2n, 
\qquad 
\lambda_\mathrm{max}(\H^{-1}) = \frac{1}{2}.
\]
Additionally,
\[
H^{\alphab\betab} =
\begin{cases}
\frac{1}{n^2}, & \alphab\cap\betab = \emptyset,\\[4pt]
-\frac{1}{2n} + \frac{1}{n^2}, & \alphab\cap\betab \neq \emptyset,\ \alphab\neq\betab,\\[4pt]
\frac{1}{2}\!\left(1-\frac{2}{n}+\frac{2}{n^2}\right), & \alphab=\betab,
\end{cases}
\]
and
\[
\Vert \wa\Vert = 2, 
\qquad 
\Vert \va\Vert = \frac{1}{2}.
\]
\end{lem}

\begin{lem}[Sum of Squares of Dual Basis]\label{lem:Form of va^2}
Let $\{\va\}_{\alphab\in\bb{I}}$ be the dual basis to $\{\wa\}_{\alphab\in\bb{I}}$. Then
\[
\sum_{\alphab\in\bb{I}} \va^2 
= \frac{n^2-2n+2}{4n}\,\J.
\]
\end{lem}

\begin{proof}
Recall that $\va = -\frac{1}{2}(\a\b^\top + \b\a^\top)$, where $\a = \e_i - \frac{1}{n}\one$ and $\b = \e_j - \frac{1}{n}\one$,
for $\alphab=(i,j)$. Then
\[
4\va^2 
= \a\b^\top\a\b^\top + \a\b^\top\b\a^\top 
+ \b\a^\top\a\b^\top + \b\a^\top\b\a^\top.
\]
Using $\a^\top\a = \b^\top\b = \frac{n-1}{n}$ and $\a^\top\b = -\frac{1}{n}$, a direct expansion yields
\[
4\va^2
= \frac{n-1}{n}(\e_{ii}+\e_{jj})
- \frac{1}{n}(\e_{ij}+\e_{ji})
+ \frac{2-n}{n^2}\big(\e_i\one^\top+\one\e_i^\top+\e_j\one^\top+\one\e_j^\top\big)
+ \frac{2(n-2)}{n^3}\oot.
\]
Summing over $\alphab\in\bb{I}$ and collecting terms gives
\[
\sum_{\alphab\in\bb{I}} 4\va^2
= \frac{n^2-2n+2}{n}\I 
- \frac{n^2-2n+2}{n^2}\oot.
\]
Since $\J = \I - \frac{1}{n}\oot$, the result follows.
\end{proof}

\begin{lem}[Subspace and Tangent Projection Perturbation Bounds {\cite{Weirecovery2016,wei2020guarantees}}]\label{lem: Projection Bounds}
Let $\X_l = \U_l \D_l \U_l^\top$ and $\X = \U\D\U^\top$ be rank-$r$ matrices, with corresponding tangent spaces $\T_l$ and $\T$ of $\mfr$. Then
\begin{gather*}
\|\U_l\U_l^\top - \U\U^\top\|
\leq \frac{\|\X_l-\X\|_\fro}{\sigma_r(\X)},
\qquad
\|\U_l\U_l^\top - \U\U^\top\|_\fro
\leq \frac{\sqrt{2}\,\|\X_l-\X\|_\fro}{\sigma_r(\X)},\\
\|(\mathcal{I} - \Ptl)\X\|_\fro
\leq \frac{\|\X_l-\X\|_\fro^2}{\sigma_r(\X)},
\qquad
\|\Ptl-\Pt\|
\leq \frac{2\,\|\X_l-\X\|_\fro}{\sigma_r(\X)}.
\end{gather*}
\end{lem}

\subsection{Operator Construction}

We begin by examining the expectation of the sampling operator $\RRo$ and showing that it is not isotropic under Bernoulli sampling. This motivates a simple reweighting that restores isotropy in expectation.

\begin{lem}[Expectation and Debiasing of the Sampling Operator]\label{lem: Expectation of Mo}
Let $\RRo(\cdot) = \sum_{\alphab,\betab\in\Omega}\wa\la\va,\vb\ra\la\cdot,\wb\ra$, and let $\Mo$ be as defined in \eqref{eqn: Mo definition}. Then
\[
\Ebb[\RRo] \neq p^2 \Id,
\qquad
\Ebb[\Mo] = p^2 \Id.
\]
\end{lem}

\begin{proof}
We first decompose $\RRo$ into diagonal and off-diagonal contributions:
\begin{align*}
\RRo(\cdot)
&= \sum_{\alphab\in\Omega}\la\cdot,\wa\ra \la\va,\va\ra \wa
+ \sum_{\alphab\in\Omega}\wa \sum_{\substack{\betab\in\Omega\\ \betab\neq\alphab}}
\la\va,\vb\ra \la\cdot,\wb\ra.
\end{align*}
This decomposition highlights an asymmetry under Bernoulli sampling: diagonal terms depend on a single index, while off-diagonal terms depend on pairs. In related operator constructions (e.g., $\mathcal{R}_{\univ}^\ast \mathcal{R}_{\univ}$ in \cite{Smith2023}),  numerical behavior can suggest an isotropic normalization. Here, however, the Bernoulli model introduces a mismatch in scaling that must be accounted for. Assuming each index is included independently with probability $p$, we obtain
\begin{align*}
\Ebb[\RRo(\cdot)]
&= p \sum_{\alphab\in\univ} \la\cdot,\wa\ra \la\va,\va\ra \wa
+ p^2 \sum_{\alphab\in\univ}\wa \sum_{\substack{\betab\in\univ\\ \betab\neq\alphab}}
\la\va,\vb\ra \la\cdot,\wb\ra.
\end{align*}
The diagonal terms scale with $p$, while the off-diagonal terms scale with $p^2$, so $\Ebb[\RRo]$ is not proportional to the identity.

To correct this imbalance, define the reweighted operator
\[
\Mo(\cdot)
= \sum_{\alphab,\betab\in\Omega} C_{\alphab\betab}\,
\wa \la\va,\vb\ra \la\cdot,\wb\ra,
\qquad
C_{\alphab\betab} =
\begin{cases}
1, & \alphab\neq\betab,\\
p, & \alphab=\betab.
\end{cases}
\]
Then
\begin{align*}
\Mo(\cdot)
&= p \sum_{\alphab\in\Omega}\la\cdot,\wa\ra \la\va,\va\ra \wa
+ \sum_{\alphab\in\Omega}\wa \sum_{\substack{\betab\in\Omega\\ \betab\neq\alphab}}
\la\va,\vb\ra \la\cdot,\wb\ra.
\end{align*}
Taking expectations yields
\begin{align*}
\Ebb[\Mo(\cdot)]
&= p^2 \sum_{\alphab\in\univ}\la\cdot,\wa\ra \la\va,\va\ra \wa
+ p^2 \sum_{\alphab\in\univ}\wa \sum_{\substack{\betab\in\univ\\ \betab\neq\alphab}}
\la\va,\vb\ra \la\cdot,\wb\ra \\
&= p^2 \sum_{\alphab,\betab\in\univ}
\wa \la\va,\vb\ra \la\cdot,\wb\ra
= p^2 \mathcal{R}_{\univ}^\ast \mathcal{R}_{\univ}
= p^2 \Id,
\end{align*}
since $\mathcal{R}_{\univ} = \Id$.
\end{proof}

\noindent
The preceding expectation calculation highlights that the diagonal and off-diagonal contributions scale differently under Bernoulli sampling. This motivates isolating the off-diagonal structure of $\H^{-1}$ via the following “hollow” matrix.

\begin{definition}[Off-diagonal part of $\H^{-1}$]\label{def: Hob}
Let $\H^{-1}\in\real^{L\times L}$ be defined by $(\H^{-1})_{\alphab\betab} = H^{\alphab\betab} = \langle \va,\vb\rangle$. We define the following matrix $\Hob$ as follows:
\begin{equation}
(\Hob)_{\alphab\betab}
=
\begin{cases}
H^{\alphab\betab}, & \alphab\neq\betab,\\
0, & \alphab=\betab.
\end{cases}
\label{eqn: Hob defn}
\end{equation}
\end{definition}

\begin{lem}[Norm Bounds for $\Hob$]\label{lem: Hob eigval bound}
For $\Hob$ as in \eqref{eqn: Hob defn},
\[
\|\Hob\|\leq 2,
\qquad
\|\Hob\|_{2\to\infty}^2 \leq \frac{2}{n},
\]
where $\|\A\|_{2\to\infty}$ denotes the maximum $\ell_2$ norm of the rows of $\A$.
\end{lem}

\begin{proof}
By Gershgorin’s theorem,
\[
\|\Hob\|
\leq \max_{\alphab}\sum_{\betab\neq\alphab} |H^{\alphab\betab}|.
\]
Using the bounds on $H^{\alphab\betab}$ from \Cref{lem: H and H^-1 eigvals}, each row contains at most $2n-3$ entries of size at most $\frac{1}{2n}+\frac{1}{n^2}$, and the remaining entries are of order $\frac{1}{n^2}$. Hence
\[
\|\Hob\|
\leq \max_{\alphab}\Big((2n-3)\Big(\tfrac{1}{2n}+\tfrac{1}{n^2}\Big) + \tfrac{L}{n^2}\Big)
\leq 2.
\]

For the second bound, note that all rows of $\Hob$ share the same multiset of entries. Thus, for any $\alphab$,
\[
\|\Hob\|_{2\to\infty}^2
= \sum_{\betab\neq\alphab} (H^{\alphab\betab})^2
\leq (2n-3)\Big(\tfrac{1}{2n}+\tfrac{1}{n^2}\Big)^2 + \frac{L}{n^4}
\leq \frac{2}{n}.
\]
\end{proof}

\section{Restricted Isometry Results}\label{appendix: RIP}
As RIP and its variants are critical to the analysis of \Cref{alg: M_omega descent} in this paper, this section is dedicated to the proofs of RIP and similar results.
\subsection{Supporting Concentration Results}
\noindent
The following results are various concentration identities used in the computation of $\Vert\Pt\Mo\Pt-p^2\Pt\Vert$, $\Vert\Pt\Mo\Vert$, and $\Vert\Mo\Vert$.

\begin{lem}[Concentration of $\hat{\Fo}$]\label{lem: Concentration of Fo}
Let $\hat{\Fo}$ denote either $\Fo$ or $\Pt\Fo\Pt$, and define $\hat{\FI}$ and $\wah$ accordingly (i.e., $\hat{\FI}=\FI$ or $\Pt\FI\Pt$, and $\wah=\wa$ or $\Pt\wa$). Let $\hat{\H} = [\la \wah,\wbh\ra]\in\real^{L\times L}$ be the associated Gram matrix. Suppose $\X$ is a rank-$r$ $\nu$-incoherent matrix with tangent space $\mathbb{T}$ on $\mfr$. Then for any $\beta>1$, if $p\geq \frac{8\beta\log n}{3n}$, we have with probability at least $1-2n^{1-\beta}$ that
\[
\|\hat{\Fo} - p\hat{\FI}\|
\leq
\sqrt{\frac{8\beta p\,\big(\max_{\alphab}\|\wah\|_\fro^2\big)\lambda_{\max}(\hat{\H})\log n}{3}}.
\]
\end{lem}

\begin{proof}
We express $\hat{\Fo}$ as a sum of independent random operators. Let $\{\xi_{\alphab}\}_{\alphab\in\bb{I}}$ be i.i.d. Bernoulli($p$) variables. Then
\[
\hat{\Fo}(\cdot) = \sum_{\alphab} \xi_{\alphab}\,\la \cdot,\wah\ra\,\wah,
\qquad
\Ebb[\hat{\Fo}] = p\hat{\FI}.
\]
Define centered operators
\[
\bm{S}_{\alphab} = (\xi_{\alphab}-p)\,\la \cdot,\wah\ra\,\wah,
\quad\text{so that}\quad
\sum_{\alphab}\bm{S}_{\alphab} = \hat{\Fo}-p\hat{\FI}.
\]
Each $\bm{S}_{\alphab}$ satisfies
\[
\|\bm{S}_{\alphab}\|
\leq \|\wah\|_\fro^2
\leq \max_{\alphab}\|\wah\|_\fro^2 =: c.
\]
For the variance term,
\begin{align*}
\left\Vert\Ebb\left[\sum_{\ai}\bm{S}_{\alphab}^2 \right]\right\Vert & =\left\Vert \Ebb\left[\sum_{\ai}\left(\xi_{\alphab}-p\right)^2\la\cdot,\wah\ra\la\wah,\wah\ra\wah\right]\right\Vert\\
        & = \left\Vert\Ebb\left[\sum_{\ai}(\xi_{\alphab}-2\xi_{\alphab} p + p^2)\la\cdot,\wah\ra\la\wah,\wah\ra\wah\right]\right\Vert\\
        & \leq p(1-p)\left(\max_{\alphab}\Vert\wah\Vert_\fro^2\right)\left\Vert\sum_{\ai}\la\cdot,\wah\ra\wah\right\Vert\\
        &\leq p\left(\max_{\alphab}\Vert\wah\Vert_\fro^2\right)\lambda_{\max}(\hat{\H}),
\end{align*}
where the final inequality follows from \Cref{lem: vectorization}. Applying matrix Bernstein with $t = \sqrt{\frac{8}{3}\sigma^2\beta\log n}$ and using $p\geq \frac{8\beta\log n}{3n}$ to ensure $t\leq \sigma^2/c$, we obtain
\[
\mathbb{P}\Big(\|\hat{\Fo}-p\hat{\FI}\|\geq t\Big)
\leq 2n\exp(-\beta\log n)
= 2n^{1-\beta},
\]
which yields the claim.
\end{proof}

\begin{lem}[Concentration of $\Fo\Pt$]\label{lem: concentration of FoPt}
Let $\X$ be a rank-$r$, $\nu$-incoherent matrix with tangent space $\T$ on $\mfr$. For any $\beta>1$, if $p\geq \frac{4}{3}\frac{\beta\log n}{n}$, then with probability at least $1-2n^{1-\beta}$,
\[
\|\Fo\Pt - p\FI\Pt\|
\leq \sqrt{\frac{128\,p\nu r\,\beta\log n}{3}}.
\]
\end{lem}

\begin{proof}
We represent $\Fo\Pt - p\FI\Pt$ as a sum of independent, mean-zero operators as follows:
\[
\Fo\Pt - p\FI\Pt
= \sum_{\alphab} (\xi_{\alphab}-p)\,\la\cdot,\Pt\wa\ra\,\wa
=: \sum_{\alphab}\J_{\alphab},
\]
Each summand satisfies
\[
\|\J_{\alphab}\|
\leq \|\Pt\wa\|_\fro \|\wa\|_\fro
\leq \sqrt{\tfrac{2\nu r}{n}}
=: c,
\]
using incoherence and $\|\wa\|_\fro=2$. For the variance terms,
\[
\left\|\sum_{\alphab}\Ebb[\J_{\alphab}\J_{\alphab}^\ast]\right\| = \left\|\sum_{\ai} p(1-p)\la \cdot,\wa\ra\la\Pt\wa,\Pt\wa\ra\wa\right\|\leq p\,\frac{\nu r}{2n}\,\lambda_{\max}(\H)
= p\nu r,
\]
and
\[
\left\|\sum_{\alphab}\Ebb[\J_{\alphab}^\ast\J_{\alphab}]\right\| = \left\|\sum_{\ai} p(1-p) \la\cdot,\Pt\wa\ra\la\wa,\wa\ra\Pt\wa\right\|
\leq 4p\,\lambda_{\max}(\tilde{\H})
\leq 16p\nu r.
\]
Thus $\sigma^2 \leq 16p\nu r$. Applying matrix Bernstein with 
\[
t = \sqrt{\tfrac{128}{3}p\nu r\,\beta\log n},
\]
and noting $t\leq \sigma^2/c$ under $p\geq \frac{4}{3}\frac{\beta\log n}{n}$, yields
\[
\mathbb{P}\!\left(\|\Fo\Pt - p\FI\Pt\|\geq t\right)
\leq 2n^{1-\beta}.
\]
\end{proof}

\begin{lem}\label{lem: Bound on So H So - E() term}
Let $\{\A_{\alphab}\}_{\ai},\{\B_{\alphab}\}_{\ai}\subset\Rnn$ be any two collections of $L$ matrices. Then
    \begin{align*}
    \left\Vert\sum_{\alphab\neq\betab}(\xa\xi_{\betab}-p^2)\A_{\alphab} H^{\alphab\betab}\langle\cdot,\B_{\betab}\rangle\right\Vert&\leq \left\Vert\sum_{\alphab\neq\betab}(\xa-p)(\xi_{\betab}-p)H^{\alphab\betab}\A_{\alphab}\langle\cdot,\B_{\betab}\rangle\right\Vert \\
        & \quad+ {\left\Vert\sum_{\alphab\neq\betab}(\xa-p)p\A_{\alphab} H^{\alphab\betab}\langle\cdot,\B_{\betab}\rangle\right\Vert}\\
        &\quad + {\left\Vert\sum_{\alphab\neq\betab}(\xi_{\betab}-p)p\A_{\alphab} H^{\alphab\betab}\langle\cdot,\B_{\betab}\rangle \right\Vert}.
    \end{align*}
\end{lem}
\begin{proof}
    This result follows from the fact that 
    \[
    \xa\xi_{\betab}-p^2 = (\xa-p)(\xi_{\betab}-p) + p(\xa-p) + p(\xi_{\betab}-p).
    \]
    As such, the expression becomes
    \begin{align*}
        \left\Vert\sum_{\alphab\neq\betab}(\xa\xi_{\betab}-p^2)\A_{\alphab} H^{\alphab\betab}\langle\cdot,\B_{\betab}\rangle\right\Vert&= \left\Vert \sum_{\alphab\neq\betab}\left((\xa-p)(\xi_{\betab}-p) + p(\xa-p) + p(\xi_{\betab}-p)\right)\A_{\alphab} H^{\alphab\betab}\langle\cdot,\B_{\betab}\rangle\right\Vert\\
        & = \bigg\Vert\sum_{\alphab\neq\betab}(\xa-p)(\xi_{\betab}-p)\A_{\alphab} H^{\alphab\betab}\langle\cdot,\B_{\betab}\rangle + \sum_{\alphab\neq\betab}p(\xa-p)\A_{\alphab} H^{\alphab\betab}\langle\cdot,\B_{\betab}\rangle \\
        & \quad +\sum_{\alphab\neq\betab}p(\xi_{\betab}-p)\A_{\alphab} H^{\alphab\betab}\langle\cdot,\B_{\betab}\rangle\bigg\Vert.
    \end{align*}
    The result then follows from the triangle inequality.
\end{proof}

\begin{lem}\label{eqn: Matrix B bound}
    Let $\X$ be a rank-$r$, $\nu$-incoherent ground truth matrix with tangent space $\T$. Let $\Va = \vec(\Pt\wa)$ or $\vec(\wa)$ for all $\ai$, and let $\G = [\Va^\top\Vb]\in\real^{L\times L}$ be the Gram matrix of $\{\Va\}_{\ai}$. Defining $\ba =\sum_{\betab\neq\alphab}H^{\alphab\betab}(\xi'_{\betab}-p)\Vb$, and letting $\B = [\b_1 \cdots \b_L]$, we have that $\B = \sum_{\alphab}(\xi_{\alphab}'-p)\Va\ha^T$, where $\ha\in\real^L$ is the $\alphab$-th column of $\Hob$, and that with probability at least $1-2n^{2-\beta}$ that if $p\geq\frac{16\beta\log{n}}{3n}$, for $\Va = \vec(\Pt\wa)$,
    \[
    \Vert\B\Vert \leq \sqrt{p\frac{32\nu r\beta\log{n}}{3n}},
    \]
and for $\Va = \vec(\wa)$,
\[
\Vert\B\Vert \leq \sqrt{p\frac{128\beta\log{n}}{3}}.
\]
\end{lem}
\begin{proof}
The first statement follows from the rules of matrix multiplication, and can be seen explicitly by defining the matrix $\W = [(\xi_1'-p)\V_1 \cdots (\xi_L'-p)\V_L]\in\real^{n^2\times L}$ and seeing that $\B = \W\Hob$. The first result then follows from the symmetry of $\Hob$. It follows that, as $\B = \sum_{\alphab}(\xi_{\alphab}'-p)\Va\ha^\top$, that $\B$ is a sum of mean-zero, independent matrices. Therefore, we can develop a high probability bound on $\Vert\B\Vert$ using \Cref{thm:Bernstein}. Defining each summand as $\B_{\alphab}$ we first bound
\begin{align*}
    \Vert\B_{\alphab}\Vert & = \Vert (\xa'-p) \Va\ha^\top\Vert\\
    &\leq \Vert\Va\Vert_2\Vert\ha\Vert_2\\
    & \leq \max_{\betab}\Vert\Vb\Vert_2\Vert\Hob\Vert_{2\to\infty}\\
    &\leq \max_{\betab}\Vert\Vb\Vert_2\sqrt{\frac{2}{n}}=: R,
\end{align*}
where the second inequality follows from \Cref{lem: Hob eigval bound}. Next, we bound the variance expressions. First, we have that
\begin{align*}
    \left\Vert\sum_{\alphab}\Ebb[\Ba\Ba^\ast]\right\Vert & = \left\Vert\sum_{\alphab}\Ebb[(\xa'-p)^2]\Va\ha^\top\ha\Va^\top\right\Vert\\
    & = p(1-p)\left\Vert\sum_{\alphab} \Vert\ha\Vert_2^2\Va\Va^\top\right\Vert\\
    &\leq p \frac{2}{n}\left\Vert\sum_{\alphab}\Va\Va^\top\right\Vert\\
    &= \frac{2p}{n}\lambda_{\max}(\G),
\end{align*}
using \Cref{lem: vectorization} for the final line, followed by
\begin{align*}
       \left\Vert\sum_{\alphab}\Ebb[\Ba^\ast\Ba]\right\Vert &=\left\Vert\sum_{\alphab}\Ebb[(\xa'-p)^2]\ha\Va^\top\Va\ha\right\Vert\\
       & = p(1-p)\left\Vert\sum_{\alphab}\Vert\Va\Vert_2^2\ha\ha^\top\right\Vert\\
       &\leq p\max_{\betab}\Vert\Vb\Vert_2^2\left\Vert\sum_{\alphab}\ha\ha^\top\right\Vert\\
       & = p\max_{\betab}\Vert\Vb\Vert_2^2\left\Vert \Hob (\Hob)^\top\right\Vert\\
       & = p\max_{\betab}\Vert\Vb\Vert_2^2\Vert\Hob\Vert^2\\
       &\leq 4p\max_{\betab}\Vert\Vb\Vert_2^2,
\end{align*}
where the last inequality follows from \Cref{lem: Hob eigval bound}. Now, if $\Va= \vec(\Pt\wa)$, it follows that taking the maximum of the variance terms, we have $\sigma^2 = p\frac{16\nu r}{n}$ and $R = \sqrt{\frac{\nu r}{n^2}}$. Now, as $\frac{\sigma^2}{R} = p\sqrt{16\nu r}\geq \sqrt{p\frac{128\nu r\beta\log{n}}{3n}}$ for $p\geq\frac{8\beta\log{n}}{3n}$, we have from \eqref{eqn:Bernstein} that
\[
\mathbb{P}\left(\Vert\B\Vert\geq \sqrt{p\frac{128\nu r\beta\log{n}}{3n}}\right) \leq 2n^{2-\beta}.
\]

\noindent If $\Va = \vec(\wa)$, then taking the maximum of the variance terms we have $\sigma^2 = 16p$ and $R  = \frac{4}{\sqrt{n}}$. Now, as $\frac{\sigma^2}{R} = 4p\sqrt{n}\geq \sqrt{p\frac{128\beta\log{n}}{3}}$ for $p\geq \frac{4\beta \log{n}}{n}$, it follows again from \eqref{eqn:Bernstein} that
\[
\mathbb{P}\left(\Vert\B\Vert\geq \sqrt{p\frac{128\beta\log{n}}{3}}\right) \leq 2n^{2-\beta}.
\]
This concludes the proof.
\end{proof}

\begin{lem}\label{lem: vector b bound}
Let $\X$ be a rank-$r$, $\nu$-incoherent ground truth matrix with tangent space $\T$. Let $\Va = \vec(\Pt\wa)$ or $\vec(\wa)$. Let $\ba = \sum_{\betab\neq\alphab}(\xi_{\betab}'-p)\Vb$. For $p\geq \frac{2\beta\log{n}}{n}$, with probability at least $1-n^{2-\beta}$ we have that for every $\ai$, if $\Va = \vec(\Pt\wa)$, 
    \[
    \Vert\ba\Vert_2\leq \sqrt{p\frac{4\beta\nu r\log{n}}{n^2}},
    \]
    and if $\Va = \vec(\wa)$
    \[
    \Vert\ba\Vert_2 = \sqrt{p\frac{16\beta\log{n}}{n}}.
    \]
\end{lem}
\begin{proof}
    This result will follow from the vector Bernstein inequality seen in \Cref{thm: Vector Bernstein}. We first note that $\ba$ is a sum of mean-zero independent random vectors. Defining the term $\y_{\betab} = (\xi'_{\betab}-p)H^{\alphab\betab}\Vb$,   we first need to bound $\Vert\y_{\betab}\Vert_2$ via
    \begin{align*}
        \Vert\y_{\betab}\Vert_2 &= \vert(\xi'_{\betab}-p)\vert H^{\alphab\betab}\Vert\Vb\Vert_2\\
        &\leq \frac{1}{n}\max_{\alphab}\Vert\Va\Vert_2=: R.
    \end{align*}
Next, to bound the variance term, we see that
\begin{align*}
   \sum_{\betab\neq\alphab}\Ebb[\Vert\yb\Vert_2^2]  & = \sum_{\betab\neq\alphab} \Ebb [(\xi_{\betab}'-p)^2](H^{\alphab\betab})^2\Vert\Vb\Vert_2^2\\
   & = p(1-p )\sum_{\alphab\neq\betab}(H^{\alphab\betab})^2\\
   &\leq p\max_{\bm{\gamma}}\Vert\bm{V_\gamma}\Vert_2^2\sum_{\alphab\neq\betab}(H^{\alphab\betab})^2\\
   & \leq \frac{2p}{n}\max_{\bm{\gamma}}\Vert\bm{V_\gamma}\Vert_2^2.
\end{align*}
When $\Vb = \vec(\Pt\wb)$, we have $\sigma^2 = p\frac{\nu r}{n^2}$ and $R = \sqrt{\frac{\nu r}{2n^3}}$. It follows then that $\frac{\sigma^2}{R} = p\sqrt{\frac{2\nu r}{n}}\geq \sqrt{p\frac{4\nu r \beta\log{n}}{n^2}}$ for $p\geq \frac{2\beta\log n}{n}$, we have from \Cref{thm: Vector Bernstein} that
\[
\mathbb{P}\left(\Vert\ba\Vert_2\geq \sqrt{p\frac{8\nu r\beta\log{n}}{n^2}}\right)\leq 2n^{-\beta}.
\]
Similarly, when $\Vb = \vec(\wb)$, we have $\sigma^2 = \frac{4p}{n}$ and $R = \frac{2}{n}$. As such, $\frac{\sigma^2}{R} = 2p\geq \sqrt{p\frac{16\beta\log{n}}{3n}}$ for $p\geq \frac{2\beta\log{n}}{n}$, that

\[
\mathbb{P}\left(\Vert\ba\Vert_2\geq \sqrt{p\frac{16\beta\log{n}}{n}}\right)\leq 2n^{-\beta}.
\]

\noindent The stated results then follow from the union bound over every $\ai$.
\end{proof}

\begin{lem}\label{lem: off-diag conditional bernstein bound}
Let $\X$ be a rank-$r$, $\nu$-incoherent ground truth matrix with tangent space $\T$. Then if $p\geq \frac{16 \beta\log{n}}{3n}$ for any $\beta> 3$, with probability at least $1-C'n^{3-\beta}$, where $C'$ is an absolute constant, we have that, for some absolute numerical constants $C_1,~C_2,$ and $C_3$ independent of $n$ or $\X$, the following three inequalities all hold:
\begin{align*} \
 \left\Vert\sum_{\substack{\ai\\\betab\neq\alphab}}(\xa-p)(\xi_{\betab}-p)H^{\alphab\betab}\langle \cdot,\Pt\wb\rangle\Pt\wa \right\Vert&\leq p\frac{C_1\nu r\sqrt{\beta}\log{n} }{n},\\
  \left\Vert\sum_{\substack{\ai\\\betab\neq\alphab}}(\xa-p)(\xi_{\betab}-p)H^{\alphab\betab}\langle \cdot,\Pt\wb\rangle\wa \right\Vert&\leq pC_2\sqrt\frac{\beta\nu r}{n}\log{n},\\
   \left\Vert\sum_{\substack{\ai\\\betab\neq\alphab}}(\xa-p)(\xi_{\betab}-p)H^{\alphab\betab}\langle \cdot,\wb\rangle\wa \right\Vert&\leq pC_3\sqrt{\beta}\log{n}.
    \end{align*}
\end{lem}

\begin{proof}
The goal is to develop a high probability bound for the operator norm of the previously mentioned operator. Prior to doing this, we define some simplifying notation. Let $\wab = \Pt\wa$ or $\wa$, and let $\wah = \Pt\wa$ or $\wa$ with them not necessarily being the same. Let $\G = [\langle \wab,\bar{\w}_{\betab}\rangle]\in\real^{L\times L}$ be the corresponding Gram matrix for $\{\wab\}_{\ai}$. 

We begin the proof by considering, for some $t>0$, the expression
    \[
    \mathbb{P}\bigg[\bigg\Vert\sum_{\alphab\neq\betab}(\xi_{\alphab}-p)(\xi_{\betab}-p)\langle\cdot, \wbh\rangle H^{\alphab\betab}\wab\bigg\Vert\geq t\bigg].
    \]
For notational simplicity, we let $\U_{\alphab} = \vec(\wab)\in\real^{n^2}$ denote the vectorized version of $\wab$. Similarly, we let $\Vb = \vec(\wbh)\in\real^{n^2}$ As such, we can see that this probabilistic expression can be equivalently represented as 
\[
\mathbb{P}\bigg[\bigg\Vert\sum_{\alphab\neq\betab}(\xi_{\alphab}-p)(\xi_{\betab}-p) H^{\alphab\betab}\Ua\Vb^\top\bigg\Vert\geq t\bigg].
\]
This term is a second order matrix chaos, and is a degenerate U-statistic---see \cite{lee2019u} for details on the theory of U-statistics, and see \cite{de2012decoupling} for a detailed investigation of the decoupling of Banach-space-valued U-statistics. As seen in \Cref{thm: Decoupling}, we can study the distribution of the operator norm for this operator by considering a new set of i.i.d. random Bernoulli variables $\xi_{\betab}'$. As such, \Cref{thm: Decoupling} states that for some absolute constant $C$ independent of $n$, $\U_{\alphab}$, $\Va$, or $\H^{-1}$, that
\[
\mathbb{P}\bigg[\bigg\Vert\sum_{\alphab\neq\betab}(\xi_{\alphab}-p)(\xi_{\betab}-p) H^{\alphab\betab}\Ua\Vb^\top\bigg\Vert\geq t\bigg]\leq C\mathbb{P}\bigg[\bigg\Vert\sum_{\alphab\neq\betab}(\xi_{\alphab}-p)(\xi'_{\betab}-p) H^{\alphab\betab}\Ua\Vb^\top\bigg\Vert\geq t/C\bigg].
\]
We note that optimal values of the constants $C$ have not been explicitly computed for each chaos of order $m$, but that the constant itself is a function of $m$. For $m=2$, which is the case we consider, this constant is usually rather small \cite{de2012decoupling}. Using the law of total probability, conditioning on the random variables $\xi'_1,...,\xi'_L$ and defining the vector $\bm{\xi'} =[\xi'_1~...~\xi_L']^\top$ that
\[
\mathbb{P}\bigg[\bigg\Vert\sum_{\alphab\neq\betab}(\xi_{\alphab}-p)(\xi'_{\betab}-p) H^{\alphab\betab}\Ua\Vb^\top\bigg\Vert\geq t\bigg] = \Ebb\bigg[\mathbb{P}\bigg(\bigg\Vert\sum_{\alphab\neq\betab}(\xi_{\alphab}-p)(\xi'_{\betab}-p) H^{\alphab\betab}\Ua\Vb^\top\bigg\Vert\geq t ~\bigg\vert~\bm{\xi}'\bigg)\bigg].
\]

\noindent We can study the conditional expression 
\[
\mathbb{P}\bigg(\bigg\Vert\sum_{\alphab\neq\betab}(\xi_{\alphab}-p)(\xi'_{\betab}-p) H^{\alphab\betab}\Ua\Vb^\top\bigg\Vert\geq t ~\bigg\vert~\bm{\xi}'\bigg)
\]using \Cref{thm:Bernstein}.  Defining the expression $\b_{\alphab} = \sum_{\betab\neq\alphab}H^{\alphab\betab}(\xi'_{\betab}-p)\Vb$, we notice that the sum can be re-written as
\[
\mathbb{P}\bigg(\bigg\Vert\sum_{\alphab\neq\betab}(\xi_{\alphab}-p)(\xi'_{\betab}-p) H^{\alphab\betab}\Ua\Vb^\top\bigg\Vert\geq t ~\bigg\vert~\bm{\xi}'\bigg) = \mathbb{P}\bigg(\bigg\Vert\sum_{\alphab}(\xi_{\alphab}-p)\Ua\bm{b_\alpha}^\top\bigg\Vert\geq t ~\bigg\vert~\bm{\xi}'\bigg).
\]

In this form, conditioned on $\bm{\xi}'$, we see that this is a sum of mean-zero, independent random variables. To study its concentration with \Cref{thm:Bernstein}, we define $\bm{T_\alpha} = (\xi_{\alphab}-p)\Ua\bm{b_\alpha}^\top$. The first task is bounding the operator norm
\[
\Vert\bm{T_\alpha}\Vert \leq \Vert\Ua\Vert_2\Vert \bm{b_\alpha}\Vert_2\leq \max_{\betab}\Vert\Ub\Vert_2\max_{\bm{\gamma}}\Vert\bm{b_\gamma}\Vert_2 =: R.
\]
Next, we need to estimate the variance terms $\left\Vert\sum_{\alphab}\Ebb[\bm{T_\alpha T_\alpha}^\ast~\vert~\bm{\xi}']\right\Vert$ and $\left\Vert\sum_{\alphab}\Ebb[\bm{T_\alpha}^\ast\bm{\bm{T_\alpha}}~\vert~\bm{\xi}']\right\Vert$. For the first term, notice that
\begin{align*}
\left\Vert\sum_{\alphab}\Ebb[\bm{T_\alpha T_\alpha}^\ast~\vert~\bm{\xi}']\right\Vert &= \left\Vert\sum_{\alphab}\Ebb[(\xi_{\alphab} - p)^2]\Ua\ba^\top\ba\Ua^\top\right\Vert\\
& = p(1-p)\left\Vert\sum_{\alphab}\Vert\ba\Vert_2^2\Ua\Ua^\top\right\Vert\\
&\leq p \max_{\alphab}\Vert\ba\Vert_2^2\left\Vert\sum_{\alphab}\Ua\Ua^\top\right\Vert\\
&\leq p\max_{\alphab} \Vert\ba\Vert_2^2 \lambda_{\max}(\G).
\end{align*}
The other variance term is given by
\begin{align*}
    \left\Vert\sum_{\alphab}\Ebb[\bm{T_\alpha}^\ast\bm{\bm{T_\alpha}}~\vert~\bm{\xi}']\right\Vert &= \left\Vert\sum_{\alphab}\Ebb[(\xi_{\alphab}-p)^2]\ba\Ua^\top\Ua\ba^\top\right\Vert\\
    &\leq p \max_{\gamma}\Vert\U_{\bm{\gamma}}\Vert_2^2\left\Vert\sum_{\alphab} \ba\ba^\top\right\Vert\\
    &  = p \max_{\gamma}\Vert\U_{\bm{\gamma}}\Vert_2^2\Vert\B\B^\top\Vert\\  
    & = p\max_{\gamma}\Vert\U_{\bm{\gamma}}\Vert_2^2\Vert\B\Vert^2,
\end{align*}
where the last matrix is defined as $\B = [\b_1 \cdots \b_L]$ and the spectral result follows from the vectorization lemma.

Defining the variance parameter $\sigma^2 = \max\left\{\left\Vert\sum_{\alphab}\Ebb[\bm{T_\alpha T_\alpha}^\ast~\vert~\bm{\xi}']\right\Vert,\left\Vert\sum_{\alphab}\Ebb[\bm{T_\alpha}^\ast\bm{\bm{T_\alpha}}~\vert~\bm{\xi}']\right\Vert\right\}$, we have from the monotonicity of expectation that 
\begin{align*}
    \Ebb\bigg[\mathbb{P}\bigg(\bigg\Vert\sum_{\alphab\neq\betab}(\xi_{\alphab}-p)(\xi'_{\betab}-p)& H^{\alphab\betab}\Ua\Vb^\top\bigg\Vert\geq t ~\bigg\vert~\bm{\xi}'\bigg)\bigg] \leq \Ebb\bigg[2n\exp\left(\frac{-t^2/2}{Rt/3 + \sigma^2}\right)\bigg].
    \end{align*}
It follows that if high probability estimates exist for $R$ and $\sigma^2$, referred heretofore as $R_{\max}$ and $\sigma^2_{\max}$ respectively, that we can effectively uncondition $R$ and $\sigma^2$ as functions of $\bm{\xi}'$ as follows:
\[
 \Ebb\bigg[2n\exp\left(\frac{-t^2/2}{Rt/3 + \sigma^2}\right)\bigg]    \leq 2\mathbb{P}(\sigma^2 \leq \sigma_{\max}^2, R\leq R_{\max})n\exp\left(\frac{-t^2/2}{R_{\max}t/3 +\sigma_{\max}^2}\right) + 2n\mathbb{P}(\sigma^2 \geq \sigma_{\max}^2, R\geq R_{\max}).
\]

Plugging in the expressions in \Cref{eqn: Matrix B bound} and \Cref{lem: vector b bound} and setting $t = \sqrt{\frac{8}{3}\sigma_{\max}^2\log{n}}$  in the expressions such that the simplification in \Cref{thm:Bernstein} is possible, we have that, if $p\geq \frac{16 \beta\log{n}}{n}$ for all three inequalities, the following inequality holds for each operator in the Lemma statement:  
    \begin{align*}
        \Ebb\bigg[2n\exp\left(\frac{-t^2/2}{Rt/3 + \sigma^2}\right)\bigg]  &  \leq 2\mathbb{P}(\sigma^2 \leq \sigma_{\max}^2, R\leq R_{\max})n\exp\left(\frac{-t^2/2}{R_{\max}t/3 +\sigma_{\max}^2}\right) + 2n\mathbb{P}(\sigma^2 \geq \sigma_{\max}^2~\text{or} ~R\geq R_{\max})\\
&\leq 2\mathbb{P}(\sigma^2 \leq \sigma_{\max}^2, R\leq R_{\max})n\exp\left(\frac{-3t^2}{8\sigma_{\max}^2}\right)+ 2n\mathbb{P}(\sigma^2 \geq \sigma_{\max}^2 ~\text{or} ~ R\geq R_{\max})\\
  & \leq 2(1-4n^{2-\beta})n^{1-\beta} + 8n^{3-\beta}\\
  &\leq 10n^{3-\beta}.
\end{align*}
Taking the union bound over the 3 inequalities and incorporating the constant $C$ from \Cref{thm: Decoupling} concludes the proof.
\end{proof}

\begin{lem}\label{lem: operator bound on recentering term}
   Let $\X$ be a ground-truth, rank-r, $\nu$-incoherent Gram matrix with tangent space $\T$. If the sampling probability $p\geq \frac{128\beta\log{n}}{3n}$ for any $\beta>1$, then the following collection of inequalities holds with probability at least $1-8n^{1-\beta}$
   \begin{small}
   \begin{align*}
    \left\Vert\sum_{\ai}\left(\xi_{\alphab}-p\right)\Pt\wa\sum_{\substack{\betab\in\univ\\\alpha\neq\beta}}H^{\alphab\betab}\langle\cdot,\Pt\wb\rangle\right\Vert&\leq \sqrt{p\frac{64\nu^2 r^2 \beta\log{n}}{3n}},&\left\Vert\sum_{\ai}\left(\xi_{\alphab}-p\right)\wa\sum_{\substack{\betab\in\univ\\\alpha\neq\beta}}H^{\alphab\betab}\langle\cdot,\wb\rangle\right\Vert&\leq \sqrt{p\frac{256 n \beta\log{n}}{3}},\\
  \left\Vert\sum_{\ai}\left(\xi_{\alphab}-p\right)\wa\sum_{\substack{\betab\in\univ\\\alpha\neq\beta}}H^{\alphab\betab}\langle\cdot,\Pt\wb\rangle\right\Vert&\leq \sqrt{p\frac{128\nu r \beta\log{n}}{3}}, & \left\Vert\sum_{\ai}\left(\xi_{\alphab}-p\right)\Pt\wa\sum_{\substack{\betab\in\univ\\\alpha\neq\beta}}H^{\alphab\betab}\langle\cdot,\wb\rangle\right\Vert&\leq \sqrt{p\frac{128\nu r \beta\log{n}}{3}}.\\
    \end{align*}
\end{small}%
\end{lem}
\begin{proof}
We first define $\wah = \Pt\wa$ or $\wa$ and $\wab = \Pt\wa$ or $\wa$, where $\wah$ and $\wab$ are not necessarily equal. We also define the corresponding Gram matrices of $\wah$ and $\wab$ as $\hat{\H}$ and $\bar{\H}$, respectively. Similarly define the dual basis vectors $\bar{\bm{v}}_{\alphab}$ and $\hat{\bm{v}}_{\alphab}$. As such, all of these inequalities can be represented through studying the expression
\[
\left\Vert\sum_{\ai}\left(\xi_{\alphab}-p\right)\wah\sum_{\substack{\betab\in\univ\\\alpha\neq\beta}}H^{\alphab\betab}\langle\cdot,\bar{\w}_{\betab}\rangle\right\Vert.
\]

Defining the vectorized version of our summands $\Uah =\vec(\wah)\in\real^{n^2}$ and $\Vab = \vec(\wab)\in\real^{n^2}$, we can see that the above expression is equivalent to 
\[
    \left\Vert\sum_{\ai}\left(\xi_{\alphab}-p\right)\Uah\sum_{\substack{\betab\in\univ\\\alpha\neq\beta}}H^{\alphab\betab}\bar{\V}_{\betab}^\top\right\Vert.
    \]
    
This expression is the operator norm bound on a sum of mean-zero, independent random matrices. We define each summand as $\G_{\alphab} = (\xa-p)\Uah\sum_{\alphab\neq\betab}H^{\alphab\betab}\bar{\V}_{\betab}^\top$.To proceed, we use \Cref{thm:Bernstein}.

    First, we bound $\Vert\G_{\alphab}\Vert$:
    \begin{align*}
        \Vert\G_{\alphab}\Vert &\leq \Vert\Uah\Vert_2 \sum_{\alphab\neq\betab}\vert H^{\alphab\betab}\vert \Vert\bar{\V}_{\betab}^\top\Vert_2\\
        & \leq 2 \max_{\betab}\Vert\Uah\Vert_2 \max_{\bm{\gamma}}\Vert\bar{\bm{V}}_{\bm{\gamma}}\Vert_2,
    \end{align*}
    following from  \Cref{lem: H and H^-1 eigvals}. Next, we bound the variance terms $\Vert\sum_{\alphab}\Ebb[\G_{\alphab}\G_{\alphab}^\ast]\Vert$ and $\Vert\sum_{\alphab}\Ebb[\G_{\alphab}^\ast\G_{\alphab}]\Vert$. It follows that, defining $\d_{\alphab} = \sum_{\alphab\neq\betab}H^{\alphab\betab}\Vb$
    \begin{align*}
      \left\Vert\sum_{\alphab}\Ebb[\G_{\alphab}\G_{\alphab}^\ast]\right\Vert & = \left\Vert\sum_{\alphab}\Ebb[(\xa-p)^2]  \Uah \da^\top\da\Uah^\top\right\Vert\\
      &\leq p\left\Vert\sum_{\alphab}\Uah \da^\top\da\Uah^\top\right\Vert.
    \end{align*}
Now,  $\da + H^{\alphab\alphab}\Vab = \vec(\bar{\bm{v}}_{\alphab})$ by duality and the linearity of both the $\vec$ operator and the projection operator $\Pt$, if present in the relevant definition of $\bar{\bm{v}}_{\alphab}$. Defining $\Qab = \vec(\bar{\bm{v}}_{\alphab})$, we see that
    \begin{align*}
          \left\Vert\sum_{\alphab}\Ebb[\G_{\alphab}\G_{\alphab}^\ast]\right\Vert & \leq p \left\Vert \sum_{\alphab} \Uah\left(\Qab-H^{\alphab\alphab}\Vab\right)^\top\left(\Qab-H^{\alphab\alphab}\Vab\right)\Uah^\top\right\Vert\\
          & = p \bigg\Vert\sum_{\alphab} \Uah\Qab^\top\Qab\Uah^\top - H^{\alphab\alphab}\sum_{\alphab}\Uah\Qab^\top\Vab\Uah^\top \\
          & \qquad - H^{\alphab\alphab}\sum_{\alphab}\Uah\Vab^\top\Qab\Uah^\top + (H^{\alphab\alphab})^2\sum_{\alphab}\Uah\Vab^\top\Vab\Uah^\top \bigg\Vert\\
          &\leq p \left\Vert\sum_{\alphab} \Uah\Qab^\top\Qab\Uah^\top\right\Vert +2p \left\Vert H^{\alphab\alphab}\sum_{\alphab}\Uah\Qab^\top\Vab\Uah^\top\right\Vert\\
          & ~+ p\left\Vert(H^{\alphab\alphab})^2\sum_{\alphab}\Uah\Vab^\top\Vab\Uah^\top \right\Vert \\
          & \leq  p \max_{\betab}\Vert\bar{\bm{Q}}_{\betab}\Vert_2^2\left\Vert\sum_{\alphab} \Uah\Uah^\top\right\Vert +2p \left\Vert H^{\alphab\alphab}\sum_{\alphab}\vert\Qab^\top\Vab\vert\Uah\Uah^\top\right\Vert\\
          & ~+ p\max_{\betab}\Vert\bar{\V}_{\betab}\Vert_2^2\left\Vert(H^{\alphab\alphab})^2\sum_{\alphab}\Uah\Uah^\top \right\Vert \\
          & \leq p\left(\max_{\betab}\Vert\bar{\bm{Q}}_{\betab}\Vert_2^2+\max_{\betab}\Vert\bar{\V}_{\betab}\Vert_2^2 +2\max_{\betab}\Vert\bar{\bm{Q}}_{\betab}\Vert_2\max_{\bm{\gamma}}\Vert\bar{\bm{V}}_{\bm{\gamma}}\Vert_2\right)\left\Vert\sum_{\alphab}\Uah\Uah^\top\right\Vert\\
          & \leq p\left(\max_{\betab}\Vert\bar{\bm{Q}}_{\betab}\Vert_2^2+\max_{\betab}\Vert\bar{\V}_{\betab}\Vert_2^2 +2\max_{\betab}\Vert\bar{\bm{Q}}_{\betab}\Vert_2\max_{\bm{\gamma}}\Vert\bar{\bm{V}}_{\bm{\gamma}}\Vert_2\right)\lambda_{\max}(\hat{\H}),
    \end{align*}
    where the first inequality follows from the triangle inequality, the second inequality follows from Cauchy-Schwarz, the third inequality follows from the fact that $\vert H^{\alphab\alphab}\vert<1$ for every $\ai$, and the last inequality follows from \Cref{lem: vectorization} and \Cref{lem: H and H^-1 eigvals}.
    
  To analyze the other variance term, we begin by defining the matrix $\E = [\d_1 \cdots\d_L]\in\real^{n^2\times L}$. We now see that
    \begin{align*}
         \left\Vert\sum_{\alphab}\Ebb[\G_{\alphab}^\ast\G_{\alphab}]\right\Vert &= \left\Vert\sum_{\alphab}\Ebb[(\xa-p)^2]\da\Uah^\top\Uah\da^\top\right\Vert\\
         & \leq p\max_{\betab}\Vert\Uah\Vert_2^2\left\Vert\sum_{\alphab}\da\da^\top\right\Vert\\
         & = p\max_{\betab}\Vert\Uah\Vert_2^2\Vert \E\E^\top\Vert.
    \end{align*}
    Similarly to \Cref{eqn: Matrix B bound}, we can see that $\E = \W\Hob$. As such, notice that
    \begin{align*}
        \Vert\E\E^\top\Vert &= \max_{\Vert\y\Vert_2 =1} \y^\top\E\E^\top\y\\
        & = \max_{\Vert\y\Vert_2 =1} \y^\top\W\Hob\Hob\W^\top\y\\
        & = \max_{\Vert \y\Vert_2=1} \Vert\Hob\W^\top\y\Vert_2^2\\
        &\leq  \Vert\Hob\Vert^2\max_{\Vert \y\Vert_2=1} \Vert\W^\top \y\Vert_2^2\\
        & \leq 4 \max_{\Vert \y\Vert_2=1} \y^\top\W\W^\top\y\\
        & = 4\Vert\W\W^\top\Vert\\
        & = 4\Vert\W^\top\W\Vert\\
        & = 4\lambda_{\max}(\bar{\H}),
    \end{align*}
    where the first inequality follows from Cauchy-Schwarz, and the second inequality follows from \Cref{lem: Hob eigval bound}. This yields a variance bound as follows:
    \begin{align*}
    \sigma^2 &= \max \left\{  \left\Vert\sum_{\alphab}\Ebb[\G_{\alphab}^\ast\G_{\alphab}]\right\Vert ,  \left\Vert\sum_{\alphab}\Ebb[\G_{\alphab}\G_{\alphab}^\ast]\right\Vert \right\} \\
    &\leq p\max\left\{4\max_{\betab}\Vert\Uah\Vert_2^2\lambda_{\max}(\bar{\H}),\left(\max_{\betab}\Vert\bar{\bm{Q}}_{\betab}\Vert_2^2+\max_{\betab}\Vert\bar{\V}_{\betab}\Vert_2^2 +2\max_{\betab}\Vert\bar{\bm{Q}}_{\betab}\Vert_2\max_{\bm{\gamma}}\Vert\bar{\bm{V}}_{\bm{\gamma}}\Vert_2\right)\lambda_{\max}(\hat{\H})\right\}.
    \end{align*}
This leads to the following variance bounds and norm bounds for the corresponding combinations of $\wab$ and $\wah$:
\begin{align*}
    \wab = \Pt\wa,~\wah = \Pt\wa & \Rightarrow \sigma^2 = 8p\frac{\nu^2 r^2}{n}, ~R = \frac{\nu r}{n},\\
    \wab = \Pt\wa,~\wah = \wa &\Rightarrow \sigma^2 = 64p\nu r,~R = \sqrt{\frac{8\nu r}{n}},\\
    \wab = \wa,~\wah=\Pt\wa & \Rightarrow \sigma^2 = 16p\nu r,~R = \sqrt{\frac{8\nu r}{n}},\\
    \wab = \wa,~\wah = \wa &\Rightarrow \sigma^2 = 32pn,~R = 4.
\end{align*}

Setting $t = \sqrt{\frac{8}{3}\sigma^2 \beta\log{n}}$ for each expression, which uniformly satisfies $t\leq \frac{\sigma^2}{R}$ for each bound as long as $p\geq \frac{128\beta\log{n}}{3n}$, yields each of the desired inequalities with failure probabilities of $1-2n^{1-\beta}$ for each expression using \Cref{thm:Bernstein}. This concludes the proof.    
    \end{proof}

\subsection{Proof of Theorem \ref{thm: RIP of Mo}}\label{subsubsec: proof of Mo RIP}
We are now ready to prove \Cref{thm: RIP of Mo}.
\begin{proof}
    First, notice that since $\Mo = \Mo^\ast$,
    \begin{align*}
        &\Vert\Pt\Mo\Pt - p^2\Pt\Vert = \left\Vert \sum_{\alphab,\betab\in\Omega} \xi_{\alphab}\xi_{\betab}C_{\alphab\betab}\Pt\wa\la\va,\vb\ra\la\cdot,\Pt\wb\ra -  p^2\Pt\right\Vert\\
        & =  \left\Vert \sum_{\alphab,\betab\in\Omega} \xi_{\alphab}\xi_{\betab}C_{\alphab\betab}\Pt\wa\la\va,\vb\ra\la\cdot,\Pt\wb\ra -  p^2\sum_{\alphab,\betab\in\univ}\Pt\wa\la\va,\vb\ra\la\cdot,\Pt\wb\ra \right\Vert\\
        & =\left\Vert p\sum_{\alphab\in\Omega} \xi_{\alphab}\la\cdot,\Pt\wa\ra \Pt\wa\la\va,\va\ra  +\sum_{\substack{\alphab, \betab \in \Omega \\ \alphab \neq \betab}}
        \xi_{\alphab}\xi_{\betab}\Pt\wa\la\va,\vb\ra\la\cdot,\Pt\wb\ra - p^2\sum_{\alphab,\betab\in\univ}\Pt\wa\la\va,\vb\ra\la\cdot,\Pt\wb\ra \right\Vert\\
        & = \left\Vert p\Vert\va\Vert_\fro^2 \sum_{\ai}(\xa-p)\langle\Y,\Pt\wa\rangle \Pt\wa+ \sum_{\alphab\neq\betab}(\xa\xi_{\betab}-p^2)\Pt\wa H^{\alphab\betab}\langle\cdot,\Pt\wb\rangle\right\Vert\\
        & \leq \left\Vert p\Vert\va\Vert_\fro^2 \sum_{\ai}(\xa-p)\langle\Y,\Pt\wa\rangle \Pt\wa\right\Vert + \left\Vert\sum_{\alphab\neq\betab}(\xa\xi_{\betab}-p^2)\Pt\wa H^{\alphab\betab}\langle\cdot,\Pt\wb\rangle\right\Vert\\
        & = {p\Vert\va\Vert_\fro^2 \left\Vert\Pt\Fo\Pt-p\Pt\FI\Pt\right\Vert} + \left\Vert\sum_{\alphab\neq\betab} (\xa\xi_{\betab}-p^2)\Pt\wa H^{\alphab\betab}\langle \cdot,\Pt\wb\rangle\right\Vert\\
        & \leq \underbrace{p\Vert\va\Vert_\fro^2 \left\Vert\Pt\Fo\Pt-p\Pt\FI\Pt\right\Vert}_{B_1} + \underbrace{\left\Vert\sum_{\alphab\neq\betab}(\xa-p)(\xi_{\betab}-p)\Pt\wa H^{\alphab\betab}\langle\cdot,\Pt\wb\rangle\right\Vert }_{B_2}\\
        & \qquad+
 \underbrace{2p\left\Vert\sum_{\ai}\left(\xi_{\alphab}-p\right)\Pt\wa\sum_{\substack{\betab\in\univ\\\alphab\neq\betab}}H^{\alphab\betab}\langle\cdot,\Pt\wb\rangle\right\Vert}_{B_3},
    \end{align*}

\noindent where the first inequality follows from the triangle inequality and the final inequality follows from \Cref{lem: Bound on So H So - E() term}. We note that the term $B_1$ is bounded by \Cref{lem: Concentration of Fo}, and by the fact that $\Vert\va\Vert_\fro$ is constant for all $\ai$ and that $\Vert\va\Vert_\fro^2\leq \frac{1}{2}$ from \Cref{lem: H and H^-1 eigvals}. We bound $B_2$ using \Cref{lem: off-diag conditional bernstein bound}, and we bound $B_3$ using \Cref{lem: operator bound on recentering term}. This concludes the proof.
\end{proof}

\begin{lem}\label{lem: norm of Mo}
    Let $\Omega\subset\univ$ be sampled with uniform Bernoulli probability $p$. If $p\geq\frac{128}{3}\frac{\beta\log{n}}{n}$ for any $\beta>3$, then with probability at least $1-C'n^{3-\beta}$, where $C'$ is an absolute constant, we have that
    \[
    \Vert\Mo\Vert \leq p^2 \left(1 + 20\sqrt{\frac{\beta n\log{n}}{3p}}\right) + Cp\beta^{1/2}\log{n},
    \]
    for some absolute constant $C>0$. Furthermore, for $p\geq \frac{C''\beta\log{n}}{n}$ for some sufficiently large $C''>0$, this simplifies to
    \[
    \Vert\Mo\Vert \leq 50p^{3/2}\sqrt{\beta n\log{n}}.
    \]
\end{lem}
\begin{proof}
    This proof follows directly from Lemmas~\ref{lem: Concentration of Fo} and~\ref{lem: Bound on So H So - E() term}. First, notice that
    \[
    \Vert\Mo\Vert = \Vert\Mo-p^2\Id+p^2\Id\Vert\leq p^2 +\Vert\Mo-p^2\Id\Vert.
    \]
    This second term can be analyzed in the same way as in the proof of \Cref{thm: RIP of Mo}, seen in \Cref{subsubsec: proof of Mo RIP}:
    \begin{align*}
        \Vert\Mo-p^2\Id\Vert\leq p\Vert\va\Vert_\fro^2\Vert\Fo-p\mathcal{F}_{\univ}\Vert + \left\Vert\sum_{\alphab\neq\betab}(\xa\xi_{\betab}-p^2)\wa H^{\alphab\betab}\langle\cdot,\wb\rangle\right\Vert.
    \end{align*}
    Applying the results of \Cref{lem: Concentration of Fo}, \Cref{lem: Bound on So H So - E() term}, \Cref{lem: off-diag conditional bernstein bound}, and \Cref{lem: operator bound on recentering term} gives the desired result. 

    For the simplification, we notice that
    \begin{align*}
            p^2 \leq 5p^{3/2}\sqrt{\beta n\log{n}}\\
        \Leftrightarrow p\leq 25 \beta n\log{n}
    \end{align*}
    and that 
    \begin{align*}
        25p^{3/2}\sqrt{\beta n\log{n}}&\geq Cp\beta^{1/2}\log{n}\\
        \Leftrightarrow p &\geq \frac{C^2}{625}\frac{\log{n}}{n},
    \end{align*}
    so for some sufficiently large $C'>0$ the stated simplification holds.
\end{proof}

\begin{lem}\label{lem: PtMo Bounds}
    Let $\Omega\subset \univ$ be sampled with uniform Bernoulli probability $p$, and let $\mathbb{T}$ be the tangent space on $\mfr$ for a rank-$r$, $\nu$-incoherent ground truth matrix $\X$. If $p\geq \frac{128}{3}\frac{\beta\log{n}}{{n}}$ for any $\beta>3$, then with probability at least $1-C'n^{3-\beta}$, where $C'$ is an absolute constant, we have that, for some absolute constant $C>0$, 
    \[
    \Vert\Mo\Pt\Vert \leq p^2 +Cp\sqrt{\frac{\nu r\beta}{n}}\log{n}+ 20p^{3/2}\sqrt{\nu r\beta\log{n}}.
    \]
    If $p\geq \frac{C''\log{n}}{n}$ for some sufficiently large constant $C''>0$ independent of $n,~\nu$ and $r$, this simplifies to
    \[
    \Vert\Mo\Pt\Vert \leq 50p^{3/2}\sqrt{\nu r\beta\log{n}}.
    \]
Furthermore, we have that
    \begin{align*}
        \Vert\Mo\Ptl\Vert \leq 2\Vert\Mo\Vert\frac{\Vert\X_l - \X\Vert_\fro}{\lambda_r(\X)} + \Vert \Mo\Pt\Vert,
    \end{align*}
    and with probability at least $1-14n^{3-\beta}$, for some sufficiently large constant $C'''>0$ independent of $n$, $\nu$ and $r$, if $p\geq \frac{C''' \log{n}}{n}$, then
    \[
    \Vert\Mo\Ptl\Vert\leq 100 p^{3/2}\sqrt{\beta n\log{n}}\frac{\Vert\X_l-\X\Vert_\fro}{\lambda_r(\X)}+50p^{3/2}\sqrt{\beta\log{n}}.
    \]
\end{lem}
\begin{proof}
    We first notice that
    \begin{align*}
        \Vert\Mo\Pt\Vert& \leq \Vert\Mo\Pt - \Ebb[\Mo\Pt]\Vert + \Vert\Ebb[\Mo\Pt]\Vert\\
        & = \Vert\Mo\Pt - \Ebb[\Mo\Pt]\Vert + p^2.
    \end{align*}
    Similarly to the proofs of \Cref{thm: RIP of Mo} and \Cref{lem: norm of Mo}, we can decompose the difference between $\Mo\Pt$ and $\Ebb[\Mo\Pt]$ as the concentration of $\Fo\Pt-\Ebb[\Fo]\Pt$ and the off-diagonal term. As such, we can see that
    \begin{align*}
        \Vert\Mo\Pt\Vert&\leq p^2 +  p\Vert\va\Vert_F^2 \left\Vert \Fo\Pt-p\FI\Pt\right\Vert+ \left\Vert\sum_{\substack{\alphab,
        \betab\in\univ\\\alphab\neq\betab}} (\xi_{\alphab}\xi_{\betab}-p^2)\wa\la\va,\vb\ra\la\cdot,\Pt\wb\ra\right\Vert\\
    \end{align*}
    From \Cref{lem: concentration of FoPt}, \Cref{lem: off-diag conditional bernstein bound}, and \Cref{lem: operator bound on recentering term}, the first result follows. The simplification follows a similar strategy as in \Cref{lem: norm of Mo}.
    For the last result, notice that
    \begin{align*}
        \Vert\Mo\Ptl\Vert &= \Vert\Mo\Ptl - \Mo\Pt + \Mo\Pt\Vert\\
        &\leq \Vert\Mo(\Ptl-\Pt)\Vert+\Vert\Mo\Pt\Vert\\
        &\leq\Vert\Mo\Vert\frac{2\Vert \X_l-\X\Vert_\fro}{\lambda_r(\X)} + \Vert\Mo\Pt\Vert.
    \end{align*}
    For the final result, the result follows from using the simplified expressions for both $\Vert\Mo\Pt\Vert$ and $\Vert\Mo\Vert$, picking a sufficiently large constant such that both simplifications hold.
\end{proof}

\begin{lem}[Local RIP of $\Mo$]\label{lem: Local RIP Mo}
    Assume that
    \begin{align}
      p^{-2}\Vert\Pt\Mo\Pt-p^2\Pt\Vert&\leq\varepsilon_0\label{eqn: Mo RIP assumption appendix} , \\
        \Vert\Mo\Pt\Vert&\leq  50p^{3/2}\sqrt{\beta \nu r\log{n}},\label{eqn: MoPt Bound} \\
      \Vert\Mo\Ptl\Vert&\leq 100 p^{3/2}\sqrt{\beta n\log{n}}\frac{\Vert\X_l-\X\Vert_\fro}{\lambda_r(\X)}+50p^{3/2}\sqrt{\beta \nu r\log{n}},\label{eqn: MoPtl Bound}\\
     \frac{\Vert\X_l-\X\Vert_\fro}{\lambda_r(\X)}&\leq\frac{\varepsilon_0 p^{1/2}}{100\left(\beta n\log{n}\right)^{1/4}}.\label{eqn: Mo Alg Local Neighborhood}
    \end{align}
    Then
    \[
    p^{-2}\Vert\Ptl\Mo\Ptl-p^2\Ptl\Vert\leq 4\varepsilon_0.
    \]    
\end{lem}
\begin{proof}
     \begin{align*}
        \Vert \Ptl - &p^{-2}\Ptl\Mo\Ptl\Vert = \Vert \Ptl - p^{-2}\Ptl\Mo\Ptl \\
        & \qquad\qquad\qquad\qquad+ \Pt - \Pt + p^{-2}\Ptl\Mo\Pt-p^{-2}\Ptl\Mo\Pt + p^{-2}\Pt\Mo\Pt - p^{-2}\Pt\Mo\Pt \Vert \\
        &\leq \Vert \Ptl-\Pt\Vert +p^{-2}\Vert \Ptl\Mo\Ptl - \Ptl\Mo\Pt\Vert 
+ p^{-2}\Vert \Ptl\Mo\Pt - \Pt\Mo\Pt\Vert + \left\Vert\Pt-p^{-2}\Pt\Mo\Pt\right\Vert \\
        &\leq \Vert \Ptl - \Pt\Vert + p^{-2}\Vert \Ptl\Mo\Vert\Vert\Ptl-\Pt\Vert \\
        &\quad+p^{-2}\Vert\Mo\Pt\Vert \Vert \Ptl-\Pt\Vert + \left\Vert \Pt - p^{-2}\Pt\Mo\Pt\right\Vert \\
       & \leq \frac{2\Vert \X_l-\X\Vert_\fro}{\lambda_r(\X)}\left(1 + p^{-2}\Vert \Ptl\Mo\Vert +p^{-2}\Vert\Mo\Pt\Vert\right) + \left\Vert\Pt-p^{-2}\Pt\Mo\Pt\right\Vert \label{eqn: Last equation local PI proof}\\
       &\leq 2\varepsilon_0 + p^{-2}\frac{2\Vert \X_l-\X\Vert_\fro }{\lambda_r(\X)}\left(100 p^{3/2}\sqrt{\beta n\log{n}}\frac{\Vert\X_l-\X\Vert_\fro}{\lambda_r(\X)}+100p^{3/2}\sqrt{\beta \nu r\log{n}}\right)\\
       &= 2\varepsilon_0 + 200p^{-1/2}\sqrt{\beta n\log{n}}\left(\frac{\Vert\X_l-\X\Vert_\fro}{\lambda_r(\X)}\right)^2 + 100p^{-1/2}\sqrt{\beta \nu r\log{n}}\frac{\Vert\X_l-\X\Vert_\fro}{\lambda_r(\X)}\\
       &\leq 4\varepsilon_0,
    \end{align*}
    where the first inequality is the triangle inequality, the second inequality is Cauchy-Schwarz, the third inequality is due to \Cref{lem: Projection Bounds} and \eqref{eqn: Mo RIP assumption appendix}, the fourth inequality is a result of \eqref{eqn: MoPt Bound} and \eqref{eqn: MoPtl Bound}, and the final inequality is due to \eqref{eqn: Mo Alg Local Neighborhood}.
\end{proof}

\section{Local Convergence Results}\label{appendix: local convergence}
We begin with the following technical lemmas used in the proof of local convergence.
\begin{lem}[\Cref{alg: M_omega descent} Stepsize Bounds]\label{lem: Mo Stepsize}
Assume that $\Vert \Ptl - p^{-2}\Ptl\Mo\Ptl\Vert\leq 4\varepsilon_0<1$. Then the stepsize $\alpha_l$ in \Cref{alg: M_omega descent} can be bounded by
\begin{equation*}
    \frac{p^{-2}}{1+4\varepsilon_0}\leq \alpha_l = \frac{\Vert \Pt\G_l\Vert_\fro^2}{\langle \Ptl\G_l,\Mo\Ptl\G_l\rangle}\leq \frac{p^{-2}}{1-4\varepsilon_0}.
\end{equation*}
\end{lem}
\begin{proof}
We will prove this by leveraging the local RIP assumption. Notice the following:
\begin{align*}
        \langle \Ptl\G_l,\Mo\Ptl\G_l\rangle &= \langle \Ptl\G_l,\Ptl\Mo\Ptl\G_l\rangle\\
        & = \left\langle \Ptl\G_l,\Ptl\Mo\Ptl\G_l-p^{2}\Ptl\G_l \right\rangle + p^{2}\la\Pt\G_l,\Pt\G_l\ra.
\end{align*}
We can now leverage the variational characterization of the spectral norm and local RIP, proven in \Cref{lem: Local RIP Mo}, to bound the following:
\[
-p^2(4\varepsilon_0)\Vert \Pt\G_l\Vert_\fro^2 \leq \left\langle \Ptl\G_l,\Ptl\Mo\Ptl\G_l-p^2\Ptl\G_l \right\rangle \leq p^2(4\varepsilon_0)\Vert \Pt\G_l\Vert_\fro^2 .
\]
As such, we can now bound the denominator as 
\[
p^2(1-4\varepsilon_0)\Vert \Pt\G_l\Vert_\fro^2\leq\langle \Ptl\G_l,\Mo\Ptl\G_l\rangle \leq p^2(1+4\varepsilon_0)\Vert \Pt\G_l\Vert_\fro^2.
\]
Rearrangement of this last expression yields the upper and lower bounds on the step size derived above.
\end{proof}

\begin{lem}[$I_1$ Bound]\label{lem: Mo Alpha spectral norm bound}
    Assume $\left\Vert \Ptl - p^{-2}\Ptl\Mo\Ptl\right\Vert\leq 4\varepsilon_0$ and $\alpha_l$ can be bounded as in \Cref{lem: Mo Stepsize}. Then the spectral norm of $\Ptl - \alpha_l\Ptl\Mo\Ptl$ can be bounded as
    \begin{equation}
        \Vert \Ptl-\alpha_l\Ptl\Mo\Ptl\Vert \leq \frac{8\varepsilon_0}{1-4\varepsilon_0}.
    \end{equation}
\end{lem}
\begin{proof}
    From direct calculation, it follows that
    \begin{align*}
        \Vert \Ptl - \alpha_l\Ptl\Mo\Ptl\Vert &\leq \left\Vert\Ptl - p^{-2}\Ptl\Mo\Ptl\right\Vert + \left\vert \alpha_l -p^{-2}\right\vert\Vert\Ptl \Mo\Ptl\Vert\\
        &\leq 4\varepsilon_0+ \left\vert \alpha_l -p^{-2}\right\vert\left(\left\Vert\Ptl \Mo\Ptl -p^2\Ptl\right\Vert +p^2\Vert\Ptl\Vert \right)\\
        &\leq 4\varepsilon_0 + \left(\frac{p^{-2}}{1-4\varepsilon_0} - \frac{p^{-2}(1-4\varepsilon_0)}{1-4\varepsilon_0}\right)\left(\left\Vert\Ptl \Mo\Ptl -p^2\Ptl\right\Vert +p^2\Vert\Ptl\Vert \right)\\
        &\leq  4\varepsilon_0 + \left(\frac{p^{-2}}{1-4\varepsilon_0} - \frac{p^{-2}(1-4\varepsilon_0)}{1-4\varepsilon_0}\right)\left(4\varepsilon_0p^2+ p^2 \right)\\
        &= 4\varepsilon_0 +  \frac{4\varepsilon_0}{1-4\varepsilon_0}(1+4\varepsilon_0)\\
        &=\frac{8\varepsilon_0}{1-4\varepsilon_0}, 
    \end{align*}
    where the first inequality comes from the triangle inequality, the second inequality comes from Local RIP in \Cref{lem: Local RIP Mo}, the third inequality comes from the stepsize bound in \Cref{lem: Mo Stepsize}, the fourth inequality again comes from \Cref{lem: Local RIP Mo}, and the remainder comes from algebraic simplification of terms. This finishes the proof.
\end{proof}

\subsection{Proof of Theorem \ref{thm: Mo Local convergence}}\label{subsubsec: Mo Local convergence}

We can now prove \Cref{thm: Mo Local convergence}.
\begin{proof}\label{proof: Mo Local convergence}
    First, it follows that
    \begin{equation*}
        \Vert \X_{l+1} - \X\Vert_\fro \leq \Vert\X_{l+1} -\W_l\Vert_\fro + \Vert \W_l - \X\Vert_\fro\leq 2\Vert \W_l - \X\Vert_\fro,
    \end{equation*}
    as $\X_{l+1}$ is the best rank-$r$ approximation of $\W_l$. Plugging in $\W_l = \X_l + \alpha_l\Ptl\G_l$, we see that
    \begin{align*}
        \Vert \X_{l+1}-\X\Vert_\fro&\leq 2\left\Vert \X_l + \alpha_l\Ptl\G_l - \X\right\Vert_\fro\\
        &=2\Vert \X_l-\X -\alpha_l\Ptl\Mo(\X_l-\X)\Vert_\fro\\
        &\leq \underbrace{2\Vert (\Ptl - \alpha_l\Ptl\Mo\Ptl)(\X_l-\X)\Vert_\fro}_{I_1}\\
        &\quad+ \underbrace{2\Vert(\mathcal{I}-\Ptl)(\X_l-\X)\Vert_\fro}_{I_2}\\
        &\quad+ \underbrace{2\vert\alpha_l\vert\Vert\Ptl\Mo(\mathcal{I}-\Ptl)(\X_l-\X)\Vert_\fro}_{I_3}.
    \end{align*}
    It remains to bound each term individually. Using \Cref{lem: Mo Alpha spectral norm bound}, we see that
    \begin{equation*}
        I_1\leq \frac{16\varepsilon_0}{1-4\varepsilon_0}\Vert\X_l-\X\Vert_\fro.
    \end{equation*}
    Next, notice that from \Cref{lem: Projection Bounds} and the fact that $\Ptl\X_l = \X_l$,
    \begin{align*}
        I_2 &= 2\Vert(\mathcal{I}-\Ptl)\X_l-(\mathcal{I}-\Ptl)\X)\Vert_\fro\\
        &= 2\Vert(\mathcal{I}-\Ptl)\X\Vert_\fro\\
        &\leq\frac{2\Vert \X_l-\X\Vert_\fro^2}{\lambda_r(\X)}\\
        &\leq \frac{\varepsilon_0 p^{1/2}}{32\left(\beta n\log{n}\right)^{1/4}}\Vert \X_l-\X\Vert_\fro\\
        &\leq \varepsilon_0\Vert\X_l-\X\Vert_\fro\\
        &\leq \frac{\varepsilon_0}{1-4\varepsilon_0}\Vert \X_l-\X\Vert_\fro,
    \end{align*}
    using \Cref{lem: Projection Bounds} and our initial local neighborhood assumption.
    Finally, we see that, following a similar argument as in the bound of $I_2$,
    \begin{align*}
        I_3&\leq 2\vert\alpha_l\vert\Vert \Ptl\Mo\Vert\Vert(I-\Ptl)\X\Vert_\fro\\
        &\leq \frac{2p^{-2}}{1-4\varepsilon_0}\left[100 p^{3/2}\sqrt{\beta n\log{n}}\frac{\Vert\X_l-\X\Vert_\fro}{\lambda_r(\X)}+25p^{3/2}\sqrt{\beta \nu r\log{n}}\right]\left(\frac{\Vert\X_l-\X\Vert_\fro}{\lambda_r(\X)}\right)\Vert\X_l-\X\Vert_\fro\\
        &\leq\frac{2}{1-4\varepsilon_0}\left[100 p^{-1/2}\sqrt{\beta n\log{n}}\left(\frac{\varepsilon_0 p^{1/2}}{100\left(\beta n\log{n}\right)^{1/4}}\right)^2+50p^{-1/2}\sqrt{\beta \nu r\log{n}}\frac{\varepsilon_0 p^{1/2}}{100\left(\beta n\log{n}\right)^{1/4}}\right]\Vert\X_l-\X\Vert_\fro\\
        &\leq\frac{\varepsilon_0}{1-4\varepsilon_0}\Vert\X_l-\X\Vert_\fro
    \end{align*}
    where the second-to-last inequality follows from the same analysis conducted in \Cref{lem: Local RIP Mo}, just divided by 2. Collecting these results, we get
    \begin{equation*}
        \Vert\X_{l+1}-\X\Vert_\fro\leq\frac{18\varepsilon_0}{1-4\varepsilon_0}\Vert\X_l-\X\Vert_\fro.
    \end{equation*}
    By the assumption of the theorem, which holds for $l=0$, and as we have a contractive sequence, it inductively follows that the assumption holds for $l\geq 0$. This concludes the proof.
\end{proof}
\section{Initialization Results (Proof of Lemma \ref{lem: Mo 1SHT initialization})}\label{appendix: initialization}

\begin{proof}
     First, notice that for $\W_0 = p^{-1}\Ro(\X)$, we get
    \begin{align*}
        \left\Vert \X_0 - \X\right\Vert&\leq \left\Vert \W_0 -\X\right\Vert + \left\Vert \W_0-\X_0\right\Vert\\
        &\leq 2\left\Vert \W_0-\X\right\Vert,
    \end{align*}
    where the first inequality follows from the triangle inequality and the second inequality follows from the fact that $\W_0$ is the best rank-$r$ approximation of $\X_0$ by Eckart-Young-Mirsky\cite{eckart1936approximation}. We now need a bound for this last term. Notice that $\W_0-\X = \sum_{\ai}(p^{-1}\xi_{\alphab}-1)\langle \X,\wa\rangle \va$ is a sum of zero-mean i.i.d.~random matrices, opening up the use of Bernstein's inequality. In order to use this, define $\Z_{\alphab} = (p^{-1}\xi_{\alphab}-1)\langle \X,\wa\rangle\va$. We need a bound on $\left\Vert \Z_{\alphab}\right\Vert$ and $\left\Vert \mathbb{E}\left[\sum_{\ai}\Z_{\alphab}\right]^2\right\Vert$.
    First, notice that 
    \begin{align*}
        \left\Vert \Z_{\alphab} \right\Vert &= \left\Vert (p^{-1}\xi_{\alphab}-1)\langle \X,\wa\rangle\va\right\Vert\\
        &\leq (p^{-1}+1)\vert\langle \X,\wa\rangle\vert \left\Vert \va\right\Vert \\
        &\leq 2p^{-1}\left(\max_{\ai}\left\vert\langle \X,\wa\right\rangle\vert \right) =:c,
    \end{align*}
    where the second inequality comes from the fact that $p\leq 1$ and $\left\Vert \va\right\Vert<1$ from \Cref{lem: H and H^-1 eigvals}.
    Next, notice that 
    \begin{align*}
        \left\Vert\mathbb{E}\left[\sum_{\ai}\Z_{\alphab}^2\right]\right\Vert &= \left\Vert\sum_{\ai}\Ebb\left[p^{-2}\xi_{\alphab}-2\xi_{\alphab}p^{-1}+1\right] \langle \X,\wa\rangle^2\va^2\right\Vert\\
        & =\left\Vert\sum_{\ai}\left(p^{-1}-1\right) \langle \X,\wa\rangle^2\va^2\right\Vert\\
        &\leq p^{-1}\left(\max_{\ai}\left\vert\langle \X,\wa\right\rangle\vert \right)^2\lambda_\text{max}\left(\sum_{\alphab}\va^2\right).
    \end{align*}
    Now, as \Cref{lem:Form of va^2}, $\sum_{\alphab} \va^2 = \frac{n^2 - 2n +2}{4n}\J$. It follows that $\lambda_\text{max}\left(\sum_{\alphab} \va^2 \right) = \frac{n^2-2n+2}{4n}\leq \frac{n}{4}$ as $\J$ is an orthogonal projection matrix. Thus,
    \begin{equation*}
        \left\Vert\mathbb{E}\left[\sum_{\ai}\Z_{\alphab}^2\right]\right\Vert\leq \frac{np^{-1}}{4}\left(\max_{\ai}\left\vert\langle \X,\wa\right\rangle\vert \right) =: \sigma^2
    \end{equation*}
    Now to determine $t$, we note that
    \begin{align*}
        \frac{\sigma^2}{c} &= \frac{np^{-1}\left(\max_{\ai}\left\vert\langle \X,\wa\right\rangle\vert \right)^2}{8p^{-1}\left(\max_{\ai}\left\vert\langle \X,\wa\right\rangle\vert \right)}\\
        &= \frac{n}{8}\left(\max_{\ai}\left\vert\langle \X,\wa\right\rangle\vert \right)\\
        &\geq \sqrt{\frac{2\beta n\log{n}}{3p}}\left(\max_{\ai}\left\vert\langle \X,\wa\right\rangle\vert \right),
    \end{align*}
    for $p\geq \frac{128\beta \log{n}}{3n}$. It follows that
    \begin{align*}
        \bb{P}\left(\left\Vert\X_0-\X\right\Vert>\sqrt{\frac{2\beta n\log{n}}{3p}}\left(\max_{\ai}\left\vert\langle \X,\wa\right\rangle\vert \right)\right)&\leq 2n \exp\left(-\beta\log(n)\right)\\
        &= 2n^{1-\beta},
    \end{align*}
    verifying the probabilistic bound. To complete the proof, we use \Cref{lem: mu_1 and nu relationship}, from which it follows that
    \begin{equation*}
        \Vert \X_0-\X\Vert_\fro\leq\sqrt{2r}\Vert\X_0-\X\Vert \leq \sqrt{\frac{2\beta nr\log{n}}{3p}}\left(\max_{\ai}\left\vert\langle \X,\wa\right\rangle\vert \right) \leq \sqrt{\frac{\beta \nu^2 r^3 \log(n)}{24pn}}\Vert \X\Vert.
    \end{equation*}
    This concludes the proof.
\end{proof}

\section{Robustness Guarantees}\label{appendix: robustness}
In this section, we will prove \Cref{thm: Robustness guarantee for Mo}.
We first control the perturbation induced by additive noise.

\begin{lem}\label{lem: X perturbation bound}
    Let $\Ph = \P+\N$ where $\N\in\real^{n\times r}$ is a matrix with independent mean-zero entries, and $\Xh = \Ph\Ph^\top$. Let $\lambda_1\geq\cdots\geq\lambda_r>0$ be the non-zero eigenvalues of $\X$ with corresponding eigenvectors $\u^i$. Assume that $\Vert\N\Vert_\infty\leq C$ for $C>0$. Additionally, let $\bm{\Sigma} = \Ebb[{\n^i}{\n^i}^\top]$ be the covariance matrix of the columns of $\N$. Then 
    \[
    \Vert\X-\Xh\Vert \leq C\left(\frac{16}{3}n^{1/2}\beta\lambda_1^{1/2}\log{n}\right) + nrC^2
    \]
    with probability at least $1-2n^{1-\beta}$.
\end{lem}
\begin{proof}
    First, notice that
    \begin{equation*}
        \Ph\Ph^\top = \P\P^\top + \N\P^\top + \P\N^\top + \N\N^\top,
    \end{equation*}
    so 
    \begin{align*}
    \Vert\X-\Xh\Vert &= \Vert\N\P^\top + \P\N^\top + \N\N^\top\Vert\\
    &\leq 2\Vert\N\P^\top\Vert + \Vert\N\N^\top\Vert.
    \end{align*}

    We will first bound the term $\Vert\N\N^\top\Vert\leq\Vert\N\N^\top\Vert_\fro$. Notice that
    \begin{align*}
        \Vert\N\N^\top\Vert_\fro^2& = \sum_{i,j=1}^n\left(\sum_{k=1}^rN_{ik}N_{jk}\right)^2\\
        &\leq\sum_{i,j=1}^n\left( \sum_{k=1}^r C^2\right)^2\\
        & = n^2r^2C^4
    \end{align*}
    where the inequality follows from the bound on $\Vert\N\Vert_\infty$.
We now bound $\|\mathbf{N}\mathbf{P}^\top\|$. Since
$\mathbf{P}=\mathbf{U}\boldsymbol{\Lambda}^{1/2}$, where
$\mathbf{U}$ has orthonormal columns and
$\boldsymbol{\Lambda}=\operatorname{diag}(\lambda_1,\ldots,\lambda_r)$, we have $\|\mathbf{p}^{(k)}\|_2^2=\lambda_k$. We write $\mathbf{N}\mathbf{P}^\top$ as follows:
\[
    \mathbf{N}\mathbf{P}^\top
    =
    \sum_{k=1}^r
    \mathbf{n}^{(k)}\mathbf{p}^{(k)\top}
    =
    \sum_{k=1}^r \mathbf{Z}_k.
\]
The matrices $\mathbf{Z}_k$ are independent, mean-zero, and rank one. Moreover,
\[
    \|\mathbf{Z}_k\|
    =
    \|\mathbf{n}^{(k)}\|_2\,\|\mathbf{p}^{(k)}\|_2
    \leq
    C\sqrt{n\lambda_1}
    =: L.
\]
\noindent We next compute the Bernstein variance parameter. Let
$\bm{\Sigma}_k
    =
    \mathbb{E}\bigl[\mathbf{n}^{(k)}\mathbf{n}^{(k)\top}\bigr].$ As the entries of $\N$ are independent, $\bm{\Sigma} = \sum_k\bm{\Sigma}_k$ is diagonal, and as we have a bound on $\Vert \N\Vert_\infty$ it follows that 
    \[
    \Vert\bm{\Sigma}\Vert\leq C^2.
    \] and
$\tr(\boldsymbol{\Sigma}) =  nC^2$. Therefore,
\begin{align*}
    \left\|
    \sum_{k=1}^r
    \mathbb{E}\bigl[\mathbf{Z}_k\mathbf{Z}_k^\top\bigr]
    \right\|
    &=
    \left\|
    \sum_{k=1}^r
    \lambda_k \boldsymbol{\Sigma}_k
    \right\| 
    \leq
    C^2\sum_{k=1}^r \lambda_k
    \leq
    rC^2\lambda_1,
\end{align*}
while
\begin{align*}
    \left\|
    \sum_{k=1}^r
    \mathbb{E}\bigl[\mathbf{Z}_k^\top\mathbf{Z}_k\bigr]
    \right\|
    =
    \left\|
    \sum_{k=1}^r
    \mathbf{p}^{(k)}\mathbf{p}^{(k)\top}
    \mathbb{E}\|\mathbf{n}^{(k)}\|_2^2
    \right\| 
    \leq
    nC^2
    \left\|
    \sum_{k=1}^r
    \mathbf{p}^{(k)}\mathbf{p}^{(k)\top}
    \right\| 
    =
    nC^2 \lambda_1.
\end{align*}
Since $n\geq r$, we set $\sigma^2 = nC^2\lambda_1$. As 
    \begin{align*}
        \frac{\sigma^2}{c} & = \frac{nC^2 \lambda_1}{n^{1/2}C\lambda_{1}^{1/2}} =  n^{1/2}C\lambda_1^{1/2},
    \end{align*}
    it follows that for $t = \frac{8}{3}n^{1/2}C\beta\lambda_1^{1/2}\log{n}$ we have that
    \begin{align*}
    \mathbb{P}\left[\left\Vert \N\P^\top\right\Vert \geq t\right]&\leq 2n\exp\left(\frac{-3t}{8c}\right)\\
    & = 2n\exp\left(\frac{n^{1/2}C\beta\lambda_1^{1/2}\log{n}}{n^{1/2}C\lambda_1^{1/2}}\right)\\
    & = 2n\exp\left(-\beta{\log{n}}\right).
    \end{align*}
This finishes the proof.
\end{proof} 

Next we will prove the following Lemma showing that bounded noise on the points does not change the incoherence of a Gram matrix substantially.

\begin{lem}\label{lem: coherence parameter after perturbation}
For $\Ph = \P+\N$, where $\N$ is a mean-zero random matrix. Let $\X=\P\P^\top$ and $\Xh = \Ph\Ph^\top$, and let $\lambda_1\geq\cdots\geq\lambda_r>0$ be the eigenvalues of $\X$. If $\Vert\N\Vert_\infty \leq C$ for some $C>0$, where $\kappa$ is the condition number of $\X$, then
\[
\left\Vert\mathcal{P}_{\hat{U}}\wa\right\Vert_\fro\leq \sqrt{\frac{\nu r}{2n}} + \frac{32}{3}C\frac{n^{1/2}}{\lambda_r^{1/2}(\X)}\kappa^{1/2}\log{n} + \frac{2nr}{\lambda_r(\X)}C^2,
\]
with probability at least $1-2n^{1-\beta}$. Furthermore, if $C\leq\frac{3}{32\sqrt{2}}\frac{\nu^{1/2}}{n}\frac{\lambda_r(\X)^{1/2}}{\beta\kappa^{1/2}\log{n}}$, then 
\[
\left\Vert\mathcal{P}_{\hat{U}}\wa\right\Vert_\fro\leq \frac{9\nu r}{2n}.
\]
\end{lem}
\begin{proof}
   This result will follow from the classic Davis-Kahan $\sin\Theta$ Theorem, seen in \Cref{thm: Davis Kahan}. Let $\mathbb{U},~\mathbb{\hat{U}}$ be the subspace spanned by the columns of $\U,\hat{\U}$ respectively. First, as $\left\Vert\Pu - \mathcal{P}_{\hat{U}}\right\Vert_\fro = \Vert\sin\Theta(\mathbb{U},\hat{\mathbb{U}})\Vert_\fro$ from \cite{bhatia2013matrix}, we can see that
    \[
    \left\Vert\Pu - \mathcal{P}_{\hat{U}}\right\Vert_\fro = \Vert\sin\Theta(\mathbb{U},\mathbb{\hat{U}})\Vert_\fro\leq\frac{\Vert\X-\Xh\Vert_\fro}{\lambda_r}.
    \]
    Next, notice that
    \begin{align*}
        \left\Vert\mathcal{P}_{\hat{U}}\wa\right\Vert_\fro &\leq \left\Vert\left(\mathcal{P}_{\hat{U}} - \Pu\right)\right\Vert\Vert\wa\Vert_\fro+\Vert\Pu\wa\Vert_\fro\\
        &\leq \frac{2\Vert\X-\Xh\Vert}{\lambda_r(\X)} + \sqrt{\frac{\nu r}{2n}}.
        \end{align*}
        The Lemma statement therefore comes from the bound in \Cref{lem: X perturbation bound}. The second expression comes as follows:
        \begin{align*}
            \sqrt{\frac{\nu r}{2n}} + \frac{32}{3}C\frac{n^{1/2}}{\lambda_r^{1/2}(\X)}\kappa^{1/2}\beta\log{n} + \frac{2nr}{\lambda_r(\X)}C^2 & \leq 2\sqrt{\frac{\nu r}{2n}} + \frac{\nu r}{n}\frac{9}{2^{11}\beta^2 \kappa\log^2{n}}\\
            &\leq 2\sqrt{\frac{\nu r}{2n}} +   \frac{\nu r}{n}\\
            &=  2\sqrt{\frac{\nu r}{2n}} + \nu^{1/2}\frac{\nu^{1/2}r}{n}\\
            &\leq  2\sqrt{\frac{\nu r}{2n}} + \frac{n^{1/2}}{r^{1/2}}\frac{\nu^{1/2}r}{n}\\
            &=  3\sqrt{\frac{\nu r}{2n}},
        \end{align*}
        where the first inequality comes from substituting in the stated bound on $C$, the second inequality follows from the fact that $\beta^2\kappa\log^2{n}\geq 1$, and the third inequality follows from the fact that $\nu\leq\frac{n}{r}$, thus concluding the proof.
\end{proof}

\section{Incoherence Results}
In this section, we provide proofs for the statements in \Cref{sec:Incoherence_Discussion}.

\begin{lem}\label{lem: Incoherence equivalence}
    If $\left\Vert\Pu\wa\right\Vert_\fro\leq\sqrt{\frac{\nu r}{8n}}$, it follows that $\left\Vert\Pu\va\right\Vert_\fro\leq\sqrt{\frac{\nu r}{2n}}$. Similarly, if $\left\Vert\Pt\wa\right\Vert_\fro\leq\sqrt{\frac{\nu r}{8n}}$, it follows that $\left\Vert\Pt\va\right\Vert_\fro\leq\sqrt{\frac{\nu r}{2n}}$.
\end{lem}
\begin{proof}
    To see this result, notice that
        \begin{align*}
        \Vert\Pu\va\Vert_\fro &= \left\Vert \Pu\left(\sum_{\betab\in\mathbb{I}}H^{\alphab\betab}\wb\right)\right\Vert_\fro\\
        &\leq\sum_{\betab\in\mathbb{I}}\left\vert H^{\alphab\betab}\right\vert\Vert\Pu\wb\Vert_\fro\\
        &\leq \sqrt{\frac{\nu r}{8n}}\sum_{\betab\in\mathbb{I}}\left\vert H^{\alphab \betab}\right\vert,
    \end{align*}
    and as $\sum_{\ai}\vert H^{\alphab\betab}\vert \leq 2$ from \Cref{lem: H and H^-1 eigvals}, the claim follows. An identical proof shows the second result, with $\Pt$ in place of $\Pu$.
\end{proof}

\begin{lem}\label{lem: mu_1 and nu relationship}
    Let $\X \succeq \bm{0}$ be a rank-$r$, $\nu$-incoherent matrix satisfying \eqref{eq: Incoherence equations} with constant $\nu$. Then
    \[
    \left(\max_{\ai}\left\vert\langle \X,\wa\right\rangle\vert \right)\leq \frac{\nu r}{4n}\Vert\X\Vert.
    \]
\end{lem}
\begin{proof}
To see the above statement, notice that for $\X = \U\bm{\Lambda}\U^\top = \P\P^\top$,
    \begin{align*}
    \displaystyle\frac{\left(\max_{\ai}\left\vert\langle \X,\wa\right\rangle\vert \right)}{\Vert\X\Vert} &= \displaystyle\frac{1}{\Vert\X\Vert} \max_{i,j} \vert X_{ii}+X_{jj} - 2X_{ij}\vert\\
    &=\displaystyle\frac{1}{\Vert\X\Vert}\max_{i,j} \Vert\p_i-\p_j\Vert_2^2\\
    & = \displaystyle\frac{1}{\Vert\X\Vert} \max_{i,j} ~(\u_i-\u_j)^\top\bm{\Lambda}(\u_i-\u_j)\\
    &\leq \max_{ij}\,\, \Vert\u_i-\u_j\Vert_2^2\\ 
    &  = \max_{\ai}\,\,\frac{1}{2}\Vert\Pu\wa\Vert_\fro^2\\
    & \leq \frac{\nu r}{4n},
\end{align*}
where the third line follows from \eqref{eqn: MDS Equation}, the fourth line follows from the fact that $\X\succeq\bm{0}$, the fifth line follows from \eqref{eqn: Puwa and u_i-u_j}, and the final line follows from \eqref{eq: Incoherence equations}, thus concluding the proof.
\end{proof}

\begin{lem}\label{lem: anisotropic incoherence}
    Let $\mu$ be an a.s. bounded, mean-zero, sub-Gaussian distribution with positive definite covariance matrix $\bm{\Sigma}\in\real^{r\times r}$. Let $n$ points $\{\p_i\}_{i=1}^n\sim \mu$ be sampled i.i.d., and let $\P = [\p_1 \dots \p_n]^\top\in\real^{n\times d}$ be the corresponding point matrix with Gram matrix $\X$, which has condition number $\kappa$. Let $\Vert\p_i\Vert_{\psi_2}\leq K$ for some $K>0$. Then with probability at least $1-Cn^{-2}$ for some absolute constant $C>0$, the incoherence parameter of $\X$ is bounded by
    \[
    \nu \leq \cO\left(\kappa+\frac{K^2\log{n}}{\sqrt{r}}\right).
    \]
\end{lem}
\begin{proof}
    This proof is much the same as the proof in \Cref{sec:Incoherence_Discussion}. 
    First, we remark that \begin{align*}
    \Ebb\left[\left(\p_i-\p_j\right)^\top\left(\p_i-\p_j\right)\right] & =  \Ebb\left[\left\Vert\p_i\right\Vert^2_2\right]-\Ebb\left[{\p_j}^\top\p_i\right]-\Ebb\left[{\p_j}^\top\p_i\right]+\Ebb\left[{\p_j}^\top\p_j\right]\\
    & = \Ebb\left[\left\Vert\p_i\right\Vert^2_2\right] + \Ebb\left[\left\Vert\p_j\right\Vert^2_2\right]- 2\Ebb\left[{\p_i}\right]^\top\Ebb\left[\p_j\right] \\
    & = \Ebb\left[\left\Vert\p_i\right\Vert^2_2\right] + \Ebb\left[\left\Vert\p_j\right\Vert^2_2\right]\\
    & = 2\Ebb\left[\left\Vert\p_i\right\Vert^2_2\right]\\
    & = 2\Ebb\left[\tr\left({\p_i}{\p_i}^\top\right)\right]\\
    & = 2~\tr\left(\Ebb\left[\p_i{\p_i}^\top\right]\right)\\
    & = 2~\tr(\bm{\Sigma}) \leq 2r\lambda_{1}(\bm{\Sigma}),
\end{align*}
where the second and fourth lines follow from the independence of $\p_i$ and $\p_j$, the third line follows from the fact that $\Ebb[\mu]=0$, and the seventh line follows from the fact that $\bm{\Sigma}$ has $r$ non-zero eigenvalues.

Next, following the argument of \Cref{lem: hwi proof for distance calc} but replacing $2r$ with $\mathbb{E}\left[\left(\p_i-\p_j\right)^\top\left(\p_i-\p_j\right)\right]$, we have that, with probability at least $1-Cn^{-2}$
\[
\left\Vert\p_i-\p_j\right\Vert_2^2\leq 2r\lambda_{1}(\bm{\Sigma}) +  4K^2\sqrt{r}\log{n}.
\]

Now, leveraging that from \Cref{lem: vershynin covariance estimate} that $\lambda_{r}(\X) \approx n\lambda_r(\bm{\Sigma})$, we have that for some $c>0$ with high probability that
\begin{align*}
    \nu&\leq \frac{n}{2r}\frac{2r\lambda_1(\bm{\Sigma})+4K^2\sqrt{r}\log{n}}{\lambda_r(\X)}\\
    & = \frac{n}{2r}\frac{2r\lambda_1(\bm{\Sigma})+4K^2\sqrt{r}\log{n}}{cn\lambda_r(\bm{\Sigma})}\\
    &\leq \frac{\kappa}{c}+\frac{2K^2 \log{n}}{c\sqrt{r}}\\
    &=\cO\left(\kappa+\frac{K^2\log{n}}{\sqrt{r}}\right).
\end{align*}
This concludes the proof.
\end{proof}

\section{Further Background}\label{appendix: background}
\subsection{Dual Bases}
In a finite dimensional vector space of matrices $\mathbb{V}$, where $\mathrm{dim}(\mathbb{V}) = n$, a basis is a linearly independent set of matrices $B = \{\X_i\}_{i=1}^{n}$ that spans $\mathbb{V}$. Any basis for a finite-dimensional vector space admits a dual, or bi-orthogonal, basis denoted $B^{*} = \{\Y_i\}_{i=1}^{n}$ that also spans $\mathbb{V}$, and admits a bi-orthogonality relationship
\[
\langle \X_i,\Y_j \ra = \delta_{ij}.
\]
Additionally, $B$ uniquely determines $B^*$. The bi-orthogonality relationship allows for the decomposition of any matrix $\Z\in \mathbb{V}$ as follows:
\[
\Z = \sum_{i=1}^{n} \langle \Z,\Y_i\rangle \X_i = \sum_{i=1}^{n} \langle \Z,\X_i\rangle \Y_i. 
\]
We define the Gram, or correlation matrix, $\H\in\Rnn$, for $B$ as $H_{ij} = \la\X_i,\X_j\ra$, and let $H^{ij} = (\H^{-1})_{ij}$. It is straightforward to show that $\Y_i = \sum_{j=1}^{n} H^{ij}\X_j$ generates $B^*$, and similarly that $\X_i = \sum_{j=1}^{n}H_{ij}\Y_j$ \cite{bhatia2013matrix}.

\subsection{Riemannian Optimization}\label{subsec: intro to RO}

The primary setting for this work is the Riemannian manifold of fixed-rank matrices. Throughout this work, we will only be considering square $n \times n$ matrices for simplicity and relevance to the problem of interest in this paper. For a fixed positive integer $r\leq n$, we denote the set $\mfr = \{\X\in\Rnn~\vert ~\mathrm{rank}(\X) = r\}$. Although not obvious at first glance, it is well-known that $\mfr$ is a smooth Riemannian manifold\cite{vandereycken2012lowrank,Boumal2023}. To make this a Riemannian manifold, we equip it with the standard trace inner product as a metric, or $\langle\A,\B\rangle = \tr(\A^\top\B)$, restricted to the tangent bundle $T\mfr$, which is the disjoint union of tangent spaces\cite{Boumal2023}.

Additionally, the tangent space at a point $\X\in\mfr$ is known and can be characterized \cite{vandereycken2012lowrank,Boumal2023,wei2020guarantees}. For notational simplicity, and of relevance in the context of optimization, assume that $\X$ is the ground truth solution to an objective function. We additionally assume that $\X = \X^\top$, as all the matrices we consider are symmetric. The following ideas can be re-stated for rectangular matrices using a singular value decomposition, but these are not the subject of this paper. As such, we denote the tangent space at $\X$ as $\T$, and for a sequence of iterates $\{\X_l\}_{l\geq0}$, we refer to their respective tangent spaces as $\T_l$. To characterize $\T$, let $\X = \U\D\U^\top$ be the thin spectral decomposition of $\X$. The tangent space $\T$ can be computed as follows:
\[
\T = \{\U \Z^\top + \Z\U^\top ~\vert~ \Z\in\real^{n\times r}\}.
\]
The tangent space can be described as the set of all possible rank-up-to-$2r$ perturbations, represented as the sum of a perturbation in the column and row space, and is computed by looking at first-order perturbations of the spectral decomposition of $\X$\cite{vandereycken2012lowrank}. Additionally, we can compute the orthogonal projection of an arbitrary $\Y\in\Rnn$ onto the tangent space at a point $T_{\X} \mfr$ as follows \cite{vandereycken2012lowrank,Boumal2023,wei2020guarantees}:
\begin{equation*}
    \Pt \Y  = \Pu\Y + \Y\Pu - \Pu\Y\Pu,
\end{equation*}
where $\Pu = \U\U^\top$ is the orthogonal projection onto the subspace spanned by the $r$ columns of $\U$.

Optimization over $\mfr$ has been investigated in detail for quite some time, and retraction-based methods are of particular interest to this work \cite{Absil2008,shalit2012,vandereycken2012lowrank,wei2016guarantees,cai2019accaltproj,wei2020guarantees,cai2021asap,hamm2022RieCUR}. First-order retraction-based methodologies rely on the general principle of taking a descent step in the tangent space, followed by a retraction onto the manifold. In the case of first-order optimization on $\mfr$, the retraction map $\mathcal{H}_r$ is given by the hard thresholding operator, which is a thin spectral decomposition that takes $\Y=\sum_{i=1}^n\lambda_i\u_i\u_i^\top\mapsto\sum_{i=1}^r \lambda_i\u_i\u_i^\top$, where $\vert\lambda_1\vert\geq \cdots\geq \vert\lambda_n\vert$ are the ordered eigenvalues of $\Y$ and $\u_i$ are the corresponding eigenvectors of $\Y$.

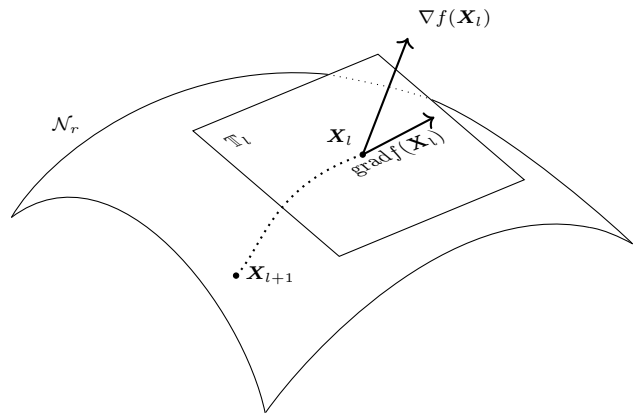
\begin{figure}
    \centering
\begin{tikzpicture}[x=1pt, y=1pt, yscale=-1]

  height of 300

  \draw (112.04,66.62) node[below right, font=\scriptsize]{$\mfr$};
  \draw (100,108) .. controls (140,78) and (188,143) 
                  .. (196,182) ;
  \draw (196.2,181.8) .. controls (211.8,159.6) and (265.71,101.03) 
                      .. (318.8,111.8) ;
  \draw (318.8,111.8) .. controls (323.67,112.42) and (330.33,115.42) 
                      .. (335,118) ;
  \draw[dotted] (219.48,53.8) .. controls (230.71,56.06) and (244.01,58.07) 
                              .. (259,63.39) ;
  \draw (259,63.39) .. controls (288.81,75.27) and (332.14,114.93) 
                    .. (335,118) ;
  \draw (100,108) .. controls (113.56,84.56) and (160.04,48.06) 
                  .. (219.48,53.8) ;

  \filldraw (232.9,84.2) circle (1pt)
                         node[above left, font=\scriptsize]{$\X_l$};

  \draw[thick,->] (232.9,84.2) --  (250,40)
                    node[above right,font=\scriptsize]{$\nabla f(\X_l)$};

  \draw[dotted,thick] (232.9,84.2) .. controls (198,96) and (195,120) 
                .. (185,130);

  \draw[thick, ->] (232.9,84.2) -- (260,70)
                        node[below left,font=\scriptsize,rotate=-337.86,
                           xslant=-0.54]{$\textrm{grad}f(\X_l)$};

    \filldraw (185,130)circle (1pt)
                         node[right, font=\scriptsize]{$\X_{l+1}$};                       

  \draw (238.68,46.2) -- (294.01,93.72) 
                      -- (224.01,122.72) 
                      -- (168.68,75.2) 
                      -- cycle ;
  \draw (173.51,76.9) node[below right, 
                           font=\scriptsize, 
                           rotate=-337.86,
                           xslant=-0.54]{$\T_l$};
\end{tikzpicture}

    \caption{A diagram of a simple first-order retraction method on $\mfr$. Again, $\nabla f(\X_l)$ is the Euclidean gradient of $f$ at $\X_l$, $\mathrm{grad}\, f(\X_l)$ is the Riemannian gradient at $\X_l$, and $\X_{l+1} = \mathcal{H}_r(\X_l - \alpha_l\mathrm{grad}\,f(\X_l))$, as in \eqref{eqn: first order retraction}.}
    \label{fig:First order retraction figure}
\end{figure}

In order to construct a first-order method on $\mfr$, we need to define the notion of a Riemannian gradient. This object can be constructed in a greater degree of generality than our approach, but for simplicity, we will assume that a function $f:\mfr\to \real$ can be smoothly extended to all of $\Rnn$. That is to say, if we consider $f:\Rnn\to\real$, the Riemannian gradient of $f\big\vert_{\mfr}$, denoted $\mathrm{grad} ~f$, for $\X_l\in\mfr$ is given by:
\[
\mathrm{grad}\,f(\X_l) = \Ptl \nabla f(\X_l),
\]
where $\nabla f$ is the Euclidean gradient of $f$. Using this approach, we can now define a Riemannian gradient descent iterate sequence using our retraction map, Riemannian gradient, and some step size sequence $\{\alpha_l\}_{l\geq 0}$ as follows:
\begin{equation}\label{eqn: first order retraction}
\X_{l+1} = \mathcal{H}_r(\X_l-\alpha_l \Ptl \nabla f(\X_l)).
\end{equation}
Intuitively, this algorithm seeks to look at changes in the objective function that lie, locally, along the manifold, followed by a retraction to stay on the desired manifold. An illustration can be seen in \Cref{fig:First order retraction figure}.

This is a simple first pass to first-order optimization on Riemannian manifolds, and is not meant to be exhaustive. Interested readers should consult \cite{Absil2008,Boumal2023} for further details on first-order methods on matrix (and other Riemannian) manifolds, along with convergence analysis for these algorithms.

\end{document}